\documentclass[12pt]{amsart}
\usepackage{amsmath}

\newtheorem{thm}{Theorem}[section]

\newtheorem{lem}[thm]{Lemma}
\newtheorem{prop}[thm]{Proposition}
\newtheorem{cor}[thm]{Corollary}

\theoremstyle{definition}
\newtheorem{defn}{Definition}[section]

\theoremstyle{remark}
\newtheorem{remark}{Remark}[section] 

\newtheorem{rem}{Remark}[section]

\theoremstyle{plain}

\newcommand{\half}{\frac{1}{2}}
\newcommand{\Z}{{\mathbf Z}}
\newcommand{\ucalr}{{\mathcal U}_{ir}}

\newcommand{\ccal}{{\mathcal C}}
\newcommand{\fcal}{{\mathcal F}}
\newcommand{\ocal}{{\mathcal O}}
\newcommand{\dcal}{{\mathcal D}}
\newcommand{\pcal}{{\mathcal P}}
\newcommand{\zcal}{{\mathcal Z}}
\newcommand{\llcal}{{\mathcal L}}
\newcommand{\lcal}{\zcal}
\newcommand{\rcal}{{\mathcal R}}

\newcommand{\X}{{\mathbf X_{\Gamma}}}
\newcommand{\HH}{\mathbf H}
\newcommand{\D}{\mathbf D}
\newcommand{\N}{\mathbf N}
\newcommand{\R}{{\mathbf R}}
\newcommand{\C}{{\mathbf C}}
\newcommand{\Lap}{\triangle}
\newcommand{\To}{\longrightarrow}
\newcommand{\cal}{\mathcal}
\newcommand{\eps}{\varepsilon}
\def\1{{{\mbox{${\mathrm{1\negthinspace\negthinspace I}}$}}}}

\numberwithin{equation}{section}

\begin{document}

\title
{Patterson-Sullivan distributions and quantum ergodicity}

\author{Nalini Anantharaman and  Steve Zelditch}
\address{\'Ecole Normale Superieure de Lyon, 46 all\'ee d'Italie, 69364 Lyon Cedex 07}
\email{nanantha@umpa.ens-lyon.fr}
\address{Johns Hopkins University, Baltimore, MD
21218, USA}\email{szelditch@jhu.edu}

\date{January 31, 2006}

\thanks{\\
 Research partially supported
  by  NSF grant \#DMS-0302518 and NSF Focussed
Research Grant  \# FRG 0354386}

\maketitle

\begin{abstract} This article gives relations between two types of
phase space distributions associated to eigenfunctions
$\phi_{ir_j}$ of the Laplacian on a compact hyperbolic surface
$X_{\Gamma}$:

\begin{itemize}

\item  Wigner distributions $\int_{S^*\X} a \;dW_{ir_j}=\langle
Op(a)\phi_{ir_j}, \phi_{ir_j}\rangle_{L^2(\X)}$, which arise in
quantum chaos. They are invariant under the wave group.

\item  Patterson-Sullivan distributions $PS_{ir_j}$, which are
the residues of the dynamical zeta-functions $\lcal(s; a): = \sum_\gamma
\frac{e^{-sL_\gamma}}{1-e^{-L_\gamma}}
\int_{\gamma_0} a$ (where the sum runs over closed geodesics) at the poles $s = \frac{1}{2} + ir_j$. They are invariant
under the geodesic flow.

\end{itemize}

We prove that these distributions  (when suitably normalized) are
asymptotically equal as $r_j \to \infty$. We also give exact
relations between them. This correspondence gives a new relation
between classical and quantum dynamics on a hyperbolic surface,
and consequently a  formulation of quantum ergodicity in terms of
classical ergodic theory.

\end{abstract}

\tableofcontents

\section{Introduction, statement of results.}

The purpose of this article is to relate two kinds of phase space
distributions which are naturally attached to the eigenfunctions
$\phi_{ir_j}$ of the Laplacian $\Lap$ on a compact hyperbolic
surface $\X$. The first kind are  the   {\it Wigner distributions}
$W_{ir_j}\in \dcal'(S^* \X)) $ (\ref{WIGDEF}) of quantum
mechanics. The second kind are what we call  {\it
Patterson-Sullivan distributions} $PS_{ir_j}  \in \dcal'(S^* \X))$
(\ref{PATSULDEF}). We prove that the Patterson-Sullivan
distributions have a purely classical mechanical definition: they
are the residues of  classical dynamical zeta functions at poles
in the `critical strip'.  We also prove that there exists an
`intertwining operator' $L_r$ (\ref{FIRSTLR}) which transforms
$PS_{ir_j} \to W_{ir_j}$ and which induces an asymptotic equality
$W_{ir_j} \sim  2^{-(1+2ir_j)}(-4i\pi r_j)^{-1/2} PS_{ir_j}$
between them. These results indicate a surprisingly close relation
between quantum and classical dynamics on a compact hyperbolic
surface. A similar relation on finite volume hyperbolic manifolds
of all dimensions is to be expected, but certainly these
relations are special properties of locally symmetric manifolds.

To state our results, we introduce some notation.  We write  $G =
 PSU(1,1), K = PSO(2)$ and identify the quotient $G/K$ with the hyperbolic disc $\D$.
 We let $\Gamma \subset G$ denote a co-compact discrete group and
 let
  $\X = \Gamma \backslash \D$ denote the associated hyperbolic
 surface. By ``phase space'' we mean the unit cotangent bundle $S^*
 \X$, which may be identified with the tangent bundle $S
 \X$ and also with the quotient
$\Gamma \backslash G$.
 We
denote by  $\lambda_0=0<\lambda_1\leq \lambda_2...$ the spectrum
of the laplacian on $\X$; with the usual parametrization
$\lambda_j=s_j(1-s_j)=\frac14+r_j^2$ ($s_j=\frac12+ir_j$), we denote by $\{\phi_{ir_j}\}_{j = 0, 1, 2,
\dots}$ an orthonormal basis of real-valued eigenfunctions: $\Lap
\phi_{ir_j}=-\lambda_j \phi_{ir_j}$.

The   Wigner distributions (microlocal lifts, microlocal defect
measures...)
 $W_{ir_j}\in \dcal'(S^* \X)) $ are defined by
\begin{equation} \label{WIGDEF}\int_{S^*\X} a(g)W_{ir_j}(dg)=\langle Op(a)\phi_{ir_j},
\phi_{ir_j}\rangle_{L^2(\X)},\;\;\; a \in C^{\infty}(S^* \X)
\end{equation}  where $Op(a)$ is a special quantization of $a$, defined
using hyperbolic Fourier analysis (Definition \ref{WIGNER}). The
Wigner distribution $W_{ir_j}$ depends quadratically on
$\phi_{ir_j}$ and has the quantum invariance property
\begin{equation} \label{QINV} \langle U_t^* Op(a) U_t \phi_{ir_j},
\phi_{ir_j} \rangle = \langle Op(a) \phi_{ir_j}, \phi_{ir_j}
\rangle, \;\;\; (U_t = \exp (i t \sqrt{\Delta})) ;\end{equation}
hence by Egorov's theorem $W_{ir_j}$ is asymptotically invariant
under the action of the geodesic flow $g^t$ on $S^*\X$, in the
large energy limit $r_j\To +\infty$. The Wigner distribution
 $W_{ir_j} $ is one of the principal objects in quantum chaos:
it determines the oscillation and concentration  of the
eigenfunction $\phi_{ir_j}$ in the classical phase space $S^*\X$
(see \S \ref{BACKGROUND}). One of the main problems in quantum
chaos is the
 {\em quantum unique ergodicity problem} of
determining which geodesic flow invariant probability measures
arise as  weak* limits of the Wigner distributions (cf. \cite{Lin,
RS, W, Sh, SV, Z2} for a few articles on hyperbolic quotients).

The  Patterson-Sullivan distributions $\{PS_{ir_j}\}$ associated
to the eigenfunctions $\{\phi_{ir_j}\}$ (cf. Definition
\ref{PSDEF})
are defined by the expression
\begin{equation} \label{PATSULDEF} PS_{ir_j}(dg)=PS_{ir_j}(db', db, dt) := \frac{T_{ir_j}(db)
T_{ir_j}(db')}{|b - b'|^{1 +  2i r_j}} \otimes |dt|.
\end{equation} In this definition, $T_{ir_j}$ is the boundary values
of $\phi_{ir_j}$ in the sense of Helgason (cf. Theorem \ref{Helga}
or \cite{He, H}.) The parameters $(b',b)$ ($b\not=b'$) vary in
$B\times B$, where $B =
\partial \D$ is the boundary of the hyperbolic disc, and $t$
varies in $\R$; $(b', b)$ parametrize the space of oriented
geodesics, $t$ is the time parameter along geodesics, and the
three parameters $(b', b, t)$ are used to parametrize the unit
tangent bundle $S\D$.  These Patterson-Sullivan distributions
$PS_{ir_j}$  are invariant under the geodesic flow $g^t$ on $S
\X$, i.e.
\begin{equation} \label{CLINV} (g^t)_* PS_{ir_j} = PS_{ir_j}. \end{equation}
 The
distributions $PS_{ir_j}$ are also $\Gamma$-invariant (cf.
Proposition \ref{INVARIANT}). Hence they define geodesic-flow
invariant distributions on $S \X$. Phase space
distributions of this kind were associated to ground state
eigenfunctions of certain infinite area hyperbolic surfaces by  by
S. J. Patterson \cite{Pat0, Pat1}, and were studied further by D.
Sullivan \cite{Su1, Su2} (see also \cite{N}). Ground state
Patterson-Sullivan distributions are positive measures, but our
analogues for higher eigenfunctions on compact (or finite area)
hyperbolic surfaces are not measures. To our knowledge, they have
not been studied for higher eigenfunctions before.

Both families $(W_{ir_j})$ and $(PS_{ir_j})$ are
$\Gamma$-invariant bilinear forms in the eigenfunctions
$\phi_{ir_j}$ with values in distributions on $S\X$. But they
possess  different invariance properties: the former are invariant
under the quantum dynamics (the wave group) while the latter are
invariant by the classical evolution (the geodesic flow). The
motivating problem in this article is to determine  how they are
related.

 Our first result
is that, after suitable normalization of $PS_{ir_j}$, the Wigner
and the Patterson-Sullivan distributions are asymptotically the
same in the semi-classical limit $r_j\To +\infty$. The exact relation
between  involves
 the operator $L_r: C_0^{\infty}(G) \to
 C^{\infty}(G)$ defined by
\begin{equation} \label{FIRSTLR} L_{r}a(g)=\int_{\R} (1+u^2)^{-(\frac12+ir)}
a(gn_u)du\end{equation} which, we will see, mediates between the
classical and quantum pictures. Here, $n_u = \left(
\begin{array}{ll} 1 & u \\ & \\
0 & 1 \end{array} \right) $ acts on the right as the horocycle flow. We further introduce a cutoff
function $\chi \in C_0^{\infty}(\D)$ which is a smooth replacement
for the characteristic function of a fundamental domain for
$\Gamma$ (called a `smooth fundamental domain cutoff', see
Definition \ref{SFD}).

\begin{thm} \label{mainintro} For any $a \in C^{\infty}(\Gamma
\backslash G)$ we have the exact formula
$$\langle Op(a)\phi_{ir_j}, \phi_{ir_j}\rangle=2^{-(1+2ir_j)}\int_{S \D} (L_{r_j}\chi  a)
dPS_{ir_j},$$ and the asymptotic formula
$$\int_{S\X} a(g)W_{ir_j}(dg)= 2^{-(1+2ir_j)}(-4i\pi r_j)^{-1/2}\int_{S\X} a(g) PS_{ir_j}(dg)+ O(r_j^{-1}).$$
\end{thm}

It follows  that the Wigner distributions are equivalent to the
Patterson-Sullivan distributions in  the study of quantum
ergodicity. The  operators $L_{r}$ in a sense intertwine classical
and quantum dynamics. We note that, although the Wigner
distributions were defined by using the special hyperbolic
pseudodifferential calculus $Op$, any other choice of $Op$ will
produce asymptotically equivalent Wigner distributions and hence
Theorem \ref{mainintro} is stable under change of quantization.

When $a$ is an automorphic eigenfunction, i.e. a joint
eigenfunction of the Casimir operator $\Omega$ and the generator
$W$ of $K$, we can evaluate the first expression in Theorem
\ref{mainintro} to obtain a very concrete relation:

\begin{thm} \label{maincont} If $\sigma$ is an eigenfunction of Casimir parameter
$\tau$ and weight $m$ in the continuous series, we have:

\noindent{\bf (i)}

$$  \begin{array}{lll}  \langle Op(\sigma)\phi_{ir_j}, \phi_{ir_j}\rangle&=&
 \mu_{m, \tau}^c(\frac{1}{2} + i r_j)  \left(\int_{\Gamma \backslash G}
 \sigma \; dPS_{ir_j}\right) + \mu_{m, \tau}^{c odd} (\frac{1}{2} + i r_j) \left(\int_{\Gamma \backslash
 G}
 X_+\sigma \; dPS_{ir_j}\right) ,
\end{array}$$
where $\mu_{m, \tau}^c$ and $\mu_{m, \tau}^{c odd}$ are defined in
\eqref{BIGDISPLAY}; $X_+$ denotes the vector field generating the
horocycle flow.

\noindent{\bf (ii)}   If  $\sigma=\psi_m$ is a lowest weight
vector in the holomorphic discrete series, we have
$$ \langle Op(\psi_m)\phi_{ir_j}, \phi_{ir_j}\rangle= \mu_m^d
(\frac{1}{2} + i r_j) \left(\int_{\Gamma \backslash G}  \sigma
dPS_{ir_j}\; \right),$$ where $\mu_m^d$ is defined in
(\ref{WTZERO}).

\end{thm}

These exact formulae are based on the identity (cf.   Proposition
\ref{CONTINUITY}),
\begin{equation} \label{CONTINUITYint} \int_{S \D} (L_{r_j}\chi  \sigma) dPS_{ir_j}  = \int_\R (1 +
u^2)^{-(\frac{1}{2}+ir_j)} I_{PS_{ir_j}}( \sigma)(u)  du,
\end{equation} where  $I_{PS_{ir_j}}: C^{\infty}(\Gamma \backslash
G) \to C(\R)$ is the operator defined by
\begin{equation} \label{IMUPSint} I_{PS_{ir_j}}(\sigma)(u):=\int_{\Gamma
\backslash G}
 \sigma(g n_u) PS_{ir_j}(dg). \end{equation}
When $\sigma$ is a joint eigenfunction of the Casimir operator
$\Omega$ and of the generator $W$ of the maximal compact subgroup
$K$, the function $ I_{PS_{ir}}(\sigma)(u)$ is a special function
of hypergeometric type depending on $r$ and the eigenvalue
parameters of $\sigma$ (cf.  \S \ref{BACKGROUND} for  a review of
the representation theory of $L^2(\Gamma \backslash G)$). The
integral on the right side of (\ref{CONTINUITYint}) can then be
evaluated to give the  explicit formulae of Theorem
\ref{maincont}.

Theorems \ref{mainintro} and \ref{maincont}  admit generalizations
to off-diagonal matrix elements, which relate off-diagonal Wigner
and Patterson-Sullivan distributions  by operators $L_{r, r'}$.
These generalizations  will be explored in our subsequent article
\cite{AZ}. It is obvious that there cannot exist any standard kind
of intertwining between the quantum dynamics (which has a discrete
spectrum) and the classical dynamics (which has a continuous
spectrum), but the operators $L_{r,r'}$ establish a kind of
intertwining on the level of invariant distributions.

Our next result gives a purely classical dynamical
interpretation of the Patterson-Sullivan distributions in terms of closed geodesics. Given $a
\in C^{\infty}(SX_{\Gamma})$, we define two closely related
dynamical zeta-functions
\begin{equation} \label{LFUNDEFa} \left\{ \begin{array}{ll} (i) & \lcal(s; a): = \sum_\gamma
\frac{e^{-s L_{\gamma}}}{1-e^{-L_{\gamma}}} \left(\int_{\gamma_0}
a \right), \;\;\; (\Re e\;s > 1), \\ &\\(ii) &  \lcal_2(a,s)=
\sum_\gamma  \frac{e^{ -(s-1)L_{\gamma}}}{|\sinh(L_{\gamma}/2)
|^2} \left(\int_{\gamma_0}a \right)
\end{array} \right. \end{equation}  where the sum runs over all
closed orbits, and $\gamma_0$ is the primitive closed orbit traced
out by $\gamma$. The sum converges absolutely for $\Re e\;s
> 1$.

\begin{thm} \label{main2}Let $a$ be a real analytic function on the unit tangent
bundle. Then $\lcal(s; a)$ and $\lcal_2(s; a)$ admit a meromorphic extension to $\C$.
The poles in the critical strip $0 <\Re e\;s < 1$, appear at
$s=1/2+ir$, where as above $1/4+r^2$ is an  eigenvalue of $\Lap$.

The residue is
$$\mu_0(\frac{1}{2} + i r) \sum_{j: r_j^2 = r^2} \int a dPS_{ir_j},$$
where $\{dPS_{ir_j}\}$ are the Patterson-Sullivan distributions
associated to an  orthonormal eigenbasis $\{\phi_{ir_j}\}$, and $\mu_0$
is defined in \eqref{BIGDISPLAY}.

\end{thm}
The normalisation factor $\mu_0(\frac{1}{2} + i r_j)$ is
 such that $\mu_0(\frac{1}{2} + i r_j)\int \1
dPS_{ir_j}=1$, for a Patterson-Sullivan distribution associated to
an eigenfunction normalized in $L^2$. Note that it does not depend
on the group $\Gamma$.

We give two quite different proofs of this theorem, each one
providing  a different perspective on the correspondence between
$W_{ir_j}$ and $PS_{ir_j}$. The first proof uses the thermodynamic
formalism. It is known that this formalism can be used  to locate
 the zeros of Selberg's zeta function \cite{Pol}. In \S
\ref{THERMO}, we apply similar ideas to prove that $\lcal_2$
admits a meromorphic continuation to $\Re e\;s > 0$, and we
describe its poles and residues in terms of ``Ruelle resonances".
In particular, Patterson-Sullivan distributions arise as the
residues. We use the methods developed by Rugh in \cite{Rugh92}
for real-analytic situations. The techniques are based on the
Anosov property of the geodesic flow, and apply in variable
curvature. However,  the relation between Wigner and
Patterson-Sullivan distributions is special to constant curvature.

The  second proof,  given in
 \S
\ref{SELBERG}, uses representation theory and the generalized
Selberg trace formula of \cite{Z} to prove the meromorphic
continuation of $\lcal$ to all of $\C$ and to determine  the poles
and residues. Like Theorems \ref{mainintro}- \ref{maincont}, the
trace formula establishes an exact  relation  between the Wigner
distributions (which appear on the `spectral side' of the trace
formula) and the geodesic periods $\int_{\gamma} a$ on the `sum
over $\Gamma$' side.  No such formula can be expected in variable
curvature, and the methods are  specific to hyperbolic surfaces.

In conclusion, the results of this paper develop to  a new level
the strangely  close relation between  classical and quantum
dynamics on hyperbolic surfaces. On the level of eigenvalues and
lengths of closed geodesics, this  close relation  is evident from
the Selberg trace formula (cf. \S \ref{SELBERG}). As is
well-known, the Selberg trace formula on a compact hyperbolic
manifold is a special case of the  general wave trace formula on a
compact Riemannian manifold  where the leading order approximation
is exact. The exactness of this stationary phase formula is
somewhat analogous to the exact stationary phase formula of
Duistermaat-Heckman for certain oscillatory integrals, but to our
knowledge no rigorous link between these exact formulae is known.
An alternative explanation of the close relation between classical
and quantum dynamics was suggested by V. Guillemin in \cite{G},
who made a formal application of the Lefschetz formula to the
action of the geodesic flow on a non-elliptic complex. The trace
on chains gave the logarithmic derivative of the (Ruelle) zeta
function, while the trace on homology gave the spectral side of
the Selberg trace formula. For later developments in this
direction (by C. Denninger, A. Deitmar, U. Bunke, M. Olbrich and
others) we refer to \cite{J}.

This paper develops the close relation on the level of
eigenfunctions and invariant distributions rather than just
eigenvalues and lengths of closed geodesics. In future work
\cite{AZ}, the intertwining of classical and quantum mechanics
will be developed further so that it applies not just to traces
(diagonal matrix elements) but to the full dynamics. It is hoped
that it will have applications related to the question of quantum
unique ergodicity. It would also be
 interesting to relate our constructions to generalizations of the
 non-elliptic Lefschetz formulae and to other representation theoretic ones
 in  \cite{SV, W}.

\medskip \noindent{\bf Acknowledgements} This work was begun while
the first author was visiting Johns Hopkins University as part of
the NSF focussed research grant  \# FRG 0354386. Much of it was
written at the  Time at Work program of the Institut Henri
Poincar\'e in June, 2005.

\section{\label{BACKGROUND} Background}

Hyperbolic surfaces are uniformized by the hyperbolic plane $\HH$
or disc $\D$. In  the disc model $\D=\{z\in\C, |z|<1\}$, the
hyperbolic  metric has the form
$$ds^2=\frac{4|dz|^2}{(1-|z|^2)^2}.$$
The group of orientation-preserving isometries can be identified
with $PSU(1,1)$ acting by Moebius transformations; the stabilizer
of $0$ is $K\simeq SO(2)$ and thus we will often identify $\D$
with $SU(1,1)/K$. Computations are sometimes simpler in the $\HH$
model, where the isometry group is $PSL(2, \R)$. We therefore use
the general notation $G$ for the isometry group, and $G/K$ for the
hyperbolic plane,  leaving it to the reader and the context to
decide whether $G = PSU(1,1)$ or $G = PSL(2, \R)$.

In hyperbolic polar coordinates centered at the origin $0$, the
Laplacian is the operator
$$\Lap=\frac{\partial^2}{\partial r^2} +\coth r\frac{\partial}{\partial r}+\frac{1}{\sinh^2 r}\frac{\partial^2}{\partial \theta^2} .$$
 The distance on $\D$ induced by the Riemannian
metric will be denoted $d_\D$. We denote the volume form by
$dVol(z)$.

Let $\Gamma \subset G$ be a co-compact  discrete subgroup, and let
us consider the automorphic eigenvalue problem on $G/K$:
\begin{equation} \left\{ \begin{array}{l} \Lap \phi = - \lambda\phi,\\ \\
\phi(\gamma z) = \phi(z) \mbox{ for all } \gamma\in\Gamma \mbox{
and for all } z. \end{array} \right. \end{equation}

In other words, we study the eigenfunctions of the Laplacian on
the compact surface $\X=\Gamma\backslash \; G \;/K$. Following
standard  notation (e.g. \cite{V, O}), the eigenvalue can be
written in the form $\lambda=\lambda_r=\frac14 +r^2$ and also
$\lambda=\lambda_s=s(1- s)$ where $s=\frac12+ir$.
\medskip

\noindent{\bf Notational remarks}

\noindent (i) We denote by $ \{\lambda_j = \frac{1}{4} + r_j^2\}$
the set of eigenvalues repeated according to multiplicity, and (in a
somewhat abusive manner) we
denote a corresponding  orthonormal basis of  eigenfunctions by
$\{\phi_{ir_j}\}$.

\noindent(ii)  We follow the notational conventions used in
\cite{N} and \cite{O}, which differ from those used in \cite{H},
where the metric is multiplied by a factor $4$. We caution that
\cite{L, Z} use the latter conventions, and there the parameter
$s$ is defined so that $4\lambda = (s - 1)(s + 1)$ and so that $s
= 2 i r$.

\subsection{\label{UTSG} Unit tangent bundle and space of geodesics}

We  denote by $B=\{z\in\C, |z|=1\}$ the boundary at infinity of
$\D$. The unit tangent bundle $S \D$ of the hyperbolic disc $\D$
is by definition the manifold of unit vectors in the tangent
bundle $T \D$ with respect to the hyperbolic metric. We may, and
will, identify $S\D$ with the unit cosphere bundle $S^*\D$ by
means of the metric. We will make a number of further
identifications:

\begin{itemize}

\item $S \D \equiv PSU(1,1)$. This comes from the fact that
$PSU(1,1)$ acts freely and transitively on $S \D$. We identify a
unit tangent vector $(z, v)$ with a group element $g$ if $g \cdot
(i, (0, 1)) = (z, v)$. Similarly, if we work with the upper half
plane model $\HH$, we have $S \HH \equiv PSL(2, \R)$. We identify
$S\D$, $S\HH$, $PSU(1,1)$, and $PSL(2,\R)$. In general, we work
with the model which simplifies the calculations best.  According
to a previous remark, $S \D$, $ PSU(1,1)$ and $PSL(2, \R)$ will
often be designated by the letter $G$.

\item $S\D \equiv \D \times B$. Here, we identify $(z, b)\in \D
\times B$ with the unit tangent vector $(z, v)$, where $v \in S_z
\D$ is the vector tangent to the unique geodesic through $z$
ending at $b$.

\end{itemize}

The geodesic flow $g^t$ on $S \D$ is defined by $g^t(z, v) =
(\gamma_v(t), \gamma_v'(t))$ where $\gamma_v(t)$ is the unit speed
geodesic with initial value $(z, v)$. The space of geodesics is
the quotient of $S \D$ by the action of $g^t$. Each geodesic has a
forward endpoint $b$ and a backward endpoint $b'$ in $B$, hence
the  space of geodesics of $\D$ may be identified with $B \times B
\setminus \Delta,$ where $\Delta$ denotes the diagonal in $B
\times B$: To $(b', b)\in B \times B \setminus \Delta$ there
corresponds a unique geodesic $\gamma_{b', b}$ whose forward
endpoint at infinity equals $b$ and whose backward endpoint equals
$b'$.

We then have the identification $$S \D \equiv (B \times B
\setminus \Delta) \times \R. $$ The choice of time parameter is
defined -- for instance -- as follows: The point $(b', b, 0)$ is
by definition the closest point to $0$ on $\gamma_{b', b}$ and
$(b', b, t)$ denotes the point $t$ units from $(b, b', 0)$ in
signed distance towards $b$.

\subsection{\label{NEFA} Non-Euclidean Fourier analysis}

Following  \cite{H}, we denote by   $ \langle z, b \rangle $  the
signed distance to $0$ of the horocycle through the points $z \in
\D, b \in B$. Equivalently,
$$e^{\langle z, b \rangle} = \frac{1 - |z|^2}{|z - b|^2} = P_\D(z, b),$$
 where  $P_\D(z, b)$ is the Poisson kernel of the unit disc. (We
 caution again that $e^{\langle z, b \rangle}$ is written $e^{2\langle z, b
 \rangle}$ in \cite{H, Z2}).
We denote Lebesgue measure on $B$ by  $|db|$, so that
the  harmonic measure issued from $0$ is given by $P_\D(z,b)
|db|$. A basic identity (cf. \cite{H}) is that \begin{equation}
\label{ID1} \langle g \cdot z, g \cdot b \rangle = \langle z, b
\rangle + \langle g \cdot 0, g \cdot b \rangle,  \end{equation}
which implies
\begin{equation} \label{ID2} P_\D(g z, g b)\, |d (gb)| = P_\D(z, b)\,
|db|.\end{equation}

The functions $e^{(\frac12+ir)\langle z, b \rangle}$ are
hyperbolic analogues of Euclidean plane waves $e^{i \langle x, \xi
\rangle}$ and are called non-Euclidean plane waves in \cite{H}.
The non-Euclidean Fourier transform is defined by
$$\fcal u(r, b) = \int_{\D} e^{(\frac12-ir)\langle z, b \rangle}
u(z) dVol(z). $$ The hyperbolic Fourier inversion formula is given
by
$$u(z) = \int_B \int_{\R} e^{( \frac12+ir )\langle z, b \rangle} \fcal u(r, b)
r \tanh (2\pi r) dr |db|. $$ As in \cite{Z3}, we  define the
hyperbolic calculus of  pseudo-differential operators $Op(a)$ on
$\D$ by
$$Op (a) e^{( \frac12+ir)\langle z, b \rangle} = a(z, b, r) e^{(\frac12 +ir)\langle z, b
\rangle}. $$ We assume that the complete symbol $a$ is a
polyhomogeneous function of $r$   in the classical sense that
$$a(z, b, r) \sim \sum_{j = 0}^{\infty} a_j(z, b) r^{-j + m}$$
for some $m$ (called its order). By asymptotics is meant that
$$a(z, b, r) - \sum_{j = 0}^{R} a_j(z, b) r^{-j + m} \in S^{m - R - 1}$$
where $\sigma \in S^k$ if $\sup (1 + r)^{j - k}| D^{\alpha}_z D^{\beta}_b
D_r^j(z)| < +\infty$ for all compact set and for all $\alpha,\beta, j$.

 The non-Euclidean Fourier
inversion formula then extends the definition of $Op(a)$ to
$C_c^{\infty}(\D)$:
$$Op(a)u(z)=\int_B \int_{\R} a(z, b, r)e^{(\frac12+ir)\langle z, b \rangle} \fcal u(r, b)
r \tanh(2 \pi r) dr |db|.$$

A key property of $Op$ is that $Op(a)$ commutes with the action of
an element $\gamma\in G$ ($T_\gamma u(z)=u(\gamma z)$)  if and
only if $a(\gamma z, \gamma b, r) = a(z, b, r)$.
$\Gamma$-equivariant pseudodifferential operators then define
operators on the quotient $\X$. This will be seen more clearly
when we discuss Helgason's representation formula for
eigenfunctions.

\subsection{\label{REP} Dynamics and group theory of  $G = SL(2, \R)$}

 We recall the group theoretic point of view towards the geodesic
 and horocycle flows on $S\X$.  As above, it is equivalent to work with $G=PSU(1,1)$ or
 $G = SL(2, \R)$; we choose the latter.
 Our notation follows
 \cite{L, Z}, save for the normalization of the metric.
The generators of $sl(2, \R)$ are denoted by
$$H = \left( \begin{array}{ll} 1 & 0 \\ & \\ 0 & - 1 \end{array}
\right), \;\;\; V = \left( \begin{array}{ll} 0 & 1 \\ & \\
1 & 0
\end{array} \right), \;\; W =  \left( \begin{array}{ll} 0 & -1 \\ & \\1 & 0 \end{array}
\right).$$ We denote the associated one parameter subgroups by $A,
A_-, K$. We denote the raising/lowering operators for $K$-weights
by \begin{equation} \label{EPM} E^+ = H + i V, \;\;\; E^- = H - i
V.  \end{equation}  The Casimir operator is then given by $4 \;
\Omega = H^2 + V^2 - W^2$; on $K$-invariant functions, the Casimir
operator acts as the laplacian $\Lap$. We also put
$$X_+ = \left( \begin{array}{ll} 0 & 1 \\ & \\ 0 & 0 \end{array}
\right),\;\;\;X_- = \left( \begin{array}{ll} 0 & 0 \\ & \\ 1 & 0
\end{array} \right), $$
and denote the associated subgroups by $N, N_-$.

In the identification $S \D\equiv PSL(2, \R)$ the geodesic flow is
given by the right action of the group of diagonal matrices, $A$:
$g^t(g)=ga_t$ where
\begin{equation*} a_t = \left( \begin{array}{ll} e^{t/2} & 0 \\ & \\
0 & e^{-t/2} \end{array} \right). \end{equation*} By a slight
abuse of notation, we sometimes  write $a$ for
$ \left( \begin{array}{ll} a & 0 \\ & \\
0 & a^{-1} \end{array} \right). $ The action of the geodesic flow
is closely related to that of the horocycle flow $(h^u)_{u\in\R}$,
defined by the right action of $N$, in other words by
$h^u(g)=gn_u$ where
\begin{equation*} n_u = \left( \begin{array}{ll} 1 & u \\ & \\
0 & 1\end{array} \right). \end{equation*}

The closed orbits of the geodesic flow $g^t$ on $\Gamma \backslash
G$ are denoted $\{\gamma\}$ and are in one-to-one correspondence
with the conjugacy classes of hyperbolic elements of $\Gamma$.  We
denote by $G_{\gamma}$, resp. $\Gamma_{\gamma}$, the centralizer
of $\gamma$ in $G$, resp. $\Gamma$. $\Gamma_{\gamma}$ is generated
by an element $\gamma_0$ which is called  a primitive hyperbolic
geodesic.  The length of $\gamma$ is denoted $L_{\gamma} > 0$ and
means that
$\gamma$ is conjugate, in $G$, to \begin{equation} \label{agamma} a_{\gamma} = \left( \begin{array}{ll} e^{L_{\gamma}/2} & 0 \\ & \\
0 & e^{-L_{\gamma}/2} \end{array} \right). \end{equation} If
$\gamma = \gamma_0^k$ where $\gamma_0$ is primitive, then
$L_{\gamma_0}$ is the primitive length of the closed geodesic
$\gamma$.

\subsection{Representation theory of $G$ and spectral theory of
$\Lap$}

Let us recall some basic facts about the representation theory of
$L^2(\Gamma \backslash G )$ in the case where the quotient is
compact (cf. \cite{K, L}).


In the compact case, we have the decomposition into irreducibles,
$$L^2(\Gamma \backslash G)  = \bigoplus_{j = 1}^S {\mathcal C}_{ ir_j} \oplus \bigoplus_{j = 0}^{\infty}
\pcal_{ir_j} \oplus \bigoplus_{m = 2, \; m \;even}^{\infty}
\mu_{\Gamma}(m) \dcal_{m}^+ \oplus \bigoplus_{m = 2, m \; even
}^{\infty} \mu_{\Gamma}(m) \dcal_{m}^-,
$$ where ${\mathcal C}_{s_j}$ denotes the complementary series
representation, resp. $\pcal_{s_j}$ denotes the unitary principal
series representation, in which $- \Omega$ equals $s_j(1-s_j ) =
\frac14+r_j^2$. In the complementary series case, $ir_j \in \R$
while in the principal series case $i r_j \in i \R^+$. The
irreducibles are indexed by their $K$-invariant vectors
$\{\phi_{ir_j}\}$, which is assumed to be the given orthonormal
basis of $\Lap$-eigenfunctions. Thus, the multiplicity of
$\pcal_{ir_j}$ is the same as the multiplicity of the
corresponding eigenvalue of $\Lap$.

 Further,
 $\dcal^{\pm}_m$ denotes the holomorphic (resp.
anti-holomorphic) discrete series representation with lowest
(resp. highest) weight $m$, and  $\mu_{\Gamma}(m)$ denotes its
multiplicity; it  depends only on the genus of $\X$. We denote by
$\psi_{m, j}$ ($j = 1, \dots, \mu_{\Gamma}(m))$ a choice of
orthonormal basis of the lowest weight vectors of $\mu_{\Gamma}(m)
\dcal_{m}^+ $ and write $\mu_{\Gamma}(m) \dcal_{m}^+ = \oplus_{j =
1}^{\mu_{\Gamma}(m)} \dcal^+_{m, j}$ accordingly.

We will also use the notations $\ccal_{ir_j}, \pcal_{ir_j}$ and
$\dcal_{m, j}^{\pm}$ for the orthogonal projection operators of
$L^2(\Gamma \backslash G)$  onto these subspaces. Thus, for $f \in
L^2$ we write
\begin{equation} f = \sum_j \ccal_{ir_j}(f) +  \sum_{j} \pcal_{ir_j}(f) + \sum_{m, j,  \pm}
\dcal_{m, j}^{\pm}(f). \end{equation}


By an automorphic $(\tau, m)$-eigenfunction, we mean
 a $\Gamma$-invariant
joint eigenfunction
\begin{equation} \label{CASIMIR} \left\{ \begin{array}{l} \Omega
\sigma_{\tau,
m} = -(\frac{1}{4} + \tau^2) \sigma_{\tau, m} \\ \\
W \sigma_{\tau, m} = i m \sigma_{\tau, m}. \end{array} \right.
\end{equation}
of the Casimir $\Omega$ and the generator $W$ of $K = SO(2)$.

We recall that the  principal series $\pcal^{+}_{ir}$ are realized
on the Hilbert space $L^2(\R)$ by the action
$$\pcal^{\pm}_{ir}\left( \begin{array}{ll} a & b \\ & \\ c & d
\end{array} \right)  f(x) = |-bx + d|^{-1 - 2 ir } f(\frac{ax - c}{- bx +
d}). $$ The  unique normalized $K$-invariant vector of
$\pcal_{ir_j}$ is a constant multiple of
$$f_{ir, 0}(x) =    (1 + x^2)^{-(\frac{1}{2} + ir)}. $$

The complementary series representations are realized on $L^2(\R,
B)$ with inner product
$$\langle B f, f \rangle = \int_{\R \times \R} \frac{f(x)
\overline{f(y)}}{|x - y|^{1 - 2 u}} dx dy $$ and with action
$${\mathcal C}_{u}\left( \begin{array}{ll} a & b \\ & \\ c & d
\end{array} \right)  f(x) = |-bx + d|^{-1 -  2u} f(\frac{ax - c}{- bx +
d}). $$ When asymptotics as $|r_j| \to \infty$ are involved, we
may ignore the  complementary series representations and therefore
do not discuss them in detail.

Let $\C_+ = \{z \in \C: \Im z > 0\}$.  We recall (see \cite{K}, \S
2.6) that $\dcal^+_m$ is realized on the Hilbert space
$${\mathcal H}_m^+ = \{ f \; \mbox{holomorphic on }\; \C_+, \;\; \; \int_{\C_+}
|f(z)|^2 y^{m-2} dx dy < \infty \} $$ with the action
\begin{equation} \label{ACTIONDISC} \dcal_m^+ \left( \begin{array}{ll} a & b \\ & \\ c & d
\end{array} \right) f(z) = (-bz + d)^{-m} f(\frac{az - c}{- bz +
d}).  \end{equation} The lowest weight vector  of $\dcal^+_m$ in
this realization is is $(z + i)^{-m}$.

We note that the $K$-weights in all irreducibles are even. Lowest
weight vectors of $\dcal^+_m$ correspond to (holomorphic)
automorphic forms of weight $m$ for $\Gamma$ in the classical
sense of holomorphic functions on $\C_+$ satisfying
$$f(\gamma \cdot z) = (c z + d)^{m} f(z), \;\;\; \gamma = \left( \begin{array}{ll} a & b \\ & \\ c & d
\end{array} \right), \;\;\; \gamma \in \Gamma. $$
Thus, $f$ is invariant under  the action of $\Gamma$.  A
holomorphic form of weight $m$ defines a holomorphic differential
of type $f(z) (dz)^{\frac{m}{2}}$ (cf. \cite{Sa2}). Forms of weight
$n$ in $\pcal_{ir}, \ccal_{u}, \dcal_{ir}^{\pm}$ always correspond
to differentials of type $(dz)^{\frac{n}{2}}$. Forms of odd weight
do not occur in $L^2(\Gamma \backslash PSL(2, \R))$.

We now consider the action of $A$, i.e. the geodesic flow, in each
irreducible.

\begin{prop}\label{CYCLIC}   The right action of $A$, i.e. the geodesic flow
$g^t$, has two invariant subspaces in each irreducible
$\ccal_{ir}, \pcal_{ir}$, namely the cyclic subspace generated by
the weight zero vector $\phi_{ir}$, and that generated by  $X_+
\phi_{ir}$. The action of $A$ is irreducible in $\dcal^{\pm}_m.$

 \end{prop}

\begin{proof}

 In the principal series we have
$$\pcal^{+}_{ir}\left( \begin{array}{ll} e^{t/2} & 0 \\ & \\ 0 &
e^{-t/2}
\end{array} \right)  f(x) = e^{ t (\frac{1}{2} + ir) } f(e^t x ). $$ The subspaces $L^2(\R_+), L^2(\R_-)$ are
invariant. The generalized eigenfunctions (eigendistributions) are
$x_+^{1/2 + i \xi}, x_-^{1/2 + i \xi}$. The action is irreducible
in each subspace.

In the discrete series we have
$$\dcal_m^+ \left( \begin{array}{ll} e^{t/2} & 0 \\ & \\ 0 &
e^{-t/2}
\end{array} \right)  f(z) = e^{ m t/2}  f(e^t z). $$
The lowest weight vector is cyclic for the action of $A$.

\end{proof}

\subsection{Time reversibility}

Time reversal refers to the involution on the  unit cosphere
bundle defined by  $\iota (x, \xi) = (x, - \xi)$.
 Under the
identification $\Gamma \backslash G \sim S^*X_{\Gamma}$, the time
reversal map takes the form $\Gamma g \to \Gamma w g w$ where $w =
\left(
\begin{array}{ll}0 & 1 \\ & \\ -1 & 0
\end{array} \right)$ is the Weyl element. For $a \in A$ one has $w a w =
a^{-1}$.

 We say that a distribution is
time-reversible if $\iota^* T = T$. The distributions of concern
in this article  all  have the property of time-reversibility,
originating in the fact that $\Lap$ is a real operator and hence
commutes with complex conjugation.  This motivates the
decomposition of $\pcal_{ir} = \pcal_{ir}^+ \oplus \pcal^-_{ir} $
into `even' and `odd' subspaces.

\begin{prop} We have:

\begin{itemize}

\item Each principal (or complementary) series irreducible
contains  a one-dimensional space of $A$-invariant and
time-reversal invariant distributions. In the realization on
$L^2(\R)$, it is  spanned by $\xi_r(x) = |x|^{-(\frac{1}{2} +
ir))}$.

\item  There exists a unique (up to scalars) $A$-invariant
time-reversal invariant distribution in $\dcal_m^+ $ when $m
\equiv 0 (\mbox{mod} 4)$ and there exists no time reversal
invariant distribution when $m \equiv 2 (\mbox{mod} 4)$. In the
realization on ${\mathcal H}_m^+$, it is $z^{-m/2}$. Similarly for
$\dcal_m^-$.
\end{itemize}

\end{prop}

\begin{proof}

\noindent{\bf (i) The complementary and principal series}

Each principal (or complementary) series irreducible contains a
two-dimensional space of $A$-invariant distributions. In the model
on $L^2(\R)$ a basis is  given by $x_+^{-(\frac{1}{2} + ir)},
x_-^{-(\frac{1}{2} + ir) }$. Indeed, $A$ invariance is equivalent
to $$ e^{- t (\frac{1}{2} + ir) } \xi_{ir} (e^t x ) = \xi_{ir}(x).
$$ Setting $x = \pm 1$ we find that
$$\xi_{i r}^{\pm} (x) = x_{\pm}^{-(\frac{1}{2} + ir)}$$ are invariant
distributions supported on $\R_{\pm}$.

The time reversal operator is given by
\begin{equation} \label{TRPR} \pcal_{ir}\left(
\begin{array}{ll} 0 & 1 \\ & \\ -1 & 0
\end{array} \right)  f(x) = |x|^{-1 - 2 ir } f(\frac{- 1}{x}).
\end{equation}
Hence, time reversal invariance is equivalent to
$$f(\frac{- 1}{x}) = |x|^{1 + 2 ir} f(x). $$

 Under time reversal
$$\pcal_{ir} (w) x_+^{-(\frac{1}{2} + ir)} =  |x|^{-1 - 2ir } x_-^{(\frac{1}{2} + ir)} =  x_-^{-(\frac{1}{2} + ir)} . $$
Hence the unique time reversal invariant distribution is
$$\xi_{ir} = |x|^{-(\frac{1}{2} + ir)}. $$

\medskip

\noindent{\bf (ii) The discrete series}

Each holomorphic (or anti-holomorphic) discrete series irreducible
$\dcal^{\pm}_m$  contains a unique (up to scalar multiple)
$A$-invariant distribution $z^{-m/2}$. Indeed, to solve
$$\dcal^{+}_{m}\left(
\begin{array}{ll} e^{t/2} & 0\\ & \\ 0 &
e^{-t/2}
\end{array} \right)  \xi_m^+(z) = e^{m t/2} \xi_m(e^t z) = \xi_m^+(z), $$ we put $z = e^{i \theta}$ and  obtain
$$  \xi_m^+(r e^{i \theta}  ) = r^{-m/2} \xi_m^+(e^{i \theta})), $$
and the only holomorphic solution is  $z^{-m/2}$.

In the holomorphic discrete series, the  time reversal operator is given by $$\dcal^{+}_{m}\left( \begin{array}{ll} 0 & 1 \\ & \\
-1 & 0
\end{array} \right)  f(z) = z^{-m} f(\frac{- 1}{z}). $$
We observe that  $z^{-m/2}$ is time-reversal invariant when $m
\equiv 0 (\mbox{mod} 4)$ and is anti-invariant when $m \equiv 2
\mbox{mod} 4$.

The anti-holomorphic discrete series is similar (by taking complex
conjugates).

\end{proof}

\begin{defn} We
denote the time reversal and geodesic flow invariant distribution
in $\dcal' (\Gamma \backslash G) \cap \pcal_{ir_j}$, normalized so
that $\langle \phi_{ir_j}, \Xi_{ir_j} \rangle = 1$,  by
$\Xi_{ir_j}$. We denote by $\Xi_{m, j}^{\pm}$ the time reversal
and geodesic flow invariant distribution in $\dcal' (\Gamma
\backslash G) \cap \dcal_{m}^{\pm} $, normalized so that $\langle
\psi_{m,j}, \Xi_{m, j} \rangle = 1$, where $||\psi_{m, j}|| = 1$.
Here, we assume $m \equiv 0 (4)$.
\end{defn}

The following is the main application of the representation
theory. By the above normalization, all denominators equal one,
but we leave them in to emphasize the normalization.

\begin{prop}\label{EXPRESSION}  Let $d\nu$ denote a time-reversal invariant and
geodesic flow invariant distribution on $\Gamma \backslash G$. Let
$f \in C^{\infty}(\Gamma \backslash G)$. Then:
$$\begin{array}{lll}\langle f, d\nu \rangle & = &\sum_j \frac{ \langle \pcal_{ir_j}(f), \Xi_{ir_j}
\rangle}{\langle \phi_{ir_j}, \Xi_{ir_j} \rangle} \; \langle
\phi_{ir_j}, d\nu \rangle \\ &&\\
&+& \sum_{\pm, m = 2, m \equiv 0(4)}^{\infty}  \sum_{j =
1}^{\mu_{\Gamma}(m)} \frac{ \langle \dcal_{m,j}^{\pm} f, \Xi_{m,
j}^{\pm} \rangle}{\langle \psi_{m,j}, \Xi_{m,j}^{\pm} \rangle} \;
\langle \psi_{m,j}, d\nu \rangle.
\end{array}
$$

\end{prop}

\begin{proof} Since $\phi_{ir}$ is cyclic in $\pcal_{ir}$, any
element $f$ in this space may be expressed in the form $\int_{\R}
\tilde{f}(t) \phi_{ir} \circ g^t dt$. If we pair with the
invariant distribution $d\nu$ we obtain $\int_{\R} \tilde{f} (t)
dt \; \langle \phi_{ir}, d\nu \rangle$. On the other hand, if we
pair $f$ with $\Xi_{ir}$ we obtain $\int_{\R} \tilde{f}(t) dt \;
\langle \phi_{ir}, \Xi_{ir} \rangle$. Similarly in the discrete
series. The statement follows immediately.

\end{proof}

To apply the Proposition, we need to understand convergence of the
series and hence to have bounds on $\langle \pcal_{ir}(f),
\Xi_{ir} \rangle$ and $ \langle \dcal_{m,j}^{\pm} f, \Xi_{m,
j}^{\pm} \rangle$ when the denominator is normalized to equal one.
Since the complementary series sum is finite, it is not necessary
to analyze these terms. The following proposition shows that the
distributions are of order one. Here, we say that a distribution
$T$ has order $r$ if $\langle T, f \rangle \leq ||f||_{W^s}$ where
$W^s(\Gamma \backslash G)$ is the Sobolev space of functions with
$s$ derivatives in $L^2$. The proposition also controls the
dependence of the norms in the Casimir parameters $ir, m$.

\begin{prop} \label{XIEST} We have:

\begin{itemize}

\item $|\langle \pcal_{ir}(f), \Xi_{r} \rangle | \leq  C r^{-
1/2}\;\; ||\pcal_{ir}(f)||_{W^1}; $

\item $|\langle \dcal_{m,j}^+ f, \Xi_{m, j} \rangle | \leq
 C m^{-1/2} \; ||\dcal_{m,j}^+ f||_{W^1}; $

\end{itemize}

\end{prop}

\begin{proof} We prove the results by conjugating to the models
above.

We begin with the continuous series and let
$${\mathcal U}_{ir} : L^2(\Gamma \backslash G) \to L^2(\R)$$
be the unitary intertwining operator from $\pcal_{ir} \subset
 L^2(\Gamma \backslash G)$ to its realization in $L^2(\R)$.
Thus, $\ucalr \Xi_{ir} = \xi_{ir}$ up to the normalizing constant.

To determine the normalizing constant, we recall (see \cite{Z}, p.
59) that
$$\begin{array}{lll} \langle \ucalr \phi_{ir}, \ucalr \Xi_{ir} \rangle& = &\int_{\R}   (1 + x^2)^{-(\frac{1}{2} + ir)}  |x|^{-(\frac{1}{2} - ir)}dx\\ &&
\\ &=&  2\; \int_0^{\infty}   (1 + x^2)^{-(\frac{1}{2} + ir)} x^{(\frac{1}{2} + ir)} \frac{dx}{x} \\ &&\\
& = &   2 \; \int_0^{\infty}
(x^{-1} + x)^{-(\frac{1}{2} + ir)} \frac{dx}{x} \\ &&\\
& = &  2\; B(\frac{1}{2} (\frac{1}{2} + ir), \frac{1}{2}
(\frac{1}{2} + ir)) :=
 2 \; \frac{\Gamma(\frac{1}{4} + \frac{ir}{2})^2}{\Gamma(\frac{1}{2} +
 ir)}
\end{array}$$
Here, $B(x, y) = \frac{\Gamma(x) \Gamma(y)}{\Gamma(x + y)}$ is the
beta-function. From the asymptotics (cf. \cite{GR} 8.328)
\begin{equation} \label{GAMASYM}  \Gamma(x + i y) \sim \sqrt{2 \pi} e^{- \frac{\pi}{2} |y|} |y|^{x -
\frac{1}{2}}\;\;\;\; (|y| \to \infty)
 \end{equation} of the $\Gamma$-function along vertical lines in $\C$,  it
follows that
$$ \left(\beta_r\right)^{-1} := \frac{\Gamma(\frac{1}{4} + \frac{ir}{2})^2}{\Gamma(\frac{1}{2} +
 ir)} \sim C r^{-1/2}, \;\; (r \to \infty). $$

Next we consider the order of $\xi_{ir}$  as a distribution in the
model. We may break up each function in $L^2(\R)$ into its even
and odd parts with respect to time reversal invariance, and then
we only need to consider $\langle f, \xi_{ir} \rangle$ for a time
reversal invariant  $f$. Let $\chi_+(x) \in C_0^{\infty}$ with
$\chi_+= 1$ for $|x| \leq \frac{1}{2}$ and $0$ for $|x| > 2$ and
with the property that   $\chi_+(x) + \chi_+(\frac{-1}{x}) = 1$.
Then $\langle f, \xi_{ir} \rangle = \langle (\chi_+ +
\chi_+(\frac{-1}{x})) f, \xi_{ir} \rangle$ and (denoting the time
reversal (\ref{TRPR}) operator by $T$)
$$\begin{array}{lll} \langle  \chi_+(\frac{-1}{x}) f, \xi_{ir} \rangle &= &
 \langle T \chi_+(\frac{-1}{x}) f, T\xi_{ir} \rangle \\ && \\
 & = &\langle (T \chi_+(\frac{-1}{x}) T^{-1} )f, \xi_{ir} \rangle
 \\ && \\
 & =  & \langle  \chi_+f, \xi_{ir} \rangle. \end{array}$$
 Hence we only need to estimate the $\chi_+$  integral. We    write $x^{-1/2 + ir} = \frac{1}{-1/2
+ ir}  \frac{d}{dx} x^{1/2 + ir}$ and integrate by parts. The
result is bounded by $C (1 + r)^{-1} (||f||_{L^2} +
||\pcal_{ir}(X_-) f||_{L^2})$. Here, we use that $X_+$ is
represented by $\frac{d}{dx}$.

 It follows that
for $f \in C^{\infty}(\Gamma \backslash G), $
$$\begin{array}{lll} |\frac{\langle f, \Xi_{ir} \rangle}{\langle \phi_{ir}, \Xi_{ir} \rangle}|  = |\langle \pcal_{ir}(f), \Xi_{ir}
\rangle | & = & \beta_r |\langle \ucalr \pcal_{ir}(f), \xi_{ir} \rangle| \\
&&\\
& \leq & \beta_r ||\frac{d}{dx} \ucalr \pcal_{ir}(f)||_{L^2(\R)} \\ && \\
&= & C \beta_r (1 + r)^{-1} ||X_- \pcal_{ir}(f)||_{L^2(\Gamma
\backslash G)} \leq r^{- 1/2}\; ||\pcal_{ir}(f)||_{W^1(\Gamma
\backslash G)}.
\end{array}$$

We now consider the discrete series. The normalizing constant is
now
$$\langle \psi_m, \Xi_m^+ \rangle = \frac{1}{||(z + i)^{-m}||}\int_{\C_+}
(z + i)^{-m} \bar{z}^{-m/2} y^{m-2} dx dy.$$ To calculate the
constant, we use the isometry
$$T_m: {\mathcal H}_m^+  \to \ocal^2(\D, d\nu_m), \;\; T_m f(w) =
 f\left( - i \frac{w + 1}{w - 1} \right) \left( \frac{-2i}{w - 1} \right)^m,$$
where $\ocal^2(\D, d\nu_m)$ are the holomorphic functions on the
unit disc which are $L^2$ with respect to the measure $d\nu_m =
\frac{4}{4^m} (1 - |w|^2)^m \frac{dw d\bar{w}}{(1 - |w|^2)^2}$
(
cf. \cite{L} IX \S 3).

We have $T_m \psi_m = 1$. Note that $1$ is not normalized to have
$L^2$ norm equal to one. It follows that
$$\langle \Xi_m^+, \psi_m \rangle = \frac{4}{4^m \; ||(z + i)^{-m}||} \int_{\D} \left( - i \frac{w + 1}{w - 1}
\right)^{-m/2} \left( \frac{-2i}{w - 1} \right)^m (1 - |w|^2)^m
\frac{dw d\bar{w}}{(1 - |w|^2)^2}.$$ We write $w = r e^{i \theta}$
and observe that the angular integral equals an $r$-independent
constant times
$$ \int_{S^1}  \left(   \frac{1 + r e^{i \theta} }{1 - r e^{i \theta}}
\right)^{-m/2} \left( \frac{2i}{1 - r e^{i \theta}} \right)^m d
\theta = \int_{|z| = 1}   \left(   \frac{1 + r z}{1 - r z}
\right)^{-m/2} \left( \frac{-2i}{1 - r z} \right)^m\frac{dz}{z} =
2 \pi i (-2i)^m,$$ since $\left( \frac{1 + r z}{1 - r z}
\right)^{-m/2}\left( \frac{-2i}{rz - 1} \right)^m$ is holomorphic
in $|z| \leq 1$ for $r < 1$. It follows that
$$\langle \psi_m, \Xi_m^+ \rangle = C  \frac{2^m }{4^m \; ||(z + i)^{-m}||} \int_{\D}  (1 - |w|^2)^m \frac{dw d\bar{w}}{(1 -
|w|^2)^2}= C (m-1)^{-1/2},$$ since the $L^2$-norm  of $T_m \psi_m
= 1$ equals $\frac{2}{2^m} \left(\int_{\D}  (1 - |w|^2)^m \frac{dw
d\bar{w}}{(1 - |w|^2)^2} \right)^{1/2}$ and $ \int_{\D}  (1 -
|w|^2)^m \frac{dw d\bar{w}}{(1 - |w|^2)^2}$ equals
$\frac{1}{m-1}.$

We then need to estimate
$$\langle \Xi_m^+, \bar{f} \rangle =   \int_{\C_+}
f(z)  z^{-m/2} y^{m-2} dx dy.$$ As above, we let $\chi_+$ be a
radial function with compact support in $\R_+$ and with $\chi_+(z)
+ \chi_+(\frac{-1}{z}) \equiv 1$. By unitary of time reversal, we
again have
$$\begin{array}{lll} \langle \chi_+(\frac{-1}{z}) \Xi^+_m, \bar{f}
\rangle &=&  \langle \chi_+ \Xi^+_m, \bar{f} \rangle,
\end{array}$$
and thus it suffices to estimate the $\chi_+$ integral.
 We note that for $m > 2$,  $z^{-m/2} =\frac{1}{1 -
m/2} \frac{d}{dx} z^{-m/2 + 1}$ and that $z^{-m/2 + 1} \in L^2(|z|
< 1, y^m \frac{dx dy}{y^2}). $ The operator $\frac{d}{dx} =
\dcal_m^+(X_-)$ is skew symmetric with respect to the inner
product. Partial integration gives the bound $\frac{1}{1 - m/2}
||f'||_{L^2} $, hence after normalizing $\Xi_m^+$ we have
$$|\frac{\langle f, \Xi_m^+ \rangle}{\langle \psi_m, \Xi^+_m \rangle} | \leq C m^{-1/2} ||f||_{L^2} +
||\dcal_m^+(X_+) f ||_{L^2}).$$

\end{proof}

\begin{rem} The paper \cite{A-P} studies related estimates in the discrete series  from a
different point of view. \end{rem}

A nice simplification is that the series $\{X_+ \phi_{ir_k}\}$
automatically has zero integrals against a time reversal invariant
distribution:

\begin{lem} \label{TRX} If $T \in \dcal'(\Gamma\backslash G)$ is
time-reversible,  then $\langle X_+ \phi_{ir_k}, T\rangle = 0$ for
all $k$. \end{lem}

\begin{proof} We have $$\begin{array}{lll} \langle X_+ \phi_{ir}, T \rangle &=& \langle X_+
\phi_{ir}, \iota^* T \rangle \\ & & \\
& = & \langle \iota^* (X_+ \phi_{ir}), T \rangle \\ &&\\ &=& -
\langle X_+ \phi_{ir}, T \rangle. \end{array}$$


\end{proof}

\section{\label{PATTSUL} Patterson-Sullivan distributions and microlocal lifts}

\subsection{Patterson-Sullivan distributions}

Let us first  recall Helgason's fundamental result about
eigenfunctions of the Laplacian on  $\D$. In the following
theorem, $\phi$ is any solution of $\Lap \phi=-\lambda \phi$
($\lambda=\frac14+r^2$ where $\lambda, r\in\C$). The function
$\phi$, defined on $\D$, is not necessarily automorphic. One says
that $\phi$ has exponential growth if there exists $C>0$ such that
$|\phi(z)|\leq C e^{Cd_\D(0, z)}$ for all $z$.

\begin{thm}\label{Helga} (\cite{H}, Theorems 4.3 and 4.29; see also \cite{He})  Let $\phi$ be an eigenfunction with exponential growth, for the
eigenvalue $\lambda=\frac14+r^2\in \C$. Then there exists a unique
distribution $T_\phi \in {\mathcal D}'(B)$ such that
$$\phi(z) = \int_B e^{(\frac12 + ir) \langle z,
b \rangle }T_\phi(db),$$ for all $z\in\D$ \end{thm} The theorem
extends the classical representation theorem for bounded harmonic
functions to the case of arbitrary eigenvalues. Note that
the kernel $e^{(\frac12 + ir) \langle z,
b \rangle }$ that
appears in the representation theorem for eigenfunctions for the
eigenvalue $\lambda_r$
is the generalized Poisson kernel, $P_\D^{(\frac12 + ir)}(z,b)$.  The distribution $T_\phi$ is called the
(radial) boundary value of $\phi$ and is defined as follows:
expand the eigenfunction into the  ``Fourier series",
\begin{equation} \label{GENSPH}  \phi(z)  = \sum_{ n \in \Z} a_{n}
\Phi_{r, n}(z), \end{equation} in terms of the generalized
spherical functions $\Phi_{r, n}$ defined by
\begin{equation}  e^{(\frac12+i r) \langle z, b \rangle } = \sum_{ n \in \Z}
\Phi_{r, n}(z) b^n,\;\;\; b \in B.
\end{equation} The boundary value of
$\phi$ is the distribution on $B$ defined  by
\begin{equation}  T_{\phi}(db) =  \sum_{n \in \Z} a_{n} b^n |db|.
\end{equation}

As observed in \cite{Z2}, when $\phi_{ir_j}$ is a
$\Gamma$-invariant eigenfunction, the boundary values
$T_{ir_j}(db)$ have the following invariance property:
\begin{equation}\label{CONFORMAL} \begin{array}{ll} \phi_{ir_j}(\gamma z) = \phi_{ir_j} (z) & \implies
e^{(\frac12+ir_j)\langle \gamma z, \gamma b \rangle} T_{ir_j}(d
\gamma b) =
e^{(\frac12+ir_j)\langle z, b \rangle} T_{ir_j} (d b)\\ &  \\
& \implies T_{ir_j}( d\gamma b) = e^{- (\frac12+ir_j) \langle
\gamma \cdot 0, \gamma \cdot b \rangle} T_{ir_j} (d b) \end{array}
\end{equation}
This follows from the  uniqueness of the Helgason representation
(\ref{Helga}) and by the identities (\ref{ID1})-(\ref{ID2}). Hence
there exist $e_{ir_j} \in \dcal'(\Gamma \backslash PSL(2, \R))$
such that
\begin{equation} \label{EJ} e^{(\frac12+ir_j)\langle z, b \rangle} T_{ir_j} (d
b)= e_{ir_j}(z, b) P(z, b) db.
\end{equation}
The distribution $e_{ir_j}$ is horocyclic-invariant and
$\Gamma$-invariant. It may be expanded in a $K$-Fourier series,
$$e_{ir_j}(z, b)= \sum_{n \in \Z} \phi_{ir_j, n}, $$
and it is easily seen (cf. \cite{Z2}) that $\phi_{ir_j, 0} =
\phi_{ir_j}$ and that $\phi_{ir_j, n}$ is obtained by applying the
$n$th normalized (i.e. unitary) raising/lowering to $\phi_{ir_j}$.
More precisely, one applies $(E^{\pm})^n$ (\ref{EPM}) and
multiplies by the normalizing  factor $\beta_{2ir_j, n} =
\frac{1}{(2ir_j + 1 \pm 2n) \cdots (2i r_j + 1 \pm 2)}$. The
regularity of these  distributions was recently studied in
\cite{FF, Co}.

 If we set
$z = 0$, the $K$-Fourier series and $B$-Fourier series coincide
and we get
\begin{equation} \label{EJTJ} T_{ir_j}(db) =   e_{ir_j}(0, b) db= \sum_{n \in \Z} \beta_{s, \pm n} \left( (E^{\pm})^n
\phi_{ir_j}(0) \right) b^n db. \end{equation} Rather than
estimating the regularity of $T_{ir_j}(db)$ using (\ref{EJTJ}),
which would take us too far afield,  we will quote some estimates
of Otal \cite{O} which suffice ( and indeed are better than
necessary) for our applications. Roughly, the say that
$T_{ir_j}(db)$ is the derivative of a H\"older continuous function
$F_{ir_j}$. Since its zeroth Fourier coefficient is non-zero,
$T_{ir_j}(db)$ is not literally the derivative of a periodic
function, but it is the derivative of a function $F_{ir_j}$ on
$\R$ satisfying $F_{ir_j}(\theta + 2 \pi) = F_{ir_j}(\theta) +
C_j$ for all $\theta \in \R$. We follow Otal in calling such a
function $2\pi$-periodic.

 For $0
\leq \delta \leq 1$ we say that a $2 \pi$-periodic function $F: \R
\to \C$ is $\delta$-H\"older if $|F(\theta) - F(\theta')| \leq C
|\theta - \theta'|^{\delta}.$ The smallest constant is denoted
$||F||_{\delta}$ and $\Lambda_{\delta}$ denotes the Banach space
of $\delta$-H\"older functions, up to additive constants.

\begin{thm} \label{O} (\cite{O} Proposition 4) Suppose that $\sigma = \delta
+ i t$ and that $\phi$ is an eigenfunction of eigenvalue
$\sigma (1 - \sigma)$ satisfying $||\nabla
\phi||_{\infty} < \infty$.  Then its Helgason boundary
value $T_{\phi}$ is the derivative of a $\delta$-H\"older
function.

\end{thm}
This estimate holds for arbitrary bounded eigenfunctions. In our
case, it says that $T_{ir_j}$ is the derivative of a H\"older
function, of H\"older exponent $\delta=\frac12$ if $\lambda_j\geq
\frac14$. Otal's proof also shows that the H\"older norm is
bounded by  a power of $r_j$. Related results can be found in
\cite{BR, C, MS, FF, Co}.

We now introduce a ``Patterson-Sullivan" distribution associated
to each automorphic eigenfunction. We denote
$\lambda_0=0<\lambda_1\leq ...$ the spectrum of the laplacian on
$\X$ ($\lambda_j=\frac14+r_j^2$), $(\phi_{ir_j})$ is a given
orthonormal basis of eigenfunctions, and their boundary values are
denoted $(T_{ir_j})$.

\begin{remark} We assume  that these eigenfunctions
are real to obtain time reversal invariant distributions. Aside
from that,  our results are  valid for complex eigenfunctions with
slight modifications.
\end{remark}

\begin{defn}
\label{PS1}The Patterson-Sullivan distribution associated to a
real eigenfunction  $\phi_{ir_j}$ is the distribution on $B \times
B \setminus \Delta$ defined by
$$ps_{ir_j}(db', db):= \frac{T_{ir_j}(db) T_{ir_j}(db')}{|b - b'|^{1 +  2i r_j}}$$
\end{defn}

If $\phi_{ir_j}$ is automorphic, it is easy to check that
$ps_{ir_j}$ is invariant under the diagonal action of $\Gamma$:

\begin{prop}\label{INVARIANT}  Suppose that $\phi_{ir_j}$ is $\Gamma$-invariant, and let $T_{ir_j}$ denote
its radial boundary values.  Then the distribution on $B \times B
\setminus\Delta$ defined by
$$ps_{ir_j}(db', db):= \frac{T_{ir_j}(db) T_{ir_j}(db')}{|b - b'|^{1 +  2i r_j}}$$
is $\Gamma$-invariant and time reversal invariant. \end{prop}
\begin{proof}
It follows from (\ref{CONFORMAL}) that
\begin{equation}\label{TJID}  T_{ir_j}(d\gamma b)T_{ir_j} (d\gamma b') =
e^{- (\frac12+ir_j)\langle \gamma \cdot 0, \gamma \cdot b \rangle}
e^{- (\frac12+ir_j)\langle \gamma \cdot 0, \gamma \cdot b'
\rangle} T_{ir_j}( db)T_{ir_j} (d b'). \end{equation}

We will also need the following identities (cf. \cite{N} (1.3.2)):
\begin{equation} \label{MOREIDS} \begin{array}{l}  |\gamma(x) - \gamma(y)| = |\gamma'(x)|^{\half} |\gamma'(y)|^{\half}
|x - y|\\ \\
1 - |\gamma(x)|^2 = |\gamma'(x)| (1 - |x|^2). \end{array}
\end{equation}
for every $x,y\in\D\cup B$, $\gamma\in\Gamma$.  Hence for $b\in B$
and $\gamma\in\Gamma$, we have
\begin{equation} |\gamma(0)-\gamma(b)|^2=|\gamma^\prime(b)|(1-|\gamma(0)|^2).
\end{equation}
Furthermore,
\begin{equation}\label{BID}  |\gamma b - \gamma b'|^2 = e^{- [ \langle \gamma \cdot 0, \gamma \cdot b \rangle +
\langle \gamma \cdot 0, \gamma \cdot b' \rangle]} |b - b'|^2.
\end{equation}
Raising (\ref{BID}) to the power $\frac{1}{2} + ir_j$, taking the
ratio with (\ref{TJID})  and simplifying completes the proof of
$\Gamma$-invariance.

Time-reversal invariance is invariance under $b \iff b'$, which is
obvious from the formula.
\end{proof}

We now construct from the  distribution  $ps_{ir_j}$ a geodesic
flow invariant distribution on $ S \D$ as follows. As reviewed in
\S \ref{BACKGROUND}, the unit tangent bundle $S\D$ can be
identified with $(B \times B \setminus \Delta)\times\R$: the set
$B \times B \setminus \Delta$ represents the set of oriented
geodesics, and $\R$ gives the time parameter along geodesics. We
then define the Radon transform:
\begin{equation} \label{RT} \rcal: C_0(S \D) \to C_0(B \times
B \setminus \Delta), \;\;\; \mbox{by}\;\; \rcal f(b', b) =
\int_{\gamma_{b',b}} f dt.
\end{equation}

Further, we need to define special cutoffs which have the property
that
\begin{equation} \int_{\dcal} f dVol(z) = \int_{\D} \chi f dVol(z)
\end{equation}
for any $f \in C(\Gamma \backslash \D)$, where $\dcal$ is a
fundamental domain for $\Gamma$ in $\dcal$. In other words, $\chi$
is a smooth replacement for the characteristic function of
$\dcal$.

\begin{defn} \label{SFD}  We say that  $\chi \in C_0^{\infty}(\D) $ is a smooth fundamental domain cutoff
if  $$\sum_{\gamma\in\Gamma}\chi(\gamma z)=1.$$\end{defn}

We then make the basic definitions:

\begin{defn} \label{PSDEF}

\begin{enumerate}

\item On $S \D$ we  define the Patterson-Sullivan  distribution
$PS_{ir_j} \in \dcal'(S \D)$ by:
$$PS_{ir_j}(db',db,d t)=ps_{ir_j}(db', db)|dt|$$
in the sense that
$$\int_{S \D} a dPS_{ir_j}  = \int_{B \times B \setminus \Delta} (\rcal a)(b',b) dps_{ir_j} .$$

\item On the quotient $S\X=\Gamma\backslash S\D=\Gamma\backslash
PSU(1,1)$, we define the automorphic Patterson-Sullivan
distributions as follows: the distribution $PS_{ir_j}^{\Gamma} \in
\dcal'(S \X)$ by
$$\int_{S\X}  a  dPS_{ir_j}^{\Gamma} = \int_{S \D} (\chi a)
dPS_{ir_j} = \int_{B \times B \setminus \Delta} \rcal(\chi
a)(b',b) dps_{ir_j},$$ where $\chi$ is a smooth fundamental domain
cutoff. To emphasize the cutoff we sometimes write $PS^{\chi}_j(a)
=  \int_{S \D} (\chi a) dPS_{ir_j}$.
\end{enumerate}

\end{defn}

The following proposition is obvious from the definition, but
important:

\begin{prop}  $PS_{ir_j}$ is a geodesic flow invariant and
$\Gamma$-invariant distribution on $S \D=\D\times B$; and
$PS_{ir_j}^\Gamma$ is geodesic flow invariant on $S\X$.
\end{prop}

The geodesic flow invariance of $PS_{ir_j}$ on $S\D$ is trivial;
for  $PS_{ir_j}^\Gamma$, it is also easy and results from the
following principle:

\begin{lem} \label{inv-principle} Let $T \in
\dcal'(\Gamma \backslash S\D)$ be a $\Gamma$-invariant
distribution. Let $a$ be a $\Gamma$-invariant smooth function on
$S\D$. Then, for any $a_1, a_2\in \dcal (S\D)$ such that
$\sum_{\gamma\in\Gamma}a_i(\gamma.(z,b))=a(z,b)$ ($i=1,2$) we have
$$\int_{S\D} a_1 dT= \int_{S\D} a_2 dT$$
\end{lem}

\begin{proof} Let $\chi$ be a function on $C_0^\infty(\D \times B)$
such that $\sum_{\gamma\in\Gamma}\chi(\gamma.(z,b))\equiv 1$ (in
general, we choose $\chi$ to be independent of $b$). For any such
$\chi$ we have
$$\begin{array}{lll} \int_{S\D}  a_i dT &= &\int_{S\D}
\left\{\sum_{\gamma\in\Gamma}\chi(\gamma
(z, b))\right\}a_i (z, b)T(dz, db) \\ & & \\
& = & \int_{S\D} \sum_{\gamma\in\Gamma} \chi( z,
b)a_i(\gamma(z, b)) T(dz, db)
\\ & & \\
& = & \int_{S\D} \chi( z, b) a(z, b) T(dz, db). \end{array}
$$
\end{proof}

If we look at the expression
\begin{equation} PS_{ir_j}(a) = \int
 |b - b'|^{-1 - 2 ir_j} \rcal(a) T_{ir_j}(db)
T_{ir_
j}(db'), \end{equation} and apply Otal's theorem saying that
$T_{ir_j}=F^\prime_{ir_j}$ for some H\"older function $F_{ir_j}$,
we easily derive:

\vspace{.3cm}
For any $a \in C^{\infty}(S \D)$ we have
$$|PS_{ir_j}(a)| \leq ||F_{ir_j}||^2_{L^\infty(B)}.\, || \frac{\partial^2}{\partial b
\partial b'} |b - b'|^{-1 - 2 i r_j} \rcal(a)||_{L^{\infty}(B \times
B\setminus\Delta)} $$
Of course, the right side may be infinite.

For future reference,
we state a sufficient condition to obtain a non-trivial estimate:

\begin{prop} \label{BASICBB}
Assume that  $|b - b'|^{-1 - 2 i r_j} \rcal(a) \in C^2(B \times
B)$. Then $$|PS_{ir_j}(a)|\leq ||F_{ir_j}||^2_{L^\infty(B)}.\, ||
\frac{\partial^2}{\partial b
\partial b'} |b - b'|^{-1 - 2 i r_j} \rcal(a)||_{L^{\infty}(B \times
B\setminus\Delta)} <
\infty.$$
\end{prop}

A simple example where the condition holds is where $a \in C^2_c(S
\D)$. In that case, there exist $C>0$ and $K>0$ such that:
\begin{equation}\label{normofPS}  |PS_{ir_j}(a)|\leq C(1+|r_j|)^K||a||_{C^2}
\end{equation} for all $j$. If
$a \in C^2(S \X)$, $|PS_{ir_j}^{\Gamma}(a) | \leq
(1+|r_j|)^K||a||_{C^2}$ for all $j.$

\subsection{Microlocal lift and Wigner distributions}

We now give a precise definition of the matrix elements $\langle
Op(a) \phi_{ir_j}, \phi_{ir_j} \rangle$ and hence of the Wigner
distributions.  When $a$ is a $\Gamma$-invariant function on
$S\D$, then in the non-Euclidean calculus \S \ref{NEFA}  we have
\begin{equation} \label{OP} Op(a)\phi_{ir_j}:=\int_B a(z,b)e^{(\frac12
+ ir_j) \langle z, b \rangle }T_{ir_j}(db).
\end{equation}

\begin{defn}\label{WIGNER} The Wigner measure of $\phi_{ir_j}$ is the distribution $W_{ir_j}$ on
$S\X=\Gamma\backslash S\D$ defined by:
$$\int_{S\X} a(g)W_{ir_j}(dg): =\langle Op(a)\phi_{ir_j}, \phi_{ir_j}\rangle_{L^2(\X)},$$
where $Op(a)$ is given by (\ref{OP}).
\end{defn}

To see that $W_{ir_j}$ is a distribution of finite order, we note
that $\langle Op(a)\phi_{ir_j}, \phi_{ir_j}\rangle_{L^2(\X)}$ is
bounded by the operator norm of $||Op(a)||$ and hence by a $C^k$
norm of $a$. In fact, Otal's regularity theorem shows that it is
of order $2$ at most.

Wigner distributions are fundamental in the theory of
quantum ergodicity. Let us
recall the basic result:

\begin{thm} \cite{Sh, Z} Let $d\mu$ denote Haar measure on $S
\X$. Then
$$\frac{1}{N(\lambda)} \sum_{j: |r_j|\leq \lambda} |\int_{S \X} a
dW_{ir_j} - \frac{1}{\mu(S \X)} \int_{S \X} a d\mu|^2 \to 0, $$
where $N(\lambda)$ is the normalization factor $\sharp\{j:
|r_j|\leq \lambda\}$.
\end{thm}

It follows that a subsequence $(W_{j_k})$ of density one of the
Wigner distributions tends to Liouville measure (which equals
normalized Haar measure in this case). The ``quantum unique
ergodicity'' problem is to know whether there exist exceptional subsequences
with other limits. Recently, E. Lindenstrauss proved that no such
exceptional sequences exist in the case of Hecke eigenfunctions on
arithmetic surfaces \cite{L}.

\section{Proof of Theorem \ref{mainintro}}

We begin the proof witha lemma giving the explicit expression of
$W_{ir_j}$:

\begin{lem} \label{COSHPROP} We have
\begin{multline}\langle Op(a) \phi_{ir_j}, \phi_{ir_j}\rangle_{L^2(\X)}\\
 = 2^{( 1+2ir_j)} \int_{B \times B} \left( \int_{\D}
\chi a(z,b) [\cosh s_{b', b}(z) ]^{-(1+2ir_j)}dVol(z)\right)
\frac{T_{ir_j}(db) T_{ir_j}(db')}{|b - b'|^{1 + 2 i r_j}},
\end{multline}
where $\cosh s_{b_1, b_2}(z)$ is given by (\ref{COSHS}). The right
hand side is independent of the choice of $\chi$.
 \end{lem}

\begin{proof} By the generalized Poisson formula and the definition of
$Op(a)$,
$$\langle Op(a) \phi_{ir_j}, \phi_{ir_j}\rangle = \int_{B \times B} \left( \int_{\D} \chi a(z,b) e^{(\frac12 +ir_j) \langle z, b
\rangle} e^{(\frac12+ir_j) \langle z, b' \rangle
}dVol(z)\right)T_{ir_j}(db) T_{ir_j}(db').$$ We then use the
following identity
\begin{lem}\cite{N} \label{COSHS} Let $z \in D$, let $b_1, b_2 \in B$ and let $s_{b_1, b_2}(z)$ denote the hyperbolic
distance from $z$ to the geodesic $\gamma_{b_1, b_2}$ defined by
$(b_1, b_2).$ Then
$$\cosh s_{b_1, b_2}(z) = \frac{2|z - b_1||z - b_2|}{|b_1 - b_2| (1 - |z|^2)}. $$ \end{lem}
Combined with (\ref{BID}) and (\ref{MOREIDS}), we get
$$ e^{ \langle z, b \rangle} e^{ \langle z, b' \rangle}
= 4[\cosh s_{b', b}(z) ]^{-2} |b- b'|^{-2}.$$ Raising both sides
to the power $\frac{1}{2} + i r_j$ completes the proof.
\end{proof}

The next step is to analyze the  integral operator
\begin{multline}\label{phase}\int_{\D} \chi a(z,b) e^{(\frac12+ir) \langle z, b
\rangle} e^{(\frac12+ir) \langle z, b' \rangle } dVol(z) \\=
2^{-(1+2ir)}\int_{\D} \chi a(z,b) [\cosh s_{b', b}(z) ]^{-(1+2ir)}
|b - b^\prime|^{-(1+2ir)}dVol(z).
\end{multline}
In this paragraph -- and probably also in the rest of the paper -- we
sometimes drop the $j$-indices of $r_j$, indexing the
eigenfunctions by $r$ instead.

If we drop the factor $|b - b^\prime|^{-(1+2ir)}$, the  right side
of \ref{phase} defines  the operator  $\llcal_r: C_c(\D) \to C(B
\times B)$ by
\begin{equation} \label{LCAL}  \llcal_r(\chi a)(b', b) =  \int_{\D} \chi
a(z,b) [\cosh s_{b', b}(z) ]^{-(1+2ir)} dVol(z). \end{equation} We
now rewrite the integral in terms of coordinates $z=(t, u)$ based
on the geodesic $\gamma_{b', b}$, after which we can relate
$\llcal_r$  with the operator in (\ref{FIRSTLR}).

Given a geodesic $\gamma_{b', b}$, we work with special
coordinates on $\D$ or $\HH$, adapted to $\gamma_{b', b}$ as
follows. We write $z=(t,u)$ where $t$ measures arclength on
$\gamma_{b',b}$ and $u$ measures arclength on horocycles centered
at $b$. More precisely, we denote by $g(b',b)$ the vector on
$\gamma_{b',b}$ which is closest to the origin, and the
coordinates $(t, u)$ parametrizing $z$ are defined by
$(z,b)=g(b',b)a_t n_u$. For any given $(b',b)$, the volume element
of $z$ is $dVol =|dtdu|$. The computation is most easily checked
in the upper half plane, with $b = \infty, b' = i$ and $g(b', b) =
e=(i, \infty).$ Then $a_t n_u i = e^{t}(i + u)$.  The area form is
$\frac{dx dy}{y^2}$. Setting $y = e^t, x = u e^t$ shows that the
area form equals $|dt du|$.

 We
obtain
\begin{equation} \label{LRINT} \llcal_r(\chi a)(b',b) = \int \cosh s_{b',b}(t,u)^{-(1+2ir)}\chi
a(g(b',b)a_t n_u)|du dt|.
\end{equation}

We further simplify as follows:

\begin{lem}\label{LEMMA}  We have
$$\llcal_r(\chi a)(b, b') = \int_{\R \times \R} (1 +
u^2)^{-(\frac12 + ir)}  \chi a(g(b,b^\prime)a_t n_u)|du dt|. $$
\end{lem}

\begin{proof}

We recall that  $s_{b',b}(t,u)$ is  the distance from the
basepoint of $ga_tn_u$ to the geodesic generated by $g$ in the
hyperbolic plane $\HH=G/K$. That distance depends only on $u$ and has
the value $
\cosh s_{b',b}(t,u) =\sqrt{1+u^2}$.

\end{proof}

Next, we further  rewrite the operator $\llcal_r$  in terms of the
operator $L_r$ in (\ref{FIRSTLR}):

\begin{lem} \label{CORPS} We have:
$$\begin{array}{lll} \langle Op(a) \phi_{ir}, \phi_{ir} \rangle
& = &2^{-(1+2ir)} \int_{G} L_r(\chi a)(g) PS_{ir}(dg).
\end{array}$$

\end{lem}

\begin{proof}

Lemma \ref{LEMMA}  states that
$$\begin{array}{lll}\llcal_r(\chi a)(b, b') &= & \int_{\R} L_r  (\chi a) (g(b,b^\prime)a_t)| dt|\\ & &
\\ & = & \rcal (L_r (\chi a))(b, b').
\end{array} $$

Integrating against $dps_{ir}$ and using the formula in Definition
\ref{PSDEF} completes the proof.

\end{proof}

The next step is to apply the  stationary phase method to
$L_r(\chi a)$. The stationary phase set of (\ref{LRINT}) is the
geodesic $\gamma_{b', b}$ from $b^\prime$ to $b$ or
equivalently it is the  set $u=0$ in the integral defining $L_r(\chi a)$.
Stationary phase method gives the asymptotic expansion

\begin{equation}
\label{statph} L_r(\chi a)(g)  = (-4i\pi r)^{-1/2}\big( \sum_{n\geq
0}r^{-n} L_{2n}( \chi a)(g)\big)
\end{equation} where $L_{2n}$ is a
differential operator of order $2n$ on $S\D$: $L_0$ is the
identity, the other $L_{2n}$ are differential operators in the
stable direction, that is, in the direction $n_u$ generated
by the vector field $X_+$.

If we now integrate \ref{statph} with respect to $PS_{ir}$, and
compare with Lemma \ref{CORPS}, we get an asymptotic expansion,

\begin{equation}\langle Op(a) \phi_{ir}, \phi_{ir}\rangle
=2^{-(1+2ir)} (-4i\pi r)^{-1/2}\big( \sum_{n\geq 0}r^{-n}
\int_{S\D}L_{2n}( \chi a)(g)PS_{ir}(dg) \big)
\end{equation}

Because the distribution on the left-hand side, $e^{(\frac12+ir)
\langle z, b \rangle} e^{(\frac12+ir) \langle z, b' \rangle }
dVol(z)T_{ir}(db)T_{ir}(db^\prime) $, is $\Gamma$-invariant (as a
distribution in the triple $(b,b^\prime, z)$), each of the
distributions obtained in the expansion,
$$f\mapsto \int_{S\D}L_{2n}(f)(g)PS_{ir}(dg) , $$
is $\Gamma$-invariant. In application of the principle
\ref{inv-principle}, the functional
$$a\mapsto\int_{G}L_{2n}(\chi a)(g)PS_{ir}(dg)$$
defines a distribution on $\Gamma\backslash G$, and the definition
does not depend on the choice of $\chi$. The first term ($n=0$) is
precisely the Patterson-Sullivan distribution $PS_{ir}$ as defined
in the quotient $S\X$.

\subsection{Completion of Proof of Theorem \ref{mainintro}}

We now turn to the relation between $W_{ir}$ and $PS_{ir}$. It
follows from the stationary phase asymptotics above,
\eqref{statph}, that
$$\int_{S\X} a(g)W_{ir_j}(dg)=2^{-(1+2ir)}\sum_{n=0}^{N}r_j^{-n}\int_{S\D} L_{2n}(\chi a)(g) \frac{PS_{ir_j}(dg)}{(-4i\pi r_j)^{1/2}}+O(r_j^{-N-1+K})$$
where $K$ was defined in \ref{normofPS}. If we choose $N>K$
then the remainder term goes to zero. Since $L_0=Id$, the operator
$L_r^{(N)}=\sum_{n=0}^{N}r^{-n}L_{2n}$ can be inverted up to
$O(r^{-N-1})$, that is, one can find  differential operators
$M_r^{(N)}=\sum_{n=0}^{N}r^{-n}M_{2n}$ (with $M_0=Id$) and
$R_r^{(N)}$ such that
$$L_r^{(N)}M_r^{(N)}=Id+r^{-N-1}R_r^{(N)}.$$

We thus get
\begin{multline*}\int_{S\X} M_{r_j}^{(N)}a(g)W_{ir_j}(dg)=\int_{S\D} L_{r_j}^{(N)}\chi M_{r_j}^{(N)}
a(g)\frac{PS_{ir_j}(dg)}{(-4i\pi
r_j)^{1/2}}+O(r_j^{-N-1+K})\\=\int_{S\D} L_{r_j}^{(N)}
M_{r_j}^{(N)}\chi
a(g)\frac{PS_{ir_j}(dg)}{(-4i\pi r_j)^{1/2}}+O(r_j^{-N-1+K})\\
=\int_{S\X} a(g)\frac{PS_{ir_j}(dg)}{(-4i\pi
r_j)^{1/2}}+O(r_j^{-N-1+K})
\end{multline*}
The second line is a consequence of \ref{inv-principle}. Since we know,
from standard estimates on pseudo-differential operators, that the Wigner
measures are uniformly bounded in $(C^k)^*$ for some $k$, we have
$$\int_{S\X} M_{r_j}^{(N)}a(g)W_{ir_j}(dg)=\int_{S\X} a(g)W_{ir_j}(dg)+O(r_j^{-1}).$$

This shows that
$$2^{-(1+2ir)}\int_{S\X} a(g)\frac{PS_{ir_j}(dg)}{(-4i\pi r_j)^{1/2}}=\int_{S\X} a(g)W_{ir_j}(dg)+O(r_j^{-1})$$
and finishes the proof of Theorem \ref{mainintro}.

\section{\label{INTOPS} Integral operators and eigenfunctions}

In this section, we give further results on the operators $L_r$
(\ref{FIRSTLR}) and $I_{PS_{ir}}$ (\ref{IMUPSint}) which will be
needed in the proof of Theorem \ref{main2}. With no extra work, we
treat general  integral operators of the form
\begin{equation} \label{IMU} I_{\mu}(\sigma)(u):=\int_{\Gamma
\backslash G}
 \sigma(g n_u) \mu(dg), \end{equation}
 where  $\sigma \in C^{\infty}(\Gamma \backslash G)$ is an automorphic form and where $d \mu$ is an
 invariant distribution for the geodesic flow
 on $\Gamma \backslash G$. In addition to $\mu = PS_{ir}$ the case
 where  $\mu$ is a periodic orbit measure is also important in this article.
 In this case,  we write  $I_{\mu} =
I_{\gamma}$ with$I_{\gamma}(\sigma)(n_u) = \int_{\langle
L_{\gamma} \rangle \backslash A} \sigma(\alpha_{\gamma}^{-1} a
n_u) da$. Here, $\alpha_{\gamma} \in G$ is an element conjugating
$\gamma \in \Gamma$ to an element of $A$. This expression arose in
the trace formulae of \cite{Z} and will arise in \S
\ref{GSELBERG}. The similarity  of these two kinds of integral
operators may be seen as one of the deus ex machina behind Theorem
\ref{main2}.

\subsection{The integral operator $I_{\mu}$}

We can view $I_{\mu}$ as an integral operator from
$C^{\infty}(\Gamma \backslash G) \to C^{\infty}(N) \simeq
C^{\infty}(\R)$ The following lemma shows that when $\sigma$ is a
joint eigenfunction of the Casimir operator and of $W$, then
$I_{\mu}(\sigma)$ solves  an ordinary differential equation in
$u$. When $\sigma$ is a $(\tau, m)$-eigenfunction in the
complementary or principal unitary series, the equation is
\begin{equation}\label{diff}(u^2+1)\frac{d^2 f}{du^2}+(2u-im)\frac{df}{du}+ (\frac14+\tau^2) f=0\end{equation}
We denote by  $F_{\tau,m}\left(\frac{u-i}{-2i}\right)$ the even
solution of (\ref{diff}) which equals $1$ at $u = 0$, and  by
$G_{\tau,m}\left(\frac{u-i}{-2i}\right)$   the odd solution whose
derivative equals $1$ at $u = 0$. In the holomorphic discrete
series, and when $\sigma$ is the lowest weight vector,   the
analogous equation is the  first order equation
 \begin{equation} \label{diff2} 2 i
\frac{df}{du}  = (- 2 \frac{d}{du} - m) f.
\end{equation}
A basis for its solutions is given by $f (u) =
(-i)^{-m/2}(u+i)^{-m/2}$.  There are similar equations for higher
weights and for the anti-holomorphic discrete series, but for
simplicity we only discuss the lowest weight case.

\begin{prop} \label{ISIGMA} Let  $\mu$ be a geodesic flow
invariant distribution on $\Gamma \backslash G$.

\begin{itemize}

\item  Let $\sigma$ be a $(\tau, m)$-eigenfunction in the
principal or complementary series.  Then $I_{\mu}(\sigma)(u)$ is a
solution of (\ref{diff}).  Hence, $$I_{\mu}(\sigma)(u)=
\big(\int_{\Gamma \backslash G}
 \sigma d\mu \big)F_{\tau,m}\left(\frac{u-i}{-2i}\right)
+\big(\int_{\Gamma \backslash G}
 X_+\sigma d\mu \big)G_{\tau,m}\left(\frac{u-i}{-2i}\right),
 $$
 where $F, G$ are the fundamental solutions of (\ref{diff}) defined in \cite{Z}
 (2.3) (see (\ref{WTZERO}) for formulae in terms of hypergeometric
 functions)).

 \item Let $\sigma$ be a $(\tau, m)$-eigenfunction in the discrete holomorphic or
 anti-holomorphic series. For simplicity, assume $\sigma=\psi_m$
 (the lowest weight vector in $\dcal_m^+$). Then:

$$I_{\mu}(\sigma)(u)=\big(\int_{\Gamma \backslash G} \sigma d\mu \big) (-i)^{-m/2}(u+i)^{-m/2}.$$

\end{itemize}

\end{prop}

\begin{proof} In the case of $I_{\mu} = I_{\gamma}$, the proof is given in
\cite{Z}, Proposition 2.3 and Corollary 2.4. We briefly verify
that the same  proof works for any invariant distribution.

First assume $\sigma $ is a $(\tau, m)$-eigenfunction in the
continuous series. Since $4 \Omega = H^2 + 4 X_+^2 - 2 H - 4 X_+
W$ we find that \begin{equation} \label{H2} \left(4
\frac{d^2}{du^2} - 4 i m \frac{d}{du} + 4 (\frac{1}{4} +
\tau^2)\right) I_{\mu}(\sigma) (u) = - \int_{\Gamma \backslash G}
\big((H^2 - 2 H) \sigma\big) (g n_u) \mu(dg). \end{equation} We
write $H \sigma(g)$ as $2 \frac{d}{dt}_{t = 0} \sigma(g a_t)$.
Using that $n_u a_t = a_t n_{u e^{- 2 t}}$ and that $\mu$ is an
$A$- invariant distribution, we find  that $ \int H\sigma(g a_t)
\mu(dg) = - 2 u \frac{d}{du} I_{\mu}(\sigma)(u). $ A similar
calculation replaces $H^2$ by the square of this operator. The
final equation is as stated above. We then evaluate
$I_{\mu}(\sigma)$ and its first derivative at $u = 0$ to obtain
the expression in terms of $F, G$.

In the discrete holomorphic series, we use that $E_- \sigma = 0$
to get $2 i X_+ \sigma = (H - m)\sigma$. This leads to equation
(\ref{diff2}) and to the solution given above.

\end{proof}

\subsection{\label{SPECIALVALUES} The integral $\int_{\R} (1 + u^2)^{-s} I_{\mu}(\sigma)(u) d u$}

In Theorems \ref{maincont}, \ref{main2} and elsewhere, we will
need explicit formulae for the integrals
\begin{equation}\label{STRANS}  \int_{\R} (1 + u^2)^{-s} I_{\mu}(\sigma)(u) d u \end{equation}  We assemble
the results here for future reference.

In view of Proposition \ref{ISIGMA}, we need explicit formulae for
the integral of $(1 + u^2)^{-s}$ against the functions
$F_{\tau,m}\left(\frac{u-i}{-2i}\right),
G_{\tau,m}\left(\frac{u-i}{-2i}\right),$ and
$(-i)^{-m/2}(u+i)^{-m/2}$. In fact, by Proposition \ref{CYCLIC}
 and \ref{TRX}, it will suffice for
Theorems \ref{maincont} and \ref{main2} to have explicit formulae
just for $F_{\tau, 0}$  and $(-i)^{-m/2}(u+i)^{-m/2}$.

We  introduce the following notation:

\begin{equation} \label{BIGDISPLAY} \left\{ \begin{array}{lll} \mu_0(s)& := & 2
\int_0^{\infty} (u^2 + 1)^{-s} du, \\ & & \\ \mu_{ir_k}^c (s) & :=
& \int_{\R}  (u^2 + 1)^{-s} F(\frac{1}{4} + \frac{2 i r_k}{4},
\frac{1}{4} - \frac{2 i r_k}{4}, \frac{1}{2}, - u^2) du,\\ && \\
\mu_{\tau, m}^c(s)& := &
\int_{-\infty}^{\infty} (u^2 + 1)^{-s} F_{\tau, m}(\frac{u-i}{-2i})du,  \;\;\; \\
&
 &\\
\mu_{\tau, m}^{c odd} (s) & := & \int_{\R}  (u^2 + 1)^{-s}
G_{\tau,m}\left(\frac{u-i}{-2i}\right) du,
\\ &&\\
\mu_m^d(s) & := & \int_{\R} (u + i)^{-m/2} (u^2 + 1)^{s} du, \;\;
;
\end{array}\right. \end{equation}

It is clear that the integrals defining $\mu_0(s)$ and
$\mu_m^d(s)$ converge absolutely for $\Re e\;s > \frac{1}{2}$ and
$\Re (2s - \frac{m}{2}) < -1$, respectively. We now show:

\begin{prop}\label{MURKEST}  The
integral defining $\mu_{ir_k}^c (s)$  converges absolutely for $
 -2\Re
s-1/2+\Re i r_k <-1  $, and in this region we have:
$$|\mu_{ir_k}^c (s)| \leq C\; \int_{-\infty}^{\infty} (|u| + 1)^{ -2\Re
s-1/2+\Re i r_k}, $$ for some constant $C$ (independent of $s,
r_k$).
\end{prop}

\begin{proof}

Indeed, as   in \cite{Z} (Proposition 2.7),  the differential
equation (\ref{diff}) is equivalent, by a change of variables, to
a hypergeometric equation, and a short calculation shows that
\begin{equation} \label{WTZERO} \left\{ \begin{array}{l}
F_{ir_k,0}\left(\frac{u-i}{-2i}\right)= F(\frac{1}{4} + \frac{2i
r_k}{4}, \frac{1}{4} - \frac{2i r_k}{4}, \frac{1}{2}, - u^2),
\\ \\
G_{ir_k,0}\left(\frac{u-i}{-2i}\right)= F(\frac{3}{4} + \frac{2i
r_k}{4}, \frac{3}{4} - \frac{2i r_k}{4}, \frac{3}{2}, - u^2).
\end{array} \right. .\end{equation}
 Classical estimates
on hypergeometric functions (see also  \cite{Z}, p. 50)
 show that there exists $C > 0$ (independent of $r_k$) such that
\begin{equation} \label{HALFDECAYR} \left\{
\begin{array}{l} F(\frac{1}{4} + \frac{2i r_k}{4}, \frac{1}{4} -
\frac{2i r_k}{4}, \frac{1}{2}, - u^2) \\ \\ F(\frac{3}{4} +
\frac{2i r_k}{4}, \frac{3}{4} - \frac{2i r_k}{4}, \frac{3}{2}, -
u^2 \end{array} \right. \leq C \; (1 + |u|)^{-1/2 + \Re i r_k},
\end{equation}  These estimates  follow immediately from the
connection formulae for hypergeometric functions:
$$\begin{array}{lll} F(a, b, c, z) &= & \frac{\Gamma(c) \Gamma(b-a)}{\Gamma(b) \Gamma (c - a)}
(-z)^{-a} F(a, 1 + a - c, 1 + a - b; z^{-1}) \\ && \\
&& +  \frac{\Gamma(c) \Gamma(a-b)}{\Gamma(a) \Gamma (c - b)}
(-z)^{-b} F(b, 1 + b - c, 1 + b - a; z^{-1}). \end{array}$$ Since
$F(0) = 1$, we obtain that (as $|u| \to \infty$)  \begin{equation}
\label{HYPER}
\begin{array}{lll} F(\frac{1}{4} + \frac{2ir_k}{4}, \frac{1}{4} -
\frac{2ir_k}{4}, \frac{1}{2}, - u^2) & \sim  &
\frac{\Gamma(\frac{1}{2}) \Gamma(- i r_k )}{\Gamma(\frac{1}{4} -
\frac{i r_k}{2})^2} |u|^{- (\frac{1}{2} +
i r_k)}  \\ && \\
&& +  \frac{\Gamma(\frac{1}{2}) \Gamma(   i
r_k)}{\Gamma(\frac{1}{4} + \frac{2i r_k}{4} )^2}
|u|^{-(\frac{1}{2} - ir_k)}.
\end{array} \end{equation}
The asymptotics (\ref{GAMASYM}) of the $\Gamma$ function on
vertical lines shows that the ratios of $\Gamma$ functions are
uniformly bounded in $r_k$. The decay rate $|u|^{-(\frac{1}{2} -
ir_k)}$ is sufficient for the absolute convergence of the integral
in (\ref{BIGDISPLAY}) as long as $\Re (\frac{1}{2} - ir_k) > 0$,
i.e. if $ir_k$ is not the parameter of the trivial representation.

\end{proof}

Although we will not need them, we note that the estimates for $G$
are similar. Each of the above functions admits a meromorphic
continuation to $\C$. Since we will not need the results for
general $\mu_{\tau, m}^{c} (s), \mu_{\tau, m}^{c odd} (s)$ we omit
them in the following.

\begin{prop}\label{MUINTS}  We have:

$$ \left\{ \begin{array}{lll} \mu_0(s)& = &
 \frac{\Gamma(\frac{1}{2}) \Gamma (s - \frac{1}{2})}{\Gamma(s)} \;\;\; (\Re e\;s > \frac{1}{2})\\ &
 &\\
\mu_{ir_k}^c (s) & = & \frac{\Gamma(\frac{1}{2}) \Gamma (s
-\frac{1}{4} + \frac{2i r_k}{4}) \Gamma (s -\frac{1}{4} -
\frac{2i r_k}{4})}{\Gamma(s)^2} \;\;\; (\Re e\;s > 0) \;\;\\ & & \\
\mu^d_m (s) & = & \frac{(-i)^{m/2} \pi 2^{2s + 2 - m/2} \Gamma(-2s
+ \frac{m}{2})}{ - (2s + 1 - \frac{m}{2}) \Gamma(-s) \Gamma(-s +
\frac{m}{2})}\;\;\; (\Re (2s - \frac{m}{2}) < -1).
\end{array}
 \right. $$

 \end{prop}

The proof is given in  \cite{Z} (see pages 50-52).

\section{ \label{SECOND} Proof of Theorem \ref{maincont}.}

The key objects in the proof of Theorem \ref{maincont} are the
closely related integrals \begin{equation} \label{IMPINTS} \left\{
\begin{array}{lll} I_r(\sigma) &= & \int_\R (1 +
u^2)^{-(\frac{1}{2}+ir)} PS_{ir}(( \sigma \chi)^u) du, \\&& \\
I^{\Gamma}_r (\sigma) &=& \int_\R (1 + u^2)^{-(\frac{1}{2}+ir)}
PS_{ir}^{\Gamma}( \sigma^u) du
\end{array} \right.
\end{equation}  where $PS_{ir_j} \in \dcal'(\D), PS_{ir_j}^{\Gamma} \in \dcal'(\Gamma \backslash G)$
are defined in Definition \ref{PSDEF} and where  $f^u(g) = f(g
n_u)$. Note that $PS_{ir}^{\Gamma}( \sigma^u) =
I_{PS_{ir}}(\sigma)(u)$ in the notation of \S \ref{INTOPS}. It
takes some work to prove that each integral is well-defined. In
Lemma \ref{CONVEQ} it is proved that the two integrals are
well-defined and  equal for $\sigma \in C^{\infty}(\Gamma
\backslash G)$.

Theorem \ref{mainintro} equates the Wigner distribution with the
distribution $\sigma \to \int_{S \D} (L_{r_j}\chi  \sigma)
dPS_{j}$ for $\sigma \in C^{\infty}(\Gamma \backslash G)$. In
Proposition \ref{CONTINUITY} we show that this functional also
equals $I_r(\sigma) = I^{\Gamma}_r (\sigma)$.  The explicit
formulae for the Wigner distributions in terms of the
Patterson-Sullivan distributions follow from the identification
with $I^{\Gamma}_r(\sigma)$, which can be explicitly evaluated
using the results of \S \ref{INTOPS}.

\subsection{Convergence and equality of the integrals}

In the following, we recall that  $\Re ir_k = 0$ in the unitary
principal series but is positive in the complementary series.

\begin{prop} \label{CONVEQ} We have:

\begin{enumerate}

\item Integral $I^{\Gamma}_r(\sigma)$ converges for all $\sigma
\in C^{\infty}(\Gamma \backslash G)$ if $\pcal_{ir}$ is in the
unitary principal series.

\item $I_r(\sigma) = I^{\Gamma}_r(\sigma)$ for all  $\sigma \in
C^{\infty}(\Gamma \backslash G)$.

\end{enumerate}

\end{prop}

\subsubsection{Proof of (1)}

We  give a representation theoretic proof that
\begin{equation} \label{BASICINT} \int_\R (1 +
u^2)^{-(\frac{1}{2}+ir)} I_{PS_{ir}}^{\Gamma}(\sigma)(u) du
\end{equation}  converges absolutely.
 We
make no attempt at a sharp estimate but only one sufficient for
the purposes of this paper.

\begin{lem} \label{HALFDECAY} Let $PS_{ir}$ be the
Patterson-Sullivan distribution corresponding to $\phi_r$. Then:
$$\left\{ \begin{array}{ll} (i) & I_{PS_{ir_j}}(\phi_{ir_k} )(u) \leq   C\; (1 + |r_k|)^4 (1+|r_j|)^K (1 + |u|)^{-1/2 + \Re (i r_j)}; \\ \\
(ii) & I_{PS_{ir_j}}(\psi_m) \leq  C \; (1 + |m|)^4 (1+|r_j|)^K (1
+ |u|)^{- m} ,
\end{array} \right.$$
where $K$ is the same as in (\ref{normofPS}).
\end{lem}

\begin{proof}

\noindent(i)  By Propositions \ref{ISIGMA} and \ref{TRX},
\begin{equation} \label{ISPRJPHIK} \begin{array}{lll}  I_{PS_{ir_j}}(\phi_{ir_k} )(u) & =
&\int_{\Gamma \backslash G}
 \phi_{ir_k}(g n_u) PS_{ir_j}(dg) \\ && \\ & = & \big(\int_{\Gamma \backslash G}
\phi_{ir_k} PS_{ir_j}(dg) \big)F_{ir_k,
0}\left(\frac{u-i}{-2i}\right)
 \end{array} \end{equation}

  By (\ref{normofPS}),
there exists $K$ so that
$$|\langle \phi_{ir_k}, PS_{ir_j} \rangle| \leq C(1+|r_j|)^K (1 + |r_k|)^4.$$
Here, we used a crude estimate  $||\phi_{ir_k}||_{C^2} \leq C (1 +
|r_k|)^4$ (in fact, $r_k^3/\log r_k$ is true, but it is not
necessary for our argument). We combine with the estimates in
Proposition \ref{MURKEST}  (cf. \ref{HYPER}) on the hypergeometric
factor to obtain the estimate stated in (i).

\medskip

\noindent(ii) We now have
\begin{equation} \label{ISPRJPSIM} \begin{array}{lll}  I_{PS_{ir_j}}(\psi_{m} )(u) & =
&\int_{\Gamma \backslash G}
 \psi_{m}(g n_u) PS_{ir_j}(dg) \\ && \\ & = & \big(\int_{\Gamma \backslash G}
\psi_{m} PS_{ir_j}(dg) \big) (u+i)^{-m/2}.
 \end{array} \end{equation}
To complete the proof we note that  $|(u+i)^{-m/2}| \leq C (1 +
|u|)^{-1/2}$ and that (by (\ref{normofPS})),
$$|\langle \psi_m , PS_{ir_j} \rangle|\leq  C(1+|r_j|)^K  (1 +
|m|)^4.$$
\end{proof}

Given a co-compact discrete group $\Gamma \subset SL(2, \R)$ we
denote by $\tau_0 = \Im (i r_0)$ the Casimir parameter
corresponding to the lowest non-zero eigenvalue of $\Delta$, i.e.
the complementary series representation closest to the trivial
representation.

\begin{lem} \label{HALFDECAYa}  If  $\sigma \in C^{\infty}(\Gamma \backslash G)$  has no
component in the trivial  representation, we have:
$$I_{PS_{ir}}(\sigma)(u) \leq  C (1 + |r|)^K  (1 + |u|)^{-1/2 + \tau_0}. $$

\end{lem}

\begin{proof}

 Since $PS_{ir}$ is geodesic flow and time reversal
invariant, we may write by Proposition \ref{EXPRESSION},
\begin{equation} \label{BIGSUM} I_{PS_{ir}} (\sigma) (u) =
\sum_{r_j} \frac{ \langle \sigma, \Xi_{ir_j} \rangle}{\langle
\phi_{ir_j}, \Xi_{ir_j} \rangle} I_{PS_{ir}}(\phi_{ir_j} )(u)+
 \sum_{m, \pm} \frac{ \langle \sigma, \Xi_m^{\pm} \rangle}{\langle \psi_m, \Xi_m^{\pm} \rangle}
 I_{PS_{ir}} ( \psi_m)(u) .   \end{equation}  It
follows by Lemma \ref{HALFDECAY} that \begin{equation}
\label{SUMM}\begin{array}{lll} |I_{PS_{ir}} (\sigma)(u)| &\leq & C
(1 + |r|)^K  \times \{\sum_{r_j} (1 + |r_j|)^{4}  \left| \frac{
\langle \sigma, \Xi_{ir_j} \rangle}{\langle \phi_{ir_j},
\Xi_{ir_j} \rangle} \right|\; (1 +
|u|)^{-1/2 +  \Re ir} \\ && \\
&+&
 \sum_m (1 + |m|)^{4} |\frac{ \langle \sigma, \Xi_m^{\pm} \rangle}{\langle \psi_m, \Xi_m^{\pm} \rangle}|\;\; (1 + |u|)^{-m}\} \end{array}.
 \end{equation}

 By Proposition \ref{XIEST},
 $$ \left| \frac{ \langle \sigma, \Xi_{ir_j}
\rangle}{\langle \phi_{ir_j}, \Xi_{ir_j} \rangle} \right| \leq
||X_- \pcal_{ir_j}(\sigma)||_{L^2}, \;\; |\frac{ \langle \sigma,
\Xi_m^{\pm} \rangle}{\langle \psi_m, \Xi_m^{\pm} \rangle}| \leq
||X_- \dcal_m^{\pm} )(\sigma)||_{L^2}. $$ It follows that  for any
$M
> 0$ there exists a constant $C_M$ so that
\begin{equation} \label{ABSCONVLEM} \left| \frac{ \langle \sigma, \Xi_{ir_j} \rangle}{\langle
\phi_{ir_j}, \Xi_{ir_j} \rangle} \right| \leq C_M (1 +
|r_j|)^{-M}, \;\;\; |\frac{ \langle \sigma, \Xi_m^{\pm}
\rangle}{\langle \psi_m, \Xi_m^{\pm} \rangle}| \leq C_M  (1 +
|m|)^{-M}.
\end{equation} Indeed,
$$\pcal_{ir}(\sigma) = \sum_{m \in \Z} \sigma_{ir, m} \phi_{ir,
m}, \;\; \mbox{with} \;\; |\sigma_{ir, m}| \leq C_M (1 + |r_j| +
|m|)^{-M},  $$  hence
$$||X_- \pcal_{ir_j}(\sigma)||_{L^2}  \leq C_M \sum_m (1 +
|r_j| + |m|)^{-M} (1 + |r_j| + |m|),$$ where we bound $||X_-
\phi_{ir, m}||_{L^2} \leq C (1 + |r_j| + |m|)$.

Similarly,
$$(\dcal_m^{\pm}  )(\sigma) = \sum_{n = 0}^{\infty} \sigma_{m, m + 2n} \psi_{m,
m + 2n}, \;\; \mbox{with} \;\; |\sigma_{m, m + 2n}| \leq C_M (1 +
|m| + |n|)^{-M},  $$ hence
$$||X_- \dcal_m^{\pm}(\sigma)||_{L^2}  \leq C_M \sum_m (1 +
 |m| + |n|)^{-M} (1  + |m| + |n|).$$

By (\ref{ABSCONVLEM}) and Lemma (\ref{HALFDECAY}),  the sum
(\ref{BIGSUM}) converges absolutely and the decay estimates in $u$
sum up to the stated rate.
\end{proof}

Completion of proof of Proposition \ref{CONVEQ} (1): It follows
from Lemma \ref{HALFDECAYa} that

\begin{equation} \label{INTINEQ} \begin{array}{lll} |I^{\Gamma}_r(\sigma)| & \leq & \int_\R (1 +
u^2)^{-(\frac{1}{2}+ \Re ir)}| |PS_{ir}^{\Gamma}( \sigma^u) |du  \\ && \\
& \leq &  C (1 + |r|)^K  \int_\R |(1 + u^2)^{-(\frac{1}{2}+ \Re
ir)}| (1 + |u|)^{-1/2 + \tau_0} du. \end{array} \end{equation}
Since $\pcal_{ir}$ is in the unitary principal series, $\Re ir =
0$ and so  $|(1 + u^2)^{-(\frac{1}{2}+ir)}| = (1 +
u^2)^{-\frac{1}{2}}$ and since $-1/2 + \tau_0 > 0$ it follows that
the last integral in (\ref{INTINEQ}) converges absolutely.

We now move on to the assertion (2) of Proposition \ref{CONVEQ}.

\subsubsection{Proof of (2)}

 By Proposition
\ref{inv-principle}, we have
\begin{equation} \label{MOVECHI} \int_{ G} \sigma (gn_u) \chi(g n_u)
dPS_{ir}(g) = \int_{ G} \sigma (gn_u) \chi(g ) dPS_{ir}(g).
\end{equation}
 Indeed,   $\chi(g)$ and o$\chi^u(g):
= \chi(g n_u)$ are both smooth fundamental cutoffs, so both sides
equal $\int_{\Gamma \backslash G} \sigma dPS_{ir}^{\Gamma}. $
Integrating against $\int_\R (1 + u^2)^{-(\frac{1}{2}+ir)} $
completes the proof.

\subsection{Continuity of $PS_{ir}$}

As mentioned above, the Wigner distribution equals the functional
$\sigma \to \int_{S \D} (L_{r_j}\chi \sigma) dPS_{ir_j}$. To prove
that this also equals $I_r(\sigma) = I^{\Gamma}_r(\sigma)$ we need
the following continuity result for the functional $PS_{ir}$.

\begin{lem}\label{CONTINUITY}  $PS_{ir} \in \dcal'(S \D)$ has the following continuity property,
$$PS_{ir} (L_r(\chi \sigma))  = \int_\R (1 +
u^2)^{-(\frac{1}{2}+ir)} PS_{ir}( \sigma  \chi)^u) du,$$ where
$f^u(g) = f(g n_u)$.
\end{lem}

\begin{proof}
 By Definition \ref{PSDEF},
$$\begin{array}{lll} PS_{ir}( (\sigma \chi)^u) & = & ps_{ir}(\rcal(
\sigma \chi)^u))\\ && \\ & = & \int_{B \times B'} \{\int_{\R_t}
(\chi
 \sigma) (g(b, b') a_t n_u) dt \} \frac{dT_{ir}(b) \otimes
dT_{ir}(b')}{|b - b'|^{1 + 2 i r}}
\end{array}
$$
We first note that for all $u$, $\rcal(\chi \sigma)^u \in
C_c^{\infty}(B \times B' \setminus \Delta)$ since $(\chi \sigma)^u
\in C_c^{\infty}(S \D)$. It follows that $ps_{ir}(\rcal( \sigma
\chi)^u)$ is well-defined and smooth in $u$.

The continuity statement is equivalent to
\begin{equation} \label{TOPROVE} \begin{array}{l}
 ps_{ir} (\rcal (L_r(\chi \sigma))  = \int_\R (1 +
u^2)^{-(\frac{1}{2}+ir)} ps_{ir}(\rcal( \sigma \chi)^u) (b, b' ))
du,
\end{array} \end{equation}
or equivalently
\begin{equation} \label{TOPROVE2}
 ps_{ir} ( \int_\R (1 + u^2)^{-(\frac{1}{2}+ir)} \rcal( \sigma
 \chi)^u(b, b')) du = \int_\R (1 +
u^2)^{-(\frac{1}{2}+ir)} ps_{ir}(\rcal( \sigma \chi)^u)(b, b') du.
\end{equation}
We must again check  that both sides of (\ref{TOPROVE})  are
well-defined. Clearly, $\rcal (L_r(\chi \sigma))$ is well-defined
because $\chi$ has compact support. The problem is to prove that
the left-hand side is well-defined, since
 that $ps_{ir}$ is only known to be a bounded linear
functional on $|b - b'|^{ 1 + 2 i r} C^2(B \times B' \setminus
\Delta)$ (cf. Proposition \ref{BASICBB}). We therefore have   to
verify that
$$\int_\R (1 + u^2)^{-(\frac{1}{2}+ir)} \rcal( \sigma
 \chi)^u(b, b')) du \in |b - b'|^{ 1 + 2 i r} C^2(B \times
B' \setminus \Delta). $$

 By  Lemma \ref{LEMMA} and (\ref{LCAL}), we have
 \begin{equation} \int_\R (1 +
u^2)^{-(\frac{1}{2}+ir)} \rcal( \sigma
 \chi)^u(b, b')) du = |b - b'|^{(1 + 2 i r)} \int_{D} (\chi
 \sigma) (z, b) e^{(\frac12+ir) \langle z, b
\rangle} e^{(\frac12+ir) \langle z, b' \rangle } dVol(z),
\end{equation}
and therefore the condition to be satisfied is that
\begin{equation}\int_{D} (\chi
 \sigma) (z, b) e^{(\frac12+ir) \langle z, b
\rangle} e^{(\frac12+ir) \langle z, b' \rangle } dVol(z) \in
C^{\infty}(B \times B).
\end{equation}
This is clear due to the compact support of $\chi \sigma$ in $z$,
which is independent of $(b, b')$.

We may then rewrite (\ref{TOPROVE2}) as:
\begin{equation} \label{TRCONT}\begin{array}{l}  \langle \big(\int_{D} (\chi
 \sigma) (z, b) e^{(\frac12+ir) \langle z, b
\rangle} e^{(\frac12+ir) \langle z, b' \rangle } dVol(z) \big), dT_{ir} \otimes dT_{ir}  \rangle  \\
\\= \int_{D} \chi
\langle  \big(\sigma (z, b) e^{(\frac12+ir) \langle z, b \rangle}
e^{(\frac12+ir) \langle z, b' \rangle }\big),  dT_{ir} \otimes
dT_{ir} \rangle dVol(z),\end{array}
\end{equation}
where $\langle, \rangle$ denotes the pairing of test functions and
distributions on $B \times B$. I.e. we need to check that we can
pass $dT_{ir} \otimes dT_{ir}$ under the $dVol(z)$ integral sign.

By Otal's regularity theorem (see Theorem \ref{O}), $dT_{ir} =
F'_r(b) db$ where  $F_{ir}$ is a continuous $2\pi$ periodic
function in the sense that $F_{ir}(\theta + 2 \pi) -
F_{ir}(\theta) = C_r$. Integration by parts then gives
$$ \langle g, T_{ir} \rangle = \int_B g(b)T_{ir}(db)= -\int_B g'(b) F_{ir}(b)db
+g(0)(F_{ir}(2\pi)-F_{ir}(0)).
$$ Applying this in each of the $(b, b')$ variables to  the pairing on $B \times B$ in
(\ref{TRCONT}) produces four terms of which three involve the
boundary term $(F_{ir}(2\pi)-F_{ir}(0))$ and the fourth is
$$\int_{B \times B} \{ \int_{D} \frac{\partial^2}{\partial b \times \partial b'} \left((\chi
 \sigma) (z, b) e^{(\frac12+ir) \langle z, b
\rangle} e^{(\frac12+ir) \langle z, b' \rangle } dVol(z) \right)\}
F_{ir}(b) F_{ir}(b') db \otimes db'.$$ By applying Fubini's
theorem to the fourth term, we obtain:
\begin{equation} \label{TRCONT2}\begin{array}{l}  \int_{B \times B} \{ \int_{D} \frac{\partial^2}{\partial b \times \partial b'} \left((\chi
 \sigma) (z, b) e^{(\frac12+ir) \langle z, b
\rangle} e^{(\frac12+ir) \langle z, b' \rangle } dVol(z) \right)\}  F_{ir}(b) F_{ir}(b') db \otimes db'  \\
\\= \int_{D}  \chi \{\int_{B \times B} \frac{\partial^2}{\partial b \times \partial b'}
  \left(\sigma (z, b) e^{(\frac12+ir) \langle z, b \rangle}
e^{(\frac12+ir) \langle z, b' \rangle } \right) T_{ir} (b)
T_{ir}(b') db db'\} dVol(z),. \end{array}
\end{equation}
Fubini's theorem applies in a similar way to the other terms.
 We then  transfer the $b$ derivatives back to $dT_{ir}$ and
obtain (\ref{TRCONT}).

\end{proof}

As a corollary of Proposition \ref{CONVEQ}, we obtain the
following explicit formula:

\begin{cor} We have:
\begin{equation} \label{BASICINT2} \begin{array}{lll} \int_\R (1 +
u^2)^{-(\frac{1}{2}+ir)} I_{PS_{ir}}^{\Gamma}(\sigma)(u) du & = &
\sum_{r_j} \frac{ \langle \sigma, \Xi_{ir_j} \rangle}{\langle
\phi_{ir_j}, \Xi_{ir_j} \rangle} \big(\int_{\Gamma \backslash G}
\phi_{ir_j} PS_{ir}(dg) \big) \mu^c_{ ir_j}(\frac{1}{2} + ir)
\\ && \\ && +
 \sum_{m, \pm} \frac{ \langle \sigma, \Xi_m^{\pm} \rangle}{\langle \psi_m, \Xi_m^{\pm} \rangle}
\big(\int_{\Gamma \backslash G} \psi_{m} PS_{ir}(dg) \big)
\mu^d_{m}(\frac{1}{2} + ir).
\end{array} \end{equation}
All integrals and series converge absolutely.
\end{cor}

\begin{proof}

In fact, by  Lemma \ref{HALFDECAYa}  we may interchange the order
of summation in (\ref{BIGSUM}) and integration in
(\ref{BASICINT}). Using (\ref{ISPRJPHIK}) and (\ref{BIGDISPLAY}),
we have
$$ \int_\R (1 + u^2)^{-(\frac{1}{2}+ir)}
I_{PS_{ir}}^{\Gamma}(\phi_{ir_j})(u) du = \big(\int_{\Gamma
\backslash G} \phi_{ir_j} PS_{ir}(dg) \big) \mu^c_{
ir_j}(\frac{1}{2} + ir),
$$ and thus obtain the first series of the stated formula. Using
(\ref{ISPRJPSIM}) and (\ref{BIGDISPLAY}), we have
$$ \int_\R (1 + u^2)^{-(\frac{1}{2}+ir)}
I_{PS_{ir}}^{\Gamma}(\psi_m)(u) du = \big(\int_{\Gamma \backslash
G} \psi_{m} PS_{ir}(dg) \big) \mu^d_{ m}(\frac{1}{2} + ir),
$$ and obtain the second series.
\end{proof}

\subsection{Completion of Proof of Theorem \ref{main2}}

To complete the proof it suffices to explicitly evaluate
$I^{\Gamma}_r(\sigma)$ for the generating automorphic forms.

\begin{lem} \label{LR2} In the special cases when $\sigma=\phi_{ir_k}, X_+\phi_{ir_k} $ or
$\psi_m$, we have the explicit formulae:

\begin{enumerate}
\item In the case $\sigma = \phi_{ir_k}$, $   \langle
Op(\phi_{ir_k})\phi_{ir_j}, \phi_{ir_j} \rangle =
 \mu_{ir_k}^c(\frac{1}{2} + i r_j)  \big(\int
\phi_{ir_k} dPS^\Gamma_{ir_j}(g)\big)$.

\item For $\sigma=X_+ \phi_{ir_k}$, $  \langle Op(X_+ \phi_{ir_k})
\phi_{ir_j}, \phi_{ir_j}\rangle= 0 $ for all $j$.

\item For $\sigma = \psi_m$, $ \langle Op(\psi_m)\phi_{ir_j},
\phi_{ir_j} \rangle = \mu_m^d (\frac{1}{2} + i r_j)  \big(\int
\psi_m dPS^\Gamma_{ir_j}(g)\big).$

\end{enumerate}

Here, the expressions $ \mu_{ir_k}^c(\frac{1}{2} + i r_j), \mu_m^d
(\frac{1}{2} + i r_j) $ are defined in (\ref{BIGDISPLAY}).

\end{lem}

\begin{proof}

All statements follow from the combination of Proposition
\ref{ISIGMA} and Proposition \ref{MUINTS}. The case
$\sigma=X_+\phi_{ir_k}$ follows from Proposition \ref{TRX} and the
fact that the Patterson-Sullivan distributions are invariant under
time-reversal (cf. Proposition \ref{INVARIANT}).

\end{proof}

We note that these explicit formula give a new proof of Theorem
\ref{mainintro}:

\begin{cor} \label{RATE}
When $\sigma$ is a joint $(\Omega, W)$-eigenfunction, we find
again that $\langle Op(\sigma)\phi_{ir}, \phi_{ir}\rangle$ is
asymptotically the same as $r^{-1/2}\int \sigma dPS_{ir}$.
\end{cor}

\begin{proof}
By definition, $F_{\tau, m}(1/2)=1$ whereas $G_{\tau, m}(1/2)=0,
G^\prime(1/2)=-2i$. The stationary phase method then shows that
$$\int c(u)^{-(1+2ir)}
F_{\tau,m}\left(\frac{u-i}{-2i}\right)du\sim r^{-1/2}$$ whereas
$$\int c(u)^{-(1+2ir)}
G_{\tau,m}\left(\frac{u-i}{-2i}\right)du \sim r^{-3/2}.$$ Here, we
use the estimates in (\ref{HALFDECAY}), which are valid in all
weights.

\end{proof}

\section{\label{THERMO} Dynamical zeta-functions: Thermodynamic formalism}

In this part, we relate the   Patterson-Sullivan distributions  to
the  Ruelle resonances introduced in \cite{Ru87}.

 Ruelle's aim was to study the Fourier transform of the correlation function,
$$\rho_{F,G}(t)=\int F( x)G(g^tx)d\omega(x)-\int Fd\omega\int Gd\omega,$$
($t\geq 0$), in the very general context of an Axiom A flow
$(g^t)$ (for instance, when $\omega$ is the measure of maximal entropy). He showed, for smooth enough functions $F$ and $G$, that
the Fourier transform $\hat\rho_{F, G}$ has a meromorphic extension
to a half-plane of the form $\{\Re es> h-\eps\}$, strictly beyond
its half-plane of absolute convergence $\{\Re e\;s> h\}$ (where $h$
represents, in a general context, the topological entropy). He used the
so-called
``thermodynamic formalism'' and showed that the poles of $\hat\rho(s)$
occured precisely for certain values $s$, linked with the existence of certain distributions obeying specific transformation rules.

 In  the case of the
geodesic flow on a compact surface of constant curvature $-1$, and
for  $C^1$ functions $F, G$ on $\Gamma\backslash G$, the Fourier
transform $\hat\rho$ is an analytic function in the half-plane
$\{s, \Re e s>1\}$ and has a meromorphic extension to $\{ \Re e
s>0\}$ with poles  at values of $s=\frac12+ir$ for which  there exists a
distribution $e_{ir}$ on $S\X$ satisfying:
\begin{itemize}

\item  $g^t.e_{ir}=e^{-(\frac12+ir)t}e_{ir}$

\item  $e_{ir}$ is invariant under the stable horocycle flow.

\end{itemize}

In the case of constant negative curvature,  it follows
that $e_{ir}$ is given by: \begin{equation} \label{ES}
e_{ir}(F)=\int F(z,b)e^{(\frac12+ir)\langle
z,b\rangle}T_{ir}(db)dVol(z)=\int_{\X}Op(F)\phi_{ir}(z)dVol(z)
\end{equation}  where $T_{ir}$ is the boundary values of an eigenfunction
$\phi_{ir}$ of $\Lap$ of eigenvalue $\frac14+r^2$ (see equation \eqref{EJ},
and \cite{Z2} ) Hence the poles of
$\hat\rho$, i.e. the Ruelle resonances, occur at $s_n=1/2+ir_n$.
If the eigenvalue is simple, the residue of $\hat\rho_{a,b}$ at
$s_n$ is given, up to multiplicative normalization, by
\begin{multline*} R_{r_n}(F,G)=\left(\int F(z,b)e^{(1/2+ir_n)\langle z,b\rangle}T_{ir_n}(db)dVol(z)\right)\left(\int G\circ\iota(z,b)e^{(1/2+ir_n)\langle z,b\rangle}T_{ir_n}(db)dVol(z)\right)\\
=\left(\int_{\X}Op(F)\phi_{ir_n}dVol(z)\right)\left(\int_{\X}Op(G\circ\iota)\phi_{ir_n}dVol(z)\right)
=e_{ir_n}(F)e_{ir_n}(G\circ\iota),
\end{multline*}
where $\iota$ denotes time reversal. To see this, we
observe that the residue  $R_{r_n}(F,G)$ is bilinear in $F$ and
$G$, and its definition implies that it satisfies the identities
$$R_{r_n}(F\circ g^t, G)=R_{r_n}(F,G\circ g^{-t})=e^{-(1/2+ir_n)t}R_{r_n}(F,G),$$
 and
$$R_{r_n}(F\circ h_+^u, G)=R_{r_n}(F, G\circ h_-^{u})=R_{r_n}(F,G),$$
where $h_+
^u$ denotes the stable horocyclic flow and $h^u_-$ the
unstable one. It follows that it must equal
$e_{ir_n}(F)e_{ir_n}(G\circ\iota)$ is the eigenvalue $1/4+r_n^2$ is simple.
 If the eigenvalue is not simple, the
expression becomes more complicated, as one has to form a linear
combination of the functionals associated to the various
eigenfunctions.

We now prove prove Theorem \ref{main2} for $\lcal_2$. In this
section,
 the Patterson-Sullivan distributions are  normalized so
that $\int \1 dPS
_{ir_n}=1$.

\begin{rem} The residue is obviously not the same as in Ruelle's paper: in our case the
residue must define a geodesic flow invariant functional, whereas in Ruelle's paper it has to be invariant under the horocycle flow.\end{rem}

\subsection{Markovian coding of the boundary.}
The proof given here relates the function $\lcal_2$ to the determinant
of certain operators, called transfer operators.
To define them, we need to recall from \cite{Se}
the construction of Markov sections, using the Bowen-Series coding
of the action of $\Gamma$ on the boundary $B$. Series used this
construction to study Poincar\'e series. We apply it to the
somewhat different objects $\lcal_2$. For this application, we
need some further discussion of Markov coding which we could not
find in the literature.

If we want to study the action of $\Gamma$ on the boundary, and
the existence of conformally invariant distributions -- by this,
we mean the property \ref{CONFORMAL} -- it is enough to consider a
set of generators of $\Gamma$. In fact, it is even enough to work
with a single, well chosen transformation of the boundary: roughly
speaking, this transformation $F^{(r)}$ will be defined by cutting
the boundary $B$ into a finite number of closed intervals $J_i$,
and will act on each $J_i$ by the action of one of the chosen
generators of $\Gamma$.

We will require the map
$F^{(r)}:J=\sqcup J_i\To J$ to have the following properties:

(i) $F^{(r)}$ is analytic, expanding (or at least, some power of
$F^{(r)}$ is expanding).

(ii) The $J_i$s form a Markov partition for $F^{(r)}$. This means
that $F^{(r)}$ sends the boundary of $J$ to itself.

(iii) The natural map $J=\sqcup J_i\To B$ gives a bijection
between periodic points of $F^{(r)}$ and points at infinity of
closed geodesics, except for the closed geodesics ending on the
boundary of an interval $J_i$, that have exactly two preimages. If
${F^{(r)}}^n x=x$, and $\gamma$ is the closed geodesic
corresponding to $x$, then $|({F^{(r)}}^n)^\prime
x|=e^{L_\gamma}$.

We recall briefly the construction of $F^{(r)}$ proposed by
Series \cite{Se}, when $\Gamma$ is cocompact: she chooses a
symmetric generating set for $\Gamma$, called $\Gamma_0$. She then
defines a notion of ``admissible representation'' of an element
$g\in\Gamma$ as a word $g=g_1...g_n$ with $g_j\in\Gamma_0$, such
that

-- an admissible word is a shortest representation of $g$ in the
alphabet $\Gamma_0$.

-- every $g\in\Gamma_0$ has a unique admissible representation.

Without going into details, admissible words are shortest word-representations;
and besides, whenever there is a choice of several such representations,
one chooses the one that ``turns the furthest possible to the right''
in the Cayley graph of $\Gamma$ with respect to $\Gamma_0$ (seen
as a subset of the hyperbolic plane).

Let us denote $\Sigma^{(r)}_f$ the set of finite admissible words;
the notation is borrowed from \cite{Se} but we are adding an $r$
to specify that we are choosing representations that turn the most
possible right in the Cayley graph -- the same convention Series used in her paper.
Replacing ``right" by ``left" one would obtain another notion of
admissible words, and we denote $\Sigma_f^{(l)}$ the set of
left-admissible words. Note that
$\Sigma_f^{(l)}={\Sigma_f^{(r)}}^{-1}$. Now define $\Sigma^{(r)}$,
the set of infinite right-admissible words, as
$$\Sigma^{(r)}=\{(g_j)\in\Gamma_0^{\Z_+}, g_{j}...g_{j+k}\in\Sigma_f^{(r)},
\forall j,k\geq 0\}.$$ Series shows in \cite{Se} that the map
\begin{equation}\Sigma^{(r)}_f \To  \HH,\;\;\; g_1...g_n \mapsto g_1...g_n.0
\end{equation} can be extended to a continuous map $j^{(r)}:\Sigma^{(r)}\To B$. She denotes $I^{(r)}(g_j)$
the set of points in $B$ that have a representation by a sequence
in $\Sigma^{(r)}$ starting with the generator $g_j$. The boundary
$B$ is thus divided into a finite number of closed intervals,
with disjoint interiors. One can define a
map $F^{(r)}$ that acts on $\Sigma^{(r)}$ by deleting the first
symbol and shifting the sequence to the left; seen as a map on
$B$, it acts as $g_j^{-1}$ on each interval $I^{(r)}(g_j)$.
Actually, the map $F^{(r)}$ is defined on $I^{(r)}:=\sqcup
I^{(r)}(g_j)$; when working on $B$ one should always remember that
its definition is ambivalent on boundary points. The partition
$B=\cup I^{(r)}(g_j)$ is not exactly a Markov partition for the
action of $F^{(r)}$: there is no reason that boundary points
should be sent to boundary points. But, by construction, the
images of these boundary points under iteration of ${F^{(r)}}$ form a finite set. Cutting
the intervals $I^{(r)}(g_j)$ at these points, one can refine the
partition $B=\cup I^{(r)}(g_j)$ into a new finite partition
$B=\cup J_j$ that is now Markov. This way we obtain a
tranformation $F^{(r)}$ satisfying all the conditions (i), (ii),
(iii). An element of $B$ may be coded either by a word in
$\Sigma^{(r)}$, as we have already seen, or by an element of the
subshift of finite type
$$\Lambda^{(r)}=\{(i_k)_{k\geq 0},F^{(r)}(J_{i_k})\cap int(J_{i_{k+1}})\not=\emptyset
\mbox { for all } k\geq 0\}.$$ Both codings are bijective except
on a countable set (in fact the coding map is at most $2$ to $1$).

To make the link with the geodesic flow, we now  extend the
expanding transformation $F^{(r)}$ to an invertible transformation
$F$ of a subset of $B\times B$; in terms of symbolic dynamics, we
want to work with two-sided subshifts. We consider
$$\Sigma^{(l)}=\{(g_j)\in\Gamma_0^{\Z_+^*},  g_{j}...g_{j+k}\in\Sigma_f^{(l)},
\forall j, k >0\}$$ and
$$\Sigma_{(l)}=\{(g_j)\in\Gamma_0^{\Z_-^*},  g_{j-k}...g_{j}\in\Sigma_f^{(r)}\forall j, k <0\}.$$
Formally, elements of $\Sigma_{(l)}$ are inverses of elements of
$\Sigma^{(l)}$. By the same consideration as before, we have a coding map
$j^{(l)}:\Sigma^{(l)}\To B$ or equivalently
$j_{(l)}:\Sigma_{(l)}\To B$; we denote $I_{(l)}(g_j)\subset B$ the
interval formed by points whose coding in $\Sigma_{(l)}$ ends with
$g_j$. This gives a partition $B=\cup I_{(l)}(g_j)$ and a map
$F_{(l)}$ on $B$, that corresponds to the shift to the right on
$\Sigma_{(l)}$. If we want we can refine the partition $B=\cup
I_{(l)}(g_j)$ into a Markov partition $B=\cup K_j$ and code the
dynamics by a one-sided subshift of finite type $\Lambda_{(l)}$.

Let us now introduce the two-sided subshift $\Sigma$,
$$\Sigma=\{(g_j)\in\Gamma_0^{\Z}, g_{j}...g_{j+k}\in\Sigma_f^{(r)},
\forall j,k\};$$ $\Sigma$ is naturally isomorphic to a subset of
$\Sigma_{(l)}\times\Sigma^{(r)}$, and thus there is a coding map
from $\Sigma$ to a subset $X\subset B\times B$:
\begin{eqnarray*}j:\Sigma&\To& X,\\
j(\sigma_{(l)}, \sigma^{(r)})&=&\left(
j_{(l)}(\sigma_{(l)}),j^{(r)}(\sigma^{(r)})\right)
\end{eqnarray*}

The shift to the left on $\Sigma$ gives an invertible transformation $F$
on $X$; note, as above, that $F$ is actually well defined on a
subset of $\sqcup I_{(l)}(g_j)\times \sqcup I^{(r)}(g_k)$ and is
defined ambivalently at certain points of $X$. If $y$ is
in $I^{(r)}(g_j)$ then $F(x,y)=(g_j^{-1}x, g_j^{-1}y)$, or in
other words $F(x,y)=({{G_{(l)}}}_j x, F^{(r)} y )$, where
${{G_{(l)}}}_j$ is the inverse branch of $F_{(l)}$ taking values
in $I_{(l)}(g_j)$.

The partition of $X$ into $ \sqcup I_{(l)}(g_j)\times\sqcup
I^{(r)}(g_k)$ is not a Markov partition for the action of $F$, but
$X=\cup  (K_j\times J_i)$ is Markov. If we want to work with a
subshift of finite type, we should consider the set
$$\Lambda=\{(m_k, n_k)_{k\in\Z} / F(K_{m_k}\times J_{n_k})\cap int(K_{m_{k+1}}
\times J_{n_{k+1}})\not=\emptyset \mbox{ for all } k\in\Z\}.$$

We can identify $X$ with a subset of the unit tangent bundle
$\Gamma\backslash G$ that forms a transversal section for the action of
the geodesic flow: we observe that, for each
$(x,y)$ in $X$, the geodesic $\gamma_{x,y}$  contains a unique
vector, denoted $v_{x,y}$, which is in the stable manifold of a
vector based at $0$. To recover the action of the geodesic flow on
the whole tangent bundle, one needs to add a time parameter
measuring the time it takes to go from $(x,y)$ to $F(x,y)$.
Because $\gamma_{x,y}$ and $\gamma_{F(x,y)}$ represent the same
geodesic in the quotient $\Gamma\backslash G$, there exists a
unique $\tau(x,y) \in\R$ such that $g^\tau v_{x,y}=v_{F(x,y)}$.
More precisely, the function $\tau$ is defined without any
ambiguity on $\sqcup \overset{\circ}{ K_j}\times
\overset{\circ}J_i$ and can then be extended to a continuous
function on $\sqcup  K_j\times J_i$. By construction,
the
function $\tau$ is actually locally constant on stable manifolds;
i.e., it depends only on the variable $y$. It is analytic on each
rectangle $K_j\times J_i$. We see $\tau$ as a return time from the
section $X$ to itself; note however that $\tau$ may change sign:
we are not exactly in the usual situation of a ``first return time".
Nevertheless, when $y$ is a
periodic point of period $n$ of $F^{(r)}$, and $\gamma$ is the
corresponding closed geodesic, we have $\sum_{k=0}^{n-1}
\tau({F^{(r)}}^k y)= L_\gamma>0$.

We have a surjection, almost one-to-one, from
$$X^\tau:=\{((x,y), t)\in X\times \R, t\in [0, \tau(x)]\}$$
to the unit tangent bundle, defined by saying that the image of
$((x,y), s)$ is $g^sv_{x,y}$, the image of $v_{x,y}$ under time $s$
of the geodesic flow. This surjection is not one-to-one
on boundary points; by definition of $\tau$, $((x,y), \tau(y))$
has the same image as $(F(x,y), 0)$. On $X^\tau$ the geodesic flow corresponds
the translation of the parameter $t$.

\subsection{Transfer operators.}
Let $a$ be an analytic function on $S\X$. Let $A$ be the real-analytic function on $X$, defined by
$$A(z_1, z_2)=\int_0^{\tau(z_2)}a(z_1, z_2, t) dt$$
if $(z_1, z_2)\in X\subset B\times B$. In other words, $A$ is the Radon transform
$A={\cal R}(\chi a)$ defined in \ref{RT}, and $\chi$ is the cut-off function
$\chi((x,y), s) :=\sum_{i,j}\1_{ J_i}(z_2)\1_{K_j}(z_1)\1_{(0,\tau(y))}(s)$.
If $z=(z_1, z_2)$ is a
periodic point of period $n$ for $F$, and if $\gamma$ is the corresponding closed geodesic, then
$$S_n A(z):=\sum_{k=0}^{n-1}A(F^k z)=\int_\gamma a.$$

Following \cite{Rugh92}, we introduce a family of transfer
operators $L_{s,z}$ ($s, z\in\C$), acting on $\chi(\omega_1,
\omega_2)$ which are analytic with respect to $\omega_2$ and
analytic functionals (distributions) with respect to $\omega_1$.
We refer to \cite{Rugh92} for the precise definition of the Banach
space to work with. The transfer operator $L_{s,z}$ acts on $\chi$
by
\begin{multline*}\int L_{s,z}\chi(z_1, z_2) \psi(z_1)dz_1=\\
\sum_{k_{-1}, k_0}g^{\prime}_{k_0}(z_2)^{1+(s/2-1)}\int_{K_{k_{-1}}}\psi(g_{k_0}^{-1}\omega_1)
\chi(\omega_1, g_{k_0}z_2)e^{zA(\omega_1, g_{k_0} z_2)}g^{-1\prime}_{k_0}(\omega_1)^{s/2}d\omega_1
\end{multline*}
($\psi$ being an arbitrary analytic test function.)
The sum runs over all $k_{-1}, k_{0}$ such that the rectangle $K_{k_{-1}}\times
J_{k_0}$ contains a point $(\omega_1, \omega_2)$ with $F(\omega_1, \omega_2)=(z_1, z_2)$. Here the notation $g_{k_0}$
means the element of $\Gamma_0$ such that $F(\omega_1, \omega_2)=(g_{k_0}^{-1}\omega_1, g_{k_0}^{-1}\omega_2)$ if $(\omega_1, \omega_2)\in K_{k_{-1}}\times
J_{k_0}$.

\subsection{\label{DETS} Determinants and zeta functions.}
Apart from the introduction of the weight $A$, our transfer operator also differs from Rugh's
by the terms $g^{\prime}_{k_0}(z_2)^{(s/2-1)}$ and $g^{-1\prime}_{k_0}(\omega_1)^{s/2}$.
All his arguments can be adapted with obvious modifications to this situation. Let us introduce the notations $\tau_1(\omega)= -\log g^{-1\prime}_{k_0}(\omega_1)$ and $\tau_2(\omega)=\log g^{\prime}_{k_0}(\omega_2)$
if $\omega=(\omega_1, \omega_2)\in X$ with $\omega_2\in J_{k_0}$.

In paragraph 4.4 of his paper, Rugh shows that $L_{s, z}$ is a nuclear
(trace class) operator.
One can take the determinant of $I-L_{s,z}$:
$$d(s,z):=det(I-L_{s,z})=\prod (1-\lambda^{(i)}_{s,z})^{m^{(i)}}$$
where the product runs over the spectrum of $L_{s,z}$, and
$m^{(i)}=m^{(i)}_{s, z}$ denotes the multiplicity of $\lambda^{(i)}$. The
eigenvalues do not necessarily depend analytically on $(s,z)$, as
the multiplicity may vary; the determinant $d(s,z)$, however, is an analytic function of $(s,z)$, as follows from the next paragraph:

For arbitrary $(s_0,z_0)$, consider, for every $i$, a neighbourhood $V_i$ of
$\lambda^{(i)}_{s_0,z_0}$, such that the $V_i$s are all pairwise
disjoint. Let $P^i_{s,z}$ be the spectral projector on $V_i$ for
the operator $L_{s,z}$: $P^i_{s,z}$ depends analytically on
$(s,z)$, in a neighbourhood of $(s_0, z_0)$. Call $B^i_{s,z}=L_{s,z}P^{(i)}_{s,z}$: they are operators of rank $m^{(i)}_{s_0,z_0}$,
depending analytically on $(s,z)$ in a neighbourhood of $(s_0, z_0)$. By definition the spectrum of $B^{(i)}(s,z)$ is contained
in $V_i$.
Of course,
$$det\,(1-
B^{(i)}_{s_0,z_0})=(1-\lambda^{(i)}_{s_0,z_0})^{m^{(i)}_{s_0,z_0}}.$$
Thus one can write, in the
neighourhood of $(s_0,z_0)$,
\begin{equation}\label{det}d(s,z)=\prod_i det\,(1- B^{(i)}_{s,z})
\end{equation}

This shows that the determinant $d(s,z)$ is an entire function,
and has zeros exactly when $L_{s,z}$ has the eigenvalue $1$.

Rugh shows that the following remarkable identities hold: for all $n$, the trace
of $L_{s, z}^n$ is
\begin{equation}\label{resonance}
Tr(L_{s, z}^n) =\sum_{z, F^n z=z} \frac{e^{z S_n A(z) -sS_n \tau_1(z)+(1-s/2)S_n \tau_2 (z)}}{|det(DF^n(z)-1)|}.
\end{equation}
It follows that
\begin{prop}
\begin{equation}\label{DET2}d(s,z):=det(I-L_{s, z})=\exp\left(- \sum_{z, F^n z=z} \frac1n\frac{e^{z S_n A(z) -sS_n \tau_1(z)+(1-s/2)S_n \tau_2 (z)}}{|det(DF^n(z)-1)|} \right).\end{equation}
\end{prop}

In particular, the function
$$\frac{\partial_z d}{d}(s, 0)=-\sum_{z, F^n z=z }\frac1n\frac{S_n A(z)}{det |DF^n-1|}$$
has poles exactly when $1$ is in the spectrum of $L_{s, 0}$.

Because periodic points of $F$ correspond to closed
geodesics, we can express \ref{DET2} in terms of periodic geodesics. If $F^n z=z$ and $\gamma$ is the corresponding closed geodesics, we have $S_n A(z)=\int_\gamma a$, $S_n \tau_1(z)=
S_n\tau_2(z)=L_\gamma$. Thus, $d(s,z)$ is more or less the same as
\begin{equation}\label{e:closedg} \exp\left(-{\sum_\gamma}^\prime\sum_{p\geq 1}\frac1p \frac{e^{p(z\int_\gamma a -(s-1) L_\gamma)}}{|\sinh(pL_\gamma/2) |^2}\right)
=\prod_{\gamma}\prod_{m, n\in\N}(1-e^{z\int_\gamma a-(s+m+n)L_\gamma}).\end{equation}
A ``prime'' following a sum or a product means we are summing over primitive closed orbits. Otherwise, we sum or take the product over all
closed geodesics.

The previous formula, however, is not exactly true, because certain periodic geodesics correspond to several
different periodic orbits of $F$; namely, those going through the boundary of $X$ (there are a finite
number of them). The precise expression of $d(s,z)$ in terms of closed geodesics
given in \cite{Rugh96}, or
\cite{Mo}:
\begin{equation}\label{CLOSED}d(s,z)= \prod_{\gamma}\prod_{m, n\in\N}(1-e^{z\int_\gamma a-(s+m+n)L_\gamma})\;\;.P(s,z)\end{equation}
where the correction term is
\begin{equation}\label{Dir} P(s,z)=\frac{\prod_{c}\prod_{m, n\geq 0}(1-e^{-(s+m+n)l(c)+z\int_{c} a})}{\prod_{c'}\prod_{m, n\geq 0}(1-e^{-(s+m+n)l(c^\prime)+z\int_{c^\prime} a})};\end{equation}
the products run over a finite number of periodic orbits
that are counted several times in the Markov coding. The correction factor
on the right is analytic and non-vanishing in $\{\Re e \;s>0\}$, thus the zeros of the two functions \eqref{DET2} and \eqref{e:closedg} are the same there.

\begin{rem}In the half-plane $\{\Re e\;s\leq 0\}$, the correction factor $P(s, z)$ is more difficult to analyze
because it seems that its singularities could depend on the choice of the Markov section $X$. It was, however
shown in \cite{Rugh96} that the apparent singularities of \eqref{e:closedg},
appearing from the identity \eqref{CLOSED}, are removable.
\end{rem}

\begin{rem} \label{rr} For $z=0$ (which is the case treated by Rugh)
we obtain the relation
\begin{equation}d(s,0)= \prod_{n\in\N}\zeta_{S}(s+n)\;\;.P(s,0)\end{equation}
where $\zeta_S$ is the Selberg zeta function. In particular, $d(s,0)$ has
the same singularities as $\zeta_S$ in $\{\Re e\; s>0\}$.
\end{rem}

We focus our attention in the region $\{\Re e\; s>0\}$. In that region,
the function $\frac{\partial _z d}{d}(s, 0)$ has the same singularities as
$${\lcal_2}(s,0):= \sum_\gamma \int_{\gamma_0}a\frac{e^{ -(s-1)L_\gamma}}{|\sinh(L_\gamma/2 |^2}.$$

This shows that the singularities of $\lcal_2$
appear when $L_{s,0}$ has $1$ as an eigenvalue. In the next paragraph,
we show that this occurs for $s=\frac12\pm ir_n$. Then we identify the residues.

\begin{rem} In our conventions, $r_n\geq 0$ and we have defined the
boundary values $T_{ir_n}$ using this choice of convention. For simplicity, we will
restrict our attention to $s=s_n=\frac12 + ir_n$, but the analysis at $s=(1-s_n)=\frac12 - ir_n$ would be similar.
\end{rem}
\subsubsection{Location of the poles.}

For $s=s_n$, one can check
directly that $1$ is in the spectrum of $L_{s,0}$: the eigenspace is spanned by the functionals
\begin{equation}\label{CHI}\chi^{ir_n}_{(l)}(z_1, z_2)=\frac{dT_{ir_n}(z_1)}{|z_1-z_2|^{s_n}},\end{equation}
where $T_{ir_n}$ are the boundary values of eigenfunctions of the laplacian.

The verification is easy : let $\psi(z_1)$ be an arbitrary analytic function. What we need to check is that
$$\int \psi(z_1)\frac{dT_{ir_n}(z_1)}{|z_1-z_2|^{s_n}}dz_1=
\sum_{k_{-1}, k_0}g^{\prime}_{k_0}(z_2)^{s_n/2}\int_{K_{k_{-1}}}\psi(g_{k_0}^{-1}\omega_1)
g^{-1\prime}_{k_0}(\omega_1)^{s_n/2} \frac{dT_{ir_n}(\omega_1)}{|\omega_1-g_{k_0}z_2|^{s_n}}.
$$

The last expression can be transformed to
\begin{multline*}\sum_{k_{-1}, k_0}g^{\prime}_{k_0}(z_2)^{1+s_n}\int_{g_{k_0}^{-1}K_{k_{-1}}}\psi(z_1)
g^{-1\prime}_{k_0}(g_{k_0}z_1)^{s_n/2} \frac{dT_{ir_n}(g_{k_0}z_1)}{|g_{k_0}z_1-g_{k_0}z_2|^{s_n} }\\=
\sum_{k_{-1}, k_0}g^{\prime}_{k_0}(z_2)^{s_n/2}\int_{g_{k_0}^{-1}K_{k_{-1}}}\psi(z_1)
g^{-1\prime}_{k_0}(g_{k_0}z_1)^{s_n/2}g_{k_0}^{\prime}(z_1)^{s_n} \frac{dT_{ir_n}(z_1)}{|z_1-z_2|^{s_n}}
g_{k_0}^{-1\prime}(z_1)^{-s_n/2}g_{k_0}^{-1\prime}(z_2)^{-s_n/2}
\\ =\int \psi(z_1)\frac{dT_{ir_n}(z_1)}{|z_1-z_2|^{s_n}}dz_1.
\end{multline*}

\begin{rem} Similarly, the functionals
\begin{equation}\label{CHI-R}\chi^{ir_n}_{(r)}(z_1, z_2)=\frac{dT_{ir_n}(z_2)}{|z_1-z_2|^{s_n}}\end{equation}
are eigenvectors for the adjoint $L_{s_n, 0}^*$.
\end{rem}

Conversely, we need to know that $1$ is in the spectrum of $L_{s, 0}$
only if $s$ is one of the $s_n$; and that the multiplicity of $1$ is exactly
the multiplicity of $s_n(1-s_n)$ in the spectrum of the laplacian.
To prove this, we can use the relation with the Selberg zeta-function (Remark \ref{rr}).
For the latter we know that the zeros occur when $s(1-s)$
is in the spectrum of the laplacian, with the same multiplicity.

It follows also from the formula \eqref{CHI} relating eigenfunctions
of the laplacian to eigenfunctions of $L_{s_n, 0}$, that the characteristic
space
for the eigenvalue $1$ of $L_{s_n, 0}$ corresponds to a diagonal
matrix (there are no Jordan blocks).

\subsubsection{The residues.}
We are interested in the singularities of $\lcal_2$ in $\{\Re e\;s >0\}$,
or equivalently in the singularities of
\begin{eqnarray}\label{log}
\frac{\partial_z d(s,0)}{d(s,0)}=\sum_i \frac{\partial_z
det\,(1-B^{(i)}_{s,0})}{det\,(1-B^{(i)}_{s,0})}
\end{eqnarray}
From the previous paragraph, we know that $d(s,0)=0$ if and only if $s=\frac12\pm ir_n$
($1/4+r_n^2$ being an eigenvalue of the laplacian). For some $i$,
the operator $B^{(i)}_{s_n,0}$ has $1$ as an eigenvalue, and its multiplicity $m_i$ is the same as the multiplicity of $1/4+r_n^2$
in the spectrum of the laplacian. As in the previous paragraph, we treat
the case of $s=s_n=\frac 12+ir_n$; the case of $s=\frac12-ir_n$ would be similar except for the choice of a different convention in the definition of boundary values. We will see that the singularity
of the function \eqref{log} at $s=s_n$ is a pole; the residue must then be given by
$$m_i\frac{\partial_z\partial_s^{m_i-1} det\,(1-B^{(i)}_{s_n,0})}{\partial_s^{m_i} det\,(1-B^{(i)}_{s_n,0})}.$$

Theorem \ref{main2} will then follow directly from:
\begin{prop}
$$m_i\frac{\partial_z\partial_s^{m-1} det\,(1-B^{(i)}_{s_n,0})}{\partial_s^{m} det\,(1-B^{(i)}_{s_n,0})}=
\sum_{r_j=r_n}\frac{\int \chi a \;dPS_{ir_j}}{\int \chi\;
dPS_{ir_j}}.
$$

\end{prop}

\begin{proof}

If
$1/4+r_n^2$ is an eigenvalue of the laplacian of multiplicity $m$,
we know that $1$ is an
eigenvalue of $L_{s_n,0}$ of multiplicity
$m$.
We also know -- and this is rather important -- that
the eigenvalue $1$ corresponds to a diagonal block for $L_{s_n,0}$.

Let $V\subset\C$ be a neighbourhood of $1$ that does not meet the rest of
the spectrum of  $L_{s_n,0}$. Let $P_{s,z}$ be the spectral
projector on $V$ for the operator $L_{s,z}$. As before, denote
$L_{s,z}P_{s,z}=B_{s,z}$.
Because we have a diagonal block, $B_{s_n,0}
=P_{s_n, 0}.$
Using the previous notations, for some $i$,
the operator $B_{s_n, 0}$ is one of the $B^{(i)}_{s_n,0}$s; it has $1$ as an eigenvalue with multiplicity $m_i=m$.

In the tensor product ${\cal H}^{\wedge m}$, the projector $P_{s_n,
0}^{\wedge m}$ is of rank $1$. Let
${\cal V}_{s_n,0}\in {\cal H}^{\wedge
m}$ be the associated eigenvector; it also belongs to $Ker(I-L_{s_n, 0})^{\wedge m}$. By perturbation theory, we can
find a family ${\cal V}_{s,z}$, depending analytically on $(s,z)$ in a
neighbourhood of $(s_n, 0)$, such that $P_{s, z}^{\wedge
m}{\cal V}_{s,z}={\cal V}_{s,z}$.

We have
\begin{equation}\label{ev}(I-L_{s,z})^{\wedge m}{\cal V}_{s,z}=\lambda_{s,z}{\cal V}_{s,z}\end{equation}
with $\lambda_{s,z}=det(I-B_{s,z})$.

Similarly there is a family ${\cal T}_{s,z}$ in the dual ${\cal
H^*}^{\wedge m}$, depending analytically on the parameters, such
that
\begin{equation}(I-L_{s,z}^*)^{\wedge m}{\cal T}_{s,z}=\lambda_{s,z}{\cal T}_{s,z}\end{equation}

Differentiating \ref{ev} once with respect to the parameters, and applying ${\cal T}_{s_n, 0}$,
 we get
\begin{multline}\label{trick}(\partial\lambda_{s_n,0})\langle {\cal V}_{s_n,0}, {\cal T}_{s_n, 0}\rangle
+\lambda_{s_n,0}\langle \partial {\cal V}_{s_n,0}, {\cal T}_{s_n, 0}\rangle=\langle(\partial(I-L_{s_n,0})^{\wedge m}){\cal V}_{s_n,0}, {\cal T}_{s_n, 0}\rangle+\langle(I-L_{s_n,0})^{\wedge m})\partial {\cal V}_{s_n,0}, {\cal T}_{s_n, 0}\rangle
\end{multline}
Because $(I-L_{s_n,0}^*)^{\wedge m}{\cal T}_{s_n,0}=\lambda_{s_n,0}{\cal T}_{s_n,0}$, the second term on each side of \ref{trick} are equal, and \ref{trick}
amounts to
\begin{equation}\label{trick2}(\partial\lambda_{s_n,0})\langle {\cal V}_{s_n,0}, {\cal T}_{s_n, 0}\rangle
=\langle(\partial(I-L_{s_n,0})^{\wedge m}){\cal V}_{s_n,0}, {\cal T}_{s_n, 0}\rangle
\end{equation}
This last term vanishes if $m>1$, and thus we see that $\partial\lambda_{s_n,0}=0$. Iterating this procedure, we see that any derivative of order $<m$
of $\lambda_{s,z}$ vanishes at $(s_n, 0)$.

This proves, in particular, that the singularity of the function \ref{log}
at $s_n$ is at most a pole, and that the residue we are interested in is
$$m\frac{\partial_z \partial_s^{m-1} \lambda_{s_n, 0}}{\partial_s^m \lambda_{s_n, 0}}$$
as announced earlier.

The same procedure already applied (differentiate \ref{ev}, then apply
${\cal T}_{s_n, 0}$) gives after $m$ steps:

\begin{multline}(\partial_z \partial_s^{m-1}\lambda_{s_n,0})\langle {\cal V}_{s_n,0}, {\cal T}_{s_n, 0}\rangle=\langle(\partial_z \partial_s^{m-1} (I-L_{s_n,0})^{\wedge m}){\cal V}_{s_n,0}, {\cal T}_{s_n, 0}\rangle\\
=(-1)^m (m-1)!\sum_{k=0}^{m-1} \langle(\partial_s L)^{\wedge
k}\wedge\partial_z L\wedge (\partial_s L)^{\wedge m-1-k}{\cal V}_{s_n,0},
{\cal T}_{s_n, 0}\rangle \end{multline}
The terms where $L$ is not been differentiated disappear, because
$1-L_{s_n, 0}^*$ vanishes on ${\cal T}_{s_n, 0}$.  Similarly,
\begin{equation*}\partial_s^{m}\lambda_{s_n,0}=(-1)^m m!\langle(\partial_s L)^{\wedge m-1}V_{s_n,0}, T_{s_n, 0}\rangle\end{equation*}

We note that $\partial_z L=L\circ M_A$ (where $M_A$ denotes
multiplication by $A$) and $\partial_s L=L\circ M_\tau$.
Remembering that $L_{s_n, 0}^*{\cal T}_{s_n, 0}={\cal T}_{s_n, 0}$ we can rewrite the last two expressions
as
\begin{equation*}\partial_z \partial_s^{m-1}\lambda_{s_n,0}=(-1)^m (m-1)!\sum_{k=0}^{m-1}
\langle M_\tau^{\wedge k}\wedge M_A \wedge M_\tau^{\wedge
m-1-k}{\cal V}_{s_n,0}, {\cal T}_{s_n, 0}\rangle
\end{equation*}
and
\begin{equation*}\partial_s^{m}\lambda_{s_n,0}=(-1)^m m!\langle M_\tau^{\wedge m-1}{\cal V}_{s_n,0}, {\cal T}_{s_n, 0}\rangle.\end{equation*}

Now, we can choose to write ${\cal T}_{s_n, 0}$ as
$${\cal T}_{s_n, 0}=\wedge_{r_j=r_n}\chi^{ir_j}_{(l)}$$
and ${\cal V}_{s_n, 0}$ as
$${\cal V}_{s_n, 0}=\wedge_{r_j=r_n}\chi^{ir_j}_{(r)},$$
where $\chi^{ir_j}_{(l)},\chi^{ir_j}_{(r)}$ are associated to $T_{ir_j}$ by the formulae \eqref{CHI}, \eqref{CHI-R}.

For $r_j=r_k=r_n$, we have
\begin{multline*}\langle \tau \chi^{ir_j}_{(l)},\chi^{ir_k}_{(r)}\rangle =\int \tau(z_2) \frac{T_{ir_j}(dz_1)T_{ir_k}(dz_2)}{|z_1-z_2|^{2s_n}}
\\=\int ({\cal R}\1)(z_2) \frac{T_{ir_j}(dz_1)T_{ir_k}(dz_2)}{|z_1-z_2|^{2s_n}}
= \mu_0(s_n)^{-1}\langle \phi_{ir_j}, \phi_{ir_k}\rangle,
\end{multline*}
by the formulae of Part 5 (which could as well be applied for two different eigenfunctions of the same eigenvalue). Because the basis $(\phi_{ir_j})$ is orthonormal, this coefficient vanishes except for $j=k$.

Similarly,
\begin{multline*}\langle A \chi^{ir_j}_{(l)},\chi^{ir_k}_{(r)}\rangle =\int A(z_1,z_2) \frac{T^{ir_j}(dz_1)T_{ir_k}(dz_2)}{|z_1-z_2|^{2s_n}}
=\int ({\cal R}a)(z_1, z_2) \frac{T_{ir_j}(dz_1)T_{ir_k}(dz_2)}{|z_1-z_2|^{2s_n}},
\end{multline*}
and if $j=k$ this is exactly the Patterson-Sullivan distribution applied to $a$.

We finally find the expression of the residue.
$$m\frac{\partial_z \partial_s^{m-1} \lambda_{s_n, 0}}{\partial_s^m \lambda_{s_n, 0}}=
\sum_{r_j=r_n}  \frac{\int ({\cal R}a)(z_1, z_2) \frac{T_{ir_j}(dz_1)T_{ir_j}(dz_2)}{|z_1-z_2|^{2s_n}}}{\int ({\cal R}1)(z_1, z_2) \frac{T_{ir_j}(dz_1)T_{ir_j}(dz_2)}{|z_1-z_2|^{2s_n}}}
$$
which is what we expected in terms of Patteron-Sullivan distributions.

\end{proof}

\section{\label{SELBERG} Classical Selberg trace formalism}

 We now begin  the Selberg trace formalism proof of Theorem
\ref{main2}.  To prepare for the proof, we review  the standard
theory of the Selberg zeta function and trace formula and then
give a non-standard proof which will be generalized in the next
section.

 As above, we  denote by $\{\phi_{ir_k}\}$ an orthonormal basis
of $\Lap$-eigenfunctions on $\Gamma \backslash G/K$, with
associated eigenvalues $\lambda_k = s_k(1 - s_k)$ with  $s_k =
\frac{1}{2} + i r_k$. In particular the trivial eigenvalue
$\lambda_0 = 0$ corresponds to $s_0 = 0, 1$ and $r_0 = \pm
\frac{i}{2}.$

\subsection{Standard Selberg zeta function}

We now review the analytic continuation and polar analysis of the
Selberg zeta function.  We refer to \cite{V} for background.

The Selbert zeta function is defined by
$$Z(s) = \Pi_{\{P\}} \Pi_{k = 0}^{\infty} (1 - N(P)^{-s- k}), \;\; \Re e\;s > 1$$ where $\{P\}$ runs over conjugacy classes of
primitive hyperbolic elements and where $N(P) = e^{ L_P}$ where
$L_P$ is the length of the corresponding geodesic.

The logarithmic derivative of the Selberg zeta function
$\frac{1}{s - 1/2}\; \frac{Z'}{Z}(s)$ is defined for $\Re e\;s >1$ by
the formula (see \cite{V}, (5.1.5))  $$\lcal(s; 1):= \frac{1}{s -
1/2}\; \frac{Z'}{Z}(s)  = \sum_{\gamma} \frac{L_{\gamma_0}}{\sinh
L_{\gamma}/2} e^{-(s - 1/2) L_{\gamma}} .
$$
In this formula, we sum over all closed orbits  $\gamma$ of the
geodesic flow and $L_{\gamma}$ is the (positive) length of
$\gamma$.

\begin{thm}
$\frac{Z'(s)}{Z(s)}$
 admits a meromorphic continuation to $\C$ with poles at the
points $s = s_n$ together with the `trivial poles' at $s = - k$,
$k = 0, 1, 2, 3, \dots$.
\end{thm}

\begin{proof} We review a few features of the standard proof to draw attention to
some important technical issues which might be confusing for the
more general versions to come. By definition, we have for $\Re e\;s
> 1$,
$$\begin{array}{l} \frac{1}{s - 1/2}\;
\frac{Z'(s)}{Z(s)} - \frac{1}{a - 1/2}\; \frac{Z'(a)}{Z(a)} =
\sum_{\gamma} \frac{L_{\gamma_0}}{\sinh L_{\gamma}/2}
\{\frac{1}{2s - 1} e^{-(s - 1/2) L_{\gamma}} - \frac{1}{2a - 1}
e^{-(a - 1/2) L_{\gamma}}\}. \end{array} $$ See for instance
\cite{V}, (5.1.5).

To analytically continue the formula,  one applies the Selberg
trace formula (cf. \cite{V} Theorem 5.5.1) with the  test function
$$h(\frac{1}{4} + r^2; s; a) = \frac{1}{(s - \frac{1}{2})^2 +
r^2} -  \frac{1}{(a - \frac{1}{2})^2 + r^2}. $$ The Fourier
transform of $h(\frac{1}{4} + r^2; s; a)$ is
$$g(u; s; a) = \frac{1}{2s - 1} e^{-(s - \frac{1}{2})|u|} - \frac{1}{2a - 1} e^{-(a - \frac{1}{2})|u|}
. $$ We  note  that the rate of decay of $h(\frac{1}{4} + r^2; s;
a)$ as $r \to \infty$ reflects the singularity of $|u|$ at $u =
0$.  In the case of a smooth compact quotient, the result is (see
\cite{V} Theorem 5.1.1; see also \cite{Sa1})
\begin{equation} \label{STF2} \begin{array}{lll}\frac{1}{(s - \frac{1}{2})} \frac{Z'(s)}{Z(s)} -
 \frac{1}{(a - \frac{1}{2})} \frac{Z'(a)}{Z(a)}&=
& \frac{Vol(\Gamma \backslash G)}{\pi} \sum_{k = 0}^{\infty}\left(
\frac{1}{s + k} -
\frac{1}{a + k}\right) \\ & & \\
& + & \sum_{n = 0}^{\infty} \left( \frac{1}{(s - \frac{1}{2})^2 +
r_n^2} - \frac{1}{(a - \frac{1}{2})^2 + r_n^2} \right).
\end{array} \end{equation}

We note that the eigenvalue series on the right side would diverge
if we only used the formula for $\frac{Z'(s)}{Z(s)}$, but it
converges (away from poles) if we subtract $\frac{Z'(a)}{Z(a)}$ or
take one derivative.

These formulae give a meromorphic  continuation of
$\frac{Z'(s)}{Z(s)}$ to $\C$ and show that the poles occur at
values of $s$ for which there exists an eigenvalue $\lambda_n$
satisfying $\lambda_n = s(1 - s)$, or at negative integers.

\subsection{Convolution operator approach}

As sketched above, the Selberg trace formula involves a Fourier
transform duality. We will need a more group theoretic approach
for the generalizations in the next section, namely the approach
in \cite{GGP} to the
 Selberg trace formula as a formula for the
trace of the convolution operator corresponding to a
$K$-bi-invariant function $\chi$.

We denote by $S_{0,0}(G)$ the continuous functions satisfying
$\chi(k_1 g k_2) = \phi(g)$ for all $k_1, k_2 \in K$. The
associated convolution operator is defined by
$$R_{\chi}  = \int_G \chi(g) R_g dg,$$
where $R_g f(x) = f(x g)$. There exists a unique (up to scalars)
eigenfunction $\Psi_s$ of $\Omega$ of eigenvalue $s (1 - s)$ in
$S_{0, 0}$.   The spherical transform $S : C_0^{\infty} (G) \cap
S_{0, 0} \to PW_m $ is defined by
$$S f(s) = \int_G f(g) \Phi_{ s} (g) dg. $$
Its range is the subspace of the  Paley-Wiener space
\begin{equation} \label{PW} PW(\C) = \{f \in \ocal(\C): \exists k\; \forall  N > 0:
|f(x + i y)| \leq C e^{k |x|} (1 + |y|)^{-N} \} \end{equation}
with a certain symmetry which we will not need to recall here (see
\cite{Z}, p.31). Here, $\ocal(\C)$ denotes the holomorphic
functions on $\C$.

We also denote the Mellin transform $M: C_0^{\infty}(A) \to
PW(\C)$ by
$$ Mf (s) = \int_0^{\infty} f(a) a^s \frac{da}{a},$$
where we identify $f(a)$ as a function of the top diagonal entry
of $a$.  Note the non-standard sign of the exponent, which is
chosen to be consistent with \cite{L, Z}.

The basic Selberg trace formula for a smooth compact quotient (in
the form stated in \cite{GGP}) states that
\begin{equation} \label{STF1} \sum_{r_k} S \chi(2i r_k) = Vol(\Gamma \backslash G)
\chi(e) + \sum_{\{\gamma\}} \int_{G_{\gamma} \backslash G}
\chi(g^{-1} \gamma g) dg, \;\;\; \chi \in S_{0,0} \end{equation}
where the sum runs over the principal and complementary series
representations (counted with multiplicity), where $G_{\gamma}$ is
the centralizer of $\gamma$ in $G$ (similarly for $\Gamma$) For
$\gamma \not= e$, $\Gamma_{\gamma} \backslash G_{\gamma}$   is a
closed geodesic.

The orbital integral on the right side of (\ref{STF1})
 may be expressed  in terms of the so-called Harish-Chandra transform as follows: If
$\chi \in S_{00}$, there exists  $\chi^D$ on $D = G/K$ such that
$\chi(g) = \chi^D(g \cdot 0)$ where $\chi^D(r e^{i \theta}) =
\chi^D(r).$
 In the proof of
\cite{Z}, Proposition 2.6, it is shown that $\chi(n_u^{-1} a n_u)
=  \chi^D( \left| \frac{u + i}{u + i \omega} \right|),$ with
$\omega = \frac{a +
a^{-1}}{a - a^{-1}}$. With some
routine manipulation (\cite{Z}, pages 55-56), we get
\begin{equation} \label{ROUTINE0} \begin{array}{lll}
H \chi(a)& = & |a - a^{-1}| |\omega|  \; \int_{\R}  \chi^D(
\left(1 + \frac{v-1}{u^2  + 1} \right)^{1/2}) du,
\end{array}
\end{equation}
and thus
\begin{equation} \label{CCT} \int_{G_{\gamma} \backslash G}
\chi(g^{-1} \gamma g) dg  = \frac{ Vol(\Gamma_{\gamma} \backslash
G_{\gamma})}{|a_{\gamma} - a_{\gamma}^{-1}|} H \chi (a_{\gamma}),
\end{equation}  where
  $ Vol(\Gamma_{\gamma} \backslash G_{\gamma}) $ is the
length of the closed geodesic. We further have $S \chi= M H \chi$,
so we finally obtain
\begin{equation} \label{STF1a} \sum_{r_k} M H  \chi (2ir_k) = Vol(\Gamma \backslash G)
 \chi^{D}(0
 ) + \sum_{\{\gamma\}} \frac{L_{\gamma_0}}{\sinh
L_{\gamma}/2} \; H \chi(a_{\gamma}).
\end{equation}
This approach  leads most naturally to the zeta function
\begin{equation} \label{RCAL00} \rcal(s; 1)  := \sum_{\gamma}
\left(\frac{L_{\gamma_0}}{\sinh L_{\gamma/2}} \right)\;
 (\cosh L_{\gamma}/2)^{-2(s - 1/2)}.
\end{equation}
In generalizations to non-constant automorphic forms, we begin
with (\ref{RCAL00}) and then  relate it (and its generalizations
(\ref{RCALMa}) to non-constant automorphic forms)  to the usual
zeta functions $\lcal(s; \sigma)$.

\section{\label{GSELBERG} Dynamical zeta functions and Selberg trace formalism}

This section is concerned with the zeta functions
\begin{equation} \label{STEVESZETA} \lcal(s; \sigma) := \sum_{\gamma}
\left(\frac{\int_{\gamma_0} \sigma ds}{\sinh L_{\gamma/2}}
\right)\; e^{-(s - 1/2) |L_{\gamma}|}, \;\;\; (\Re e\;s > 1).
\end{equation}

\begin{thm} \label{main2a} For each automorphic form $\sigma = \phi_{ir_k}, X_+ \phi_{ir_k}, \psi_m$,
$\lcal(s; \sigma)$ is absolutely convergent in $\Re e\;s > 1$ and
admits a meromorphic continuation to $\C$. Except for the trivial
representation $\sigma \equiv 1$, the only poles in $\Re e\;s
> 0$ occur at  values $s = \frac{1}{2} + i r$ for which $\frac{1}{4} + r^2 $ is an eigenvalue
of $\Lap$, and the residue is given by
$$\mbox{Res}_{s = \frac{1}{2} + i r}\; \lcal(s; \sigma) = \mu_0(\frac{1}{2} + ir)
\sum_{j:\; r_j^2 = r^2} \int_{\Gamma \backslash G} \sigma dPS_{i
r_{j}}.
$$

\end{thm}

This proves a  special case of Theorem  \ref{main2} in which the
function $a$ has components in a finite number of irreducible
representations. We briefly sketch the extension to analytic
symbols in the final section.

The proofs are based on a generalized Selberg trace formula
introduced in \cite{Z} for  the traces $Tr \sigma R_{\chi}$ on
$L^2(\Gamma \backslash G)$ of the composition of $R_{\chi}$ with
multiplication by $\sigma$.   Here, $\sigma$ is a Casimir
eigenfunction of weight $m$ and $R_{\chi}$ is a convolution
operator with kernel $\chi \in S_{m, n}(G)$, where $S_{m, n}$
denotes the functions $\chi(g)$ on $G$ satisfying $\chi(k_{\psi} g
k_{\theta}) =  = e^{i m \psi} e^{i n \theta} \chi(g)$, where
$k_{\theta} = exp \theta W \in K$.  The eigenspaces of $\Omega$ on
$S_{m, n}(G)$
are one-dimensional, spanned by the {\it spherical function}
$\Phi_{m, n, s}$ of $\Omega$-eigenvalue $s(1 - s)$.  We will only
be considering the case $n=0$, and denote the associated
normalized spherical function by $\Phi_{m, s}$. Our normalization
follows \cite{H, Z}. The spherical transform $S_m : C_0^{\infty}
(G) \cap S_{m, 0} \to PW_m $ is defined by
$$S_m f(s) = \int_G f(g) \Phi_{-m, s} (g) dg. $$
Its range is the subspace of the  Paley-Wiener space (\ref{PW})
with a  symmetry depending on $m$ which we will not need to recall
here (see \cite{Z}, p.31).

We will also need a variety of Harish-Chandra transforms which
depend on the weight $m$ and also on the type of representation
$\pcal_{ir}, \dcal^+_m$. There is a canonical one, defined as
follows: Let  $\chi \in S_{m, 0}$ and let (see \cite{Z} page 57
for (i) and page 49 for (ii)):
\begin{equation}\label{HCHa}
\begin{array}{lll}  H_m \chi (a) &= & |a - a^{-1}|
\int_{-\infty}^{\infty} \chi^D( \frac{u + i}{u + i \omega})
 du,\;\;\;(\omega = \frac{a +
a^{-1}}{a - a^{-1}})
\end{array} \end{equation}
Here,
 if $\chi \in S_{m, 0}$ then there exists $\chi^D$ on $D =
G/K$ such that $\chi(g) = \chi^D(g \cdot 0)$ where $\chi^D(r e^{i
\theta}) = e^{i \frac{m}{2} \theta} \chi_D(r).$
 In the proof of
\cite{Z}, Proposition 2.6, it is shown that
$$\begin{array}{lll}\chi(n_u^{-1} a n_u) = \chi^D( \frac{u + i}{u + i \omega})= e^{i \frac{m}{2}
\theta(a, u)} \chi^D( \left| \frac{u + i}{u + i \omega} \right|).
\end{array}$$
 with
\begin{equation} \label{EITHETA} e^{i \theta(a, u)} = \frac{(u + i) (u^2 + \omega^2)^{1/2}}{(u
+ i \omega) (u^2 + 1)^{1/2}},\;\; \mbox{where} \;\; \omega =
\frac{a + a^{-1}}{a - a^{-1}}, \;\; v = \omega^{-2}.
\end{equation} With some routine manipulation (see \cite{Z}, pages
55-56), we get
\begin{equation} \label{ROUTINE} \begin{array}{lll}
H_m \chi(v)& = & |a - a^{-1}| |\omega| v^{m /4} \; \int_{\R}  (1 +
\frac{v-1}{u^2 + 1})^{-m/4} \chi^D( \left(1 + \frac{v-1}{u^2  + 1}
\right)^{1/2}) du.
\end{array}
\end{equation}

From the Selberg trace formalism viewpoint, it turns out to be
most natural  to work first with auxiliary  dynamical
zeta-functions $\rcal(s; \sigma)$  that do not seem to arise in
the thermodynamic formalism. When $\sigma_m$ has weight $m$ we put
\begin{equation} \label{RCALMa} \rcal(s; \sigma_m)  := \sum_{\gamma}
\left(\frac{\int_{\gamma_0} \sigma_m}{\sinh L_{\gamma/2}}
\right)\; \left(\tanh L_{\gamma}/2  \right)^{m/2}
 (\cosh L_{\gamma}/2)^{-2(s - 1/2)}.
\end{equation}
We then  express $\lcal(s; \sigma)$ in terms of $\rcal(s; \sigma)$
to obtain results on the analytic continuation of the latter. This
somewhat circuitous route comes about because the trace formula is
on the `quantum level' and therefore does not quite produce the
`classical' zeta-function.

\subsection{Forms of weight $0$ in $\pcal_{ir}$}

In this section, we prove Theorem \ref{main2} for the case $\sigma
= \lcal(s; \phi_{ir_k})$.

In this case the auxiliary zeta function has the form
\begin{equation} \label{RCAL0} \rcal(s; \phi_{ir_k})  := \sum_{\gamma}
\left(\frac{\int_{\gamma_0} \phi_{ir_k}}{\sinh L_{\gamma/2}}
\right)\;
 (\cosh L_{\gamma}/2)^{-2(s - 1/2)}.
\end{equation}

\begin{thm} \label{RCAL0THM} $\rcal(s; \phi_{ir_k})$ admits a meromorphic continuation to
$\C$ with poles at $s = \frac{1}{2} + i r_j - k, k = 0, 1, 2,
\dots$, where $\frac{1}{4} + r_j^2$ is an eigenvalue of $\Lap$,
and with
$$\mbox{Res}_{s = \frac{1}{2} + i r}\; \rcal(s; \phi_{ir_k}) =  \mu_0(ir + \frac{1}{2}) \sum_{j:\; r_j^2 = r}
\frac{1}{2}\langle \phi_{ir_k}, PS_{r_j}\rangle \; .$$

\end{thm}

\begin{proof} We assume throughout that $\phi_{ir_k} \perp 1$, so
that the identity term on the $\sum_{\gamma \in \Gamma}$ side of
the trace formula vanishes and so that the trivial representation
term with $r = \frac{i}{2}$ also vanishes.  After the proof, we
remark on the case $\phi_{ir_k} \equiv 1$.

By    Proposition 2.12 of \cite{Z} (applied in the continuous
series case), we have
\begin{equation}\label{DSTFa}  \sum_{n = 0}^{\infty} \langle Op(\phi_{ir_k}) \phi_{ir_n},
\phi_{ir_n} \rangle M H_0 \chi(2ir_n) = \sum_{\gamma}
\left(\frac{\int_{\gamma_0} \psi_m }{\sinh L_{\gamma/2}} \right)
H_{ir_k}^c \chi(a_{\gamma}),
\end{equation} where $H_0$ is defined by (\ref{HCHa})-(\ref{ROUTINE}),
and where  (see \cite{Z} page 57 for (i) and page 49 for (ii)):
\begin{equation}\label{HCHRK} \begin{array}{lll}
H_{ir_k}^c \chi(a) & : = & |a - a^{-1}| \int_{-\infty}^{\infty}
F_{ir_k,0}(\frac{u - i}{-2 i})  \chi^D( \frac{u + i}{u + i
\omega})
 du.
\end{array} \end{equation}
 Here,
 $a = e^{L/2}$ and $
F_{ir_k,0}(\frac{u - i}{-2 i}) $ is defined in (\ref{WTZERO}). We
note that the identity term on the right side vanishes by
orthogonality.

\begin{rem} (i)  We note that we do not use Proposition 2.10 of
\cite{Z}, which gives a less convenient zeta function. Although
Proposition 2.12 of \cite{Z} is only stated for symbols in the
discrete series, it is valid for the continuous series as long as
we use the corresponding expressions (given in \cite{Z} Corollary
2.4)  for the integrals $I_{\gamma}(\sigma)(n_u)$ in
\cite{Z}(2.2). \medskip

\noindent(ii) A priori, the right side of (\ref{HCHRK}) should
also contain the term
$$ |a - a^{-1}| \int_{-\infty}^{\infty}
G_{ir_k,0}(\frac{u - i}{-2 i})  \chi^D( \frac{u + i}{u + i
\omega})
 du,$$
 but $\chi^D$ is a radial function since it has weight zero and $\chi^D( \left| \frac{u + i}{u + i \omega}
 \right|)$
 is even in $u$ while $G_{ir_k,0}(\frac{u - i}{-2 i})$ is odd.
 Hence this integral vanishes (cf. Proposition 2.7 of \cite{Z}).
 \end{rem}

By (\ref{ROUTINE}), we have
\begin{equation} \label{ROUTINEa} \left\{ \begin{array}{lll}
H_0 \chi(v)& = & |a - a^{-1}| |\omega| \; \int_{\R} \chi^D(
\left(1
+ \frac{v-1}{u^2  + 1} \right)^{1/2}) du, \\ && \\
H_{2ir_k}^c \chi(v)& = & |a - a^{-1}| |\omega| \; \int_{\R}
F_{ir_k,0}(\frac{u - i}{-2 i}) \chi^D( \left(1 + \frac{v-1}{u^2 +
1} \right)^{1/2}) du.
\end{array} \right.
\end{equation}
  We now define $\chi_s(g) \in S_{0, 0}$ by the rule
that
$$\chi_s^D(r) : = \frac{(1 - r^2)^{s}}{\mu_{ir_k}^c(s)},\;\;\; 0 \leq r \leq 1.$$
Using (\ref{ROUTINEa}) and the fact that $|a - a^{-1}| |\omega| =
(a + a^{-1})$, we obtain
 \begin{equation}
\label{MOREHMa} \left\{
\begin{array}{lll}(i) \;\;H_{2ir_k}^c \chi_s(a)  & = &   (a + a^{-1})^{-2(s - 1/2)}, \\ & & \\(ii)
H_0 \chi_s(a) = & = & \frac{\mu_0(s)}{\mu_{ir_k}(s)} \;  \left( a
+ a^{-1} \right)^{-2(s - 1/2)}
\end{array} \right..
\end{equation}
If we substitute $\chi_s$ into  the  right side of the trace
formula (\ref{DSTFa}), we obtain  the desired zeta-function
$\rcal(s; \phi_{ir_k})$. Therefore, the  left side of the trace
formula (\ref{DSTF}) gives a meromorphic continuation of $\rcal(s;
\phi_{ir_k}
)$.  By  Theorem \ref{main2}, we have

$$\begin{array}{lll}\langle Op(\phi_{ir_k})\phi_{ir}, \phi_{ir}\rangle &=
&
 \langle
\phi_{ir_k},  PS_{ir} \rangle  \mu_{ir_k}^c(\frac{1}{2} + ir),
\end{array}$$ hence
\begin{equation} \rcal(s;
\phi_{ir_k}) = \sum_{n = 0}^{\infty} \langle \phi_{ir_k},
PS_{ir_n}\rangle \mu_{ir_k}^c(\frac{1}{2} + ir_n)\; M H_0
\chi_s(2ir_n).
\end{equation}
By  (\ref{MOREHMa}(ii)) we have
$$  \begin{array}{lll} M H_0 \chi_s(2ir) &= & \frac{\mu_0(s) }{\mu_{ir_k}^c(s)}  \int_0^{\infty} a^{2ir}  (a + a^{-1})^{-2(s -
1/2)} \frac{da}{a}\\&& \\
& = & \frac{\mu_0(s) }{\mu_{ir_k}^c(s)} \int_{-\infty}^{\infty}
e^{2ir t}  (\cosh t)^{-2(s - 1/2)} dt\\ && \\
& = & \frac{\mu_0(s) }{\mu_{ir_k}^c(s)} \frac{\Gamma(s -
(\frac{1}{2} + ir )
 ) \Gamma(s - (\frac{1}{2} - ir))}{\Gamma(2s - 1)}.
\end{array}$$ For the last line we refer to  \cite{Z} (p. 60).

In conclusion, we obtain (at least formally)
\begin{equation}\label{RCALCONT}  \rcal(s;
\phi_{ir_k}) = \sum_{n = 0}^{\infty} \langle \phi_{ir_k},
PS_{ir_n}\rangle \frac{\mu_0(s)\mu_{ir_k}^c(\frac{1}{2} + ir_n)
}{\mu_{ir_k}^c(s)} \frac{\Gamma(s - (\frac{1}{2} + ir_n )
 ) \Gamma(s - (\frac{1}{2} - ir_n))}{\Gamma(2s - 1)}
\end{equation}
We note that $ \frac{\Gamma(s - (\frac{1}{2} + ir )
 ) \Gamma(s - (\frac{1}{2} - ir))}{\Gamma(2s - 1)} = B(s - (\frac{1}{2} +
 ir), s - (\frac{1}{2} - ir)). $
As above, we assume that $\phi_{ir_k} \perp 1$, so  that the
trivial representation term vanishes. Regarding the convergence of
the right side, we note that by (\ref{BIGDISPLAY}) and
(\ref{GAMASYM}), as $|r_n| \to \infty$,
$$ \left\{ \begin{array}{l} \mu_{ir_k}^c(\frac{1}{2} + ir_n)
\sim r_n^{-1/2} \\ \\
\Gamma(s - (\frac{1}{2} + ir_n )
 ) \Gamma(s - (\frac{1}{2} - ir_n))
\sim e^{- \frac{\pi}{2} (|\Im s + r_n| + |\Im s - r_n|)} |r_n +
\Im s|^{-  \Re s - 1 }|r_n - \Im s|^{-  \Re s - 1 }.
\end{array} \right.
$$
Since $\langle \phi_{ir_k}, PS_{ir_n}\rangle  =  O_{r_k}(1)$ as $n
\to \infty$, it follows that the series converges absolutely in
the critical strip away from the poles and defines a meromorphic
function.

There are simple poles at $s = \frac{1}{2} \pm  i r_n$ where
$\frac{1}{4} + r_n^2$ is an eigenvalue of $\Lap$. In the case
where the multiplicity of the eigenvalue equals one, the  residue
at $s = \frac{1}{2} + i r_n$ equals
$$\langle \phi_{ir_k},
PS_{ir_n}\rangle \frac{\mu_0(\frac{1}{2} +
ir_n)\mu_{ir_k}^c(\frac{1}{2} + ir_n) }{\mu_{ir_k}^c(\frac{1}{2} +
ir_n))} \frac{ \Gamma(2ir_n))}{\Gamma(2ir_n)} = \mu_0(\frac{1}{2}
+ ir_n) \langle \phi_{ir_k}, PS_{ir_n}\rangle, $$ as stated. In
the case of a multiple eigenvalue one sums over an orthonormal
basis of the eigenspace.

\end{proof}

  \subsection{$\lcal(s;\phi_{ir_k})$}

Now we deduce properties of $\lcal(s; \phi_{ir_k})$ from those of
$\rcal(s; \phi_{ir_k})$.

  We
introduce the measure
$$d \Theta(L; \phi_{ir_k}) = \sum_{\gamma}
\left(\frac{\int_{\gamma_0} \phi_{ir_k} ds}{\sinh L_{\gamma}/2}
\right)\; \;  \delta (L - L_{\gamma}). $$ We note that
\begin{equation} \label{LAPTa} \left\{
\begin{array}{l} \lcal(s; \phi_{ir_k}) = \int_0^{\infty}
e^{-(s - \frac{1}{2}) L} \; \;d \Theta(L; \phi_{ir_k}), \\ \\
\rcal(s; \phi_{ir_k}) = \int_0^{\infty} \left((\frac{e^{L/2} +
e^{-L/2}}{2})^2\right)^{-(s - \frac{1}{2})} d\Theta(L;
\phi_{ir_k}).
\end{array} \right.
\end{equation}

\begin{lem}\label{ZSUMa}  We have:
$$\lcal(s; \phi_{ir_k}) =  \sum_{n = 0}^{\infty} B_m (s, n)  \rcal(s + n; \phi_{ir_k}),$$
where
$$B_m(s, n)= 2^{-s + \frac{1}{2}}  2^n
\{\sum_{m, k_1, \dots, k_m = 0; k_1 + \cdots k_m = n}^{\infty} {2s
- 1 \choose m} {\frac{1}{2} \choose k_1 + 1} \cdots {\frac{1}{2}
\choose k_m + 1} \}$$

\end{lem}

\begin{proof} By elementary manipulation, we have
 \begin{equation} \label{LCAL1m0} \lcal(s,\phi_{ir_k})
= \int_0^{\infty} \left( 1 + e^{- L} \right)^{2s - 1} d \Theta (L;
s; \phi_{ir_k}),
\end{equation}
where
$$d \Theta(L, s; \phi_{ir_k}) = \sum_{\gamma}
\left(\frac{\int_{\gamma_0} \phi_{ir_k} }{\sinh L_{\gamma}/2}
\right)\; \left((\cosh L_{\gamma}/2)^2 \right)^{-(s -
\frac{1}{2})}\;  \delta (L - L_{\gamma}). $$ We then change
variables to  $y =  (\cosh L/2)^2$, and note that $e^{-L/2} =
\sqrt{y} - \sqrt{y - 1}$ to obtain,
\begin{equation} \label{LCAL1d0} \begin{array}{lll} \lcal(s,\phi_{ir_k}) &=& \int_0^{\infty}
\left( 1 + ( \sqrt{y} - \sqrt{y - 1})^2\right)^{2s - 1} d \Psi (y;
s; \phi_{ir_k}) \\ &&\\
 &=& \int_0^{\infty} (2 y)^{2s - 1}
\left( 1 - \sqrt{1 - \frac{1}{y}}\right)^{2s - 1} d \Psi (y; s;
\phi_{ir_k})
\end{array} ,
\end{equation}
where
$$d\Psi (y; s;  \phi_{ir_k}) = \sum_{\gamma} \left(\frac{\int_{\gamma_0} \phi_{ir_k}}{\sinh L_{\gamma}/2}
\right)\; \left(y_{\gamma} \right)^{-(s - \frac{1}{2})}\;  \delta
(y - y_{\gamma}).$$

By  repeated use of  the binomial theorem,  there exist
coefficients $B_m(s, n)$ such that
$$\left( 1 - \sqrt{1 - \frac{1}{y}}\right)^{2s - 1}  = y^{- (2s - 1)}  \sum_{n=0}^{\infty} B_m(s, n) y^{-n}. $$
Cancelling the factors of $y^{\pm (2s - 1)}$, we thus have
\begin{equation} \label{LCAL1ba} \begin{array}{lll} \lcal(s,\phi_{ir_k}) &=& \sum_{n=0}^{\infty} B_m(s, n) \int_0^{\infty}
y^{-n} d \Psi (y; s; \phi_{ir_k}) \\ &&\\   &=& \sum_{n =
0}^{\infty} B(s, n) \rcal(s + n; \phi_{ir_k}).
\end{array} \end{equation}

 Since the poles of  $\rcal(s + n,\phi_{ir_k})$ are the shifts by
$-n$ of the poles of $\rcal(s, \phi_{ir_k})$, and since the
non-trivial poles of $\rcal(s, \psi_m)$ in $\Re e\;s > 0$ lie only at
the points $s = \frac{1}{2} + ir $, only the term $n  = 0$ in the
series contributes non-trivial poles to the critical strip.

Writing out $B_m(s, 0)$ as a sum $\sum_{m, k_1, \dots, k_m = 0;
k_1 + \cdots k_m = n}$, we see that the only term has    $m = 0 =
k_j$ (for all $j =1, \dots, m)$. Thus,
$$Res_{s = \frac{1}{2} + ir} \lcal(s; \phi_{ir_k}) =  Res_{s = \frac{1}{2} + i r} \rcal(s;
\phi_{ir_k}). $$ This completes the proof of Theorem \ref{main2a}
in the case $\sigma = \phi_{ir_k}.$

\end{proof}

\end{proof}

\begin{rem} As a check on the calculations, we observe that in the
case $\phi_{ir_k} \equiv 1$, $ \langle \phi_{ir_k},
PS_{ir_n}\rangle = 1$ for all $n$, $\mu_{ir_k}(s) = \mu_0(s)$ and
we get
\begin{equation}\label{RCALSEL}  \begin{array}{lll} \rcal(s;
1) &= & \frac{Vol(\Gamma \backslash G)}{2\pi}\int_{\R}
\frac{\Gamma(s - (\frac{1}{2} + ir )
 ) \Gamma(s - (\frac{1}{2} - ir))}{\Gamma(2s - 1)} r (\tanh \pi r) dr  \\ && \\ && + \sum_{n =
0}^{\infty} \frac{\Gamma(s - (\frac{1}{2} + ir_n )
 ) \Gamma(s - (\frac{1}{2} - ir_n))}{\Gamma(2s - 1)}.
\end{array} \end{equation}  The series converges rapidly
 to a meromorphic function with simple poles at $s = \frac{1}{2} \pm r_n
 - k$ ($k = 0, 1, 2, \dots $), the residue at $s = \frac{1}{2} \pm r_n
 - k$ being $\frac{(-1)^k}{k!}.$ Thus, Lemma \ref{ZSUMa} shows
 that $\lcal(s; 1)$ has simple poles in the critical strip with
 residues equal to $1$. The formula (\ref{RCALSEL}) also follows from
 the standard (Fourier transform duality) Selberg trace formula \cite{V} Theorem 4.3.6) by using the integral formula
 (cf. \cite{WW}, Exercise 24)
 $$B(s - \frac{1}{2} - ir_n, s - \frac{1}{2} + ir_n) =
 \frac{1}{4^{s - \frac{1}{2}}} \int_{\R} \frac{\cos (2 i r_n u)
 du}{\cosh^{2s - 1}(u)}$$
 and the fact noted above that  $ \frac{\Gamma(s - (\frac{1}{2} +
ir )
 ) \Gamma(s - (\frac{1}{2} - ir))}{\Gamma(2s - 1)} = B(s - (\frac{1}{2} +
 ir), s - (\frac{1}{2} - ir)). $

\end{rem}

\subsection{Forms of weight $\pm 2$ in $\pcal_{ir}$}

In this case both sides of the trace formula equal zero due to
time reversibility. By Propositions \ref{TRX} and \ref{INVARIANT},
each side of the trace formula equals zero, noting that  $
\left(\frac{\int_{\gamma_0}
 X_+ \phi_{ir_k} ds}{\sinh L_{\gamma/2}} \right) + \left(\frac{\int_{\gamma_0^{-1}}
 X_+ \phi_{ir_k} ds}{\sinh L_{\gamma/2}} \right) = 0.
$

\subsection{Weight $m$ in $\dcal_m^{\pm}$}

We now prove Theorem \ref{main2a} for $\sigma = \psi_m \in
\dcal_m^+$. The anti-holomorphic discrete series case is simply
the complex conjugate and is omitted.

The proof is similar to the case $\lcal(s; \phi_{ir_k})$ but
involves the higher weight analogue zeta-function:

\begin{equation} \label{RCALM} \rcal(s; \psi_m)  := \sum_{\gamma}
\left(\frac{\int_{\gamma_0} \psi_m ds}{\sinh L_{\gamma/2}}
\right)\; \left(\tanh L_{\gamma}/2  \right)^{m/2}
 (\cosh L_{\gamma}/2)^{-2(s - 1/2)}.
\end{equation}

We begin the proof with an analysis of its meromorphic
continuation.

\subsubsection{Meromorphic continuation of $\rcal(s; \psi_m)$}

\begin{thm}\label{NONAME} $\rcal(s; \psi_m)$ admits a meromorphic continuation to
$\C$. In the critical strip, its poles occur at $s = \frac{1}{2} +
ir$ such that $\frac{1}{4} + r^2$ is an eigenvalue of $\Lap$ ,
with residue $\mu_0(\frac{1}{2} + ir) \sum_{j: \; r_j^2 = r^2}
\langle \psi_m, PS_{ir_j} \rangle.$
\end{thm}

\begin{proof}

We  study $\rcal(s; \psi_m)$ using the trace formula given in
\cite{Z},  Proposition 2.12:
\begin{equation}\label{DSTF}  \sum_{n = 0}^{\infty} \langle Op(\psi_m) \phi_{ir_n},
\phi_{ir_n} \rangle M H _m \chi(2ir_n) = \sum_{\gamma}
\left(\frac{\int_{\gamma_0} \psi_m }{\sinh L_{\gamma/2}} \right)
H_m^d \chi(a_{\gamma}),
\end{equation} where (see \cite{Z} page 57 for (i) and page 49 for (ii)): \begin{equation}\label{HCH} \left\{ \begin{array}{lll}
(i)\; H_m \chi (a) &= & |a - a^{-1}| \int_{-\infty}^{\infty}
\chi^D( \frac{u + i}{u + i \omega})
 du, \\ & & \\(ii)\;
H_m^d \chi(a) & = & |a - a^{-1}| \int_{-\infty}^{\infty} (u +
i)^{-m/2} \chi^D( \frac{u + i}{u + i \omega})
 du.
\end{array} \right.\end{equation}
 We
caution that in the definition of $H_m^d$ (\ref{HCH})(ii) we
follow a slightly different notation convention in \cite{Z}
whereby we multiply the integral by $|a - a^{-1}|$ as for $H_m$

The integral uses the notation of (\ref{HCHa})-(\ref{EITHETA}). We
simplify the expressions in (\ref{HCH}) by further using these
identities to obtain   (see also \cite{Z}, pages 55-56)
\begin{equation} \left\{ \begin{array}{lll}(i) \;\;
H_m \chi(v)& = & |a - a^{-1}| |\omega| v^{m /4} \; \int_{\R}  (1 +
\frac{v-1}{u^2 + 1})^{-m/4} \chi^D( \left(1 + \frac{v-1}{u^2  + 1}
\right)^{1/2}) du,
\\ &&\\(ii)\;\;
H_m^d \chi(v)& = & |a - a^{-1}| |\omega| v^{m /4} \; \int_{\R} (u
+ i)^{-m/2} (1 + \frac{v-1}{u^2  + 1})^{-m/4} \chi^D(\left(1 +
\frac{v-1}{u^2  + 1} \right)^{1/2}) du,
\end{array} \right..
\end{equation}

We now define $\chi_s(g) \in S_{m, 0}$ by the rule that
$$r^{-m/2} \chi_s^D(r) : = \frac{(1 - r^2)^{s}}{\mu_m^d(s)},\;\;\; 0 \leq r \leq 1,$$
where (see \cite{Z}, Proposition 3.6)
$$\mu_m^d(s) = \int_{\R} (u + i)^{-m/2} (u^2 + 1)^{- s} du =   \frac{(-i)^{m/2} \pi 2^{2s + 2 - m/2} \Gamma(-2s
+ \frac{m}{2})}{ - (2s + 1 - \frac{m}{2}) \Gamma(-s) \Gamma(-s +
\frac{m}{2})}. $$ Since  $ (1 - v)^s = (a + a^{-1})^{-2s}, $ and
$|a - a^{-1}| | \omega| v^{m /4} = (a + a^{-1}) \left( \frac{a -
a^{-1}}{a+ a^{-1}} \right)^{m/2}$,  we have \begin{equation}
\label{MOREHM} \left\{
\begin{array}{lll}(i) \;\;H_m^d \chi_s(a)  & = & \left( \frac{a - a^{-1}}{a+
a^{-1}} \right)^{m/2} (a + a^{-1})^{-2(s - 1/2)}, \\ & & \\(ii)
H_m \chi_s(a) = & = & \frac{\mu_0(s)}{\mu_m^d(s)} \; \left(
\frac{a - a^{-1}}{a+ a^{-1}} \right)^{m/2} (a + a^{-1})^{-2(s -
1/2)}
\end{array} \right..
\end{equation}

 It follows first that if we substitute $\chi_s$ into  the  right side of the trace
formula (\ref{DSTF}) is the desired zeta-function $\rcal(s;
\psi_m)$. Therefore, the  left side of the trace formula
(\ref{DSTF}) gives a meromorphic continuation of $\rcal(s;
\psi_m)$.  By  Theorem \ref{main2}, we have
$$\begin{array}{lll}\langle Op(\psi_m)\phi_{ir}, \phi_{ir}\rangle &=
&
 \langle
\psi_m,  PS_{ir} \rangle  \mu_{m}^d(\frac{1}{2} + ir),
\end{array}$$ hence
\begin{equation} \rcal(s;
\psi_m) = \sum_{n = 0}^{\infty} \langle \psi_m, PS_{ir_n}\rangle
\mu_{m}^d(\frac{1}{2} + ir_n)\; M H_m \chi_s(2ir_n).
\end{equation}
By  (\ref{MOREHM}(ii)) we have
$$  \begin{array}{lll} M H_m \chi_s(2ir) &= & \frac{\mu_0(s) }{\mu_m^d(s)}  \int_0^{\infty} a^{2ir} \left(
\frac{a - a^{-1}}{a+ a^{-1}} \right)^{m/2} (a + a^{-1})^{-2(s -
1/2)} \frac{da}{a}\\&& \\
& = & \frac{\mu_0(s) }{\mu_m^d(s)} \int_{-\infty}^{\infty} e^{2ir
t} \left( \tanh t \right)^{m/2} (\cosh t)^{-2(s - 1/2)} dt.
\end{array} $$
Putting things together, we obtain the discrete series analogue of
(\ref{RCALCONT}),
\begin{equation}\label{RCALDISCRETE}  \rcal(s;
\psi_m) = \sum_{n = 0}^{\infty} \langle \psi_m, PS_{ir_n}\rangle
\mu_{m}^d(\frac{1}{2} + ir_n)\; \frac{\mu_0(s) }{\mu_m^d(s)}
\int_{-\infty}^{\infty} e^{2ir_n t} \left( \tanh t \right)^{m/2}
(\cosh t)^{-2(s - 1/2)} dt.
\end{equation}
The integral is more complicated than its zero weight analogue,
but as  $\tanh t = 1 + r(t)$ with $r(t) = O(e^{- 2 |t|})$, we may
write \begin{equation} \label{TWOTERMS} \int_{-\infty}^{\infty}
e^{2ir_n t} \left( \tanh t \right)^{m/2} (\cosh t)^{-2(s - 1/2)}
dt = \int_{-\infty}^{\infty} e^{2ir_n t} (\cosh t)^{-2(s - 1/2)}
dt + R_2(s, r_n), \end{equation} where
\begin{equation} \label{R2}  R_2(s, r_n) = \int_{-\infty}^{\infty}
e^{2ir_n t} r(t) (\cosh t)^{-2(s - 1/2)} dt. \end{equation}  The
first term of (\ref{TWOTERMS}) gives the expression in the weight
zero case analyzed above. Hence, the sum over $r_n$ with this term
converges, and the poles and residues of $\rcal(s; \psi_m)$ on
$\Re e\;s = 1/2$ due to this term are the same as for
$$\begin{array}{l}  \langle
\psi_m,  PS_{ir_n} \rangle  \frac{\mu_0(s)\mu_{m}^d(\frac{1}{2} +
ir_n) }{\mu_m^d(s)}  \frac{\Gamma(s - ( \frac{1}{2} + ir_n)
\Gamma(s - ( \frac{1}{2} - ir_n) )}{\Gamma(2s - 1)}.
\end{array}$$   There
are simple poles at $s = \frac{1}{2} + i r_n$ and the residue is $
 \langle
\psi_m,  PS_{ir_n} \rangle  \frac{\mu_0(\frac{1}{2} + i r_n))
\mu_{m}^d(\frac{1}{2} + ir_n) }{\mu_m^d(\frac{1}{2} + i r_n))} =
 \langle
\psi_m,  PS_{ir_n} \rangle  \mu_0(\frac{1}{2} + i r_n). $ Summing
over an orthonormal basis of lowest weight vectors of $\dcal^+_m$
gives the stated expression.

To complete the proof, it is only necessary to observe that the
second integral $R_2(s, r_n)$  is holomorphic in the region $\Re s
> - \frac{1}{2}. $ It is also rapidly decaying in $r_n$.
Therefore it does not contribute any poles or residues to
$\rcal(s; \psi_m)$ in the critical strip.

\end{proof}

  \subsection{$\lcal(s; \psi_m)$}

Now we deduce properties of $\lcal(s; \psi_m)$ from those of
$\rcal(s; \psi_m)$.

 As in the weight zero case, we
introduce the measure
$$d \Theta(L; \psi_m) = \sum_{\gamma}
\left(\frac{\int_{\gamma_0} \psi_m ds}{\sinh L_{\gamma}/2}
\right)\;  \left(\tanh L_{\gamma}/2  \right)^{m/2}\;  \delta (L -
L_{\gamma}). $$ We note that
\begin{equation} \label{LAPT} \left\{
\begin{array}{l} \lcal(s; \psi_m) = \int_0^{\infty}
e^{-(s - \frac{1}{2}) L} \; \left( \tanh L/2 \right)^{-m/2}\;d \Theta(L; \psi_m), \\ \\
\rcal(s; \psi_m) = \int_0^{\infty} \left((\frac{e^{L/2} +
e^{-L/2}}{2})^2\right)^{-(s - \frac{1}{2})} d\Theta(L; \psi_m).
\end{array} \right.
\end{equation}

Because the factor $\left( \tanh L/2 \right)^{-m/2}$ is somewhat
inconvenient, we also consider
\begin{equation} \label{LAPTWID} \left\{
\begin{array}{lll} \tilde{\lcal}(s; \psi_m) & = & \int_0^{\infty}
e^{-(s - \frac{1}{2}) L} \;\;d \Theta(L; \psi_m), \\ &&\\
& = & \sum_{\gamma} \left(\frac{\int_{\gamma_0} \psi_m ds}{\sinh
L_{\gamma}/2} \right)\;  \left(\tanh L_{\gamma}/2  \right)^{m/2}\;
e^{- (s - \frac{1}{2}) L_{\gamma}}.
\end{array} \right.
\end{equation}


\begin{lem}\label{ZSUM}  We have:
$$\tilde{\lcal}(s; \psi_m) =  \sum_{n = 0}^{\infty} B_m (s, n)  \rcal(s + n; \psi_m),$$
where $B_m(s, n)$ is the same as in Lemma \ref{ZSUMa}.

\end{lem}

\begin{proof} We use similar manipulations as in the weight zero case.
We now have
 \begin{equation} \label{LCAL1m} \tilde{\lcal}(s,\psi_m)
= \int_0^{\infty} \left( 1 + e^{- L} \right)^{2s - 1} d \Theta (L;
s; \psi_m),
\end{equation}
where
$$d \Theta(L, s; \psi_m) = \sum_{\gamma}
\left(\frac{\int_{\gamma_0} \psi_m ds}{\sinh L_{\gamma}/2}
\right)\;  \left(\tanh L_{\gamma}/2  \right)^{m/2}\; \left((\cosh
L_{\gamma}/2)^2 \right)^{-(s - \frac{1}{2})}\;  \delta (L -
L_{\gamma}). $$ We  change variables as before  to  $y =  (\cosh
L/2)^2$,  and obtain as in (\ref{LCAL1d0}),
\begin{equation} \label{LCAL1d} \begin{array}{lll} \tilde{\lcal}(s,\psi_m) &=& \int_0^{\infty} (2 y)^{2s - 1}
\left( 1 - \sqrt{1 - \frac{1}{y}}\right)^{2s - 1} d \Psi (y; s;
\psi_m)
\end{array} ,
\end{equation}
where
$$d\Psi (y; s;  \psi_m) = \sum_{\gamma} \left(\frac{\int_{\gamma_0} \psi_m ds}{\sinh L_{\gamma}/2}
\right)\;  \left(\tanh L_{\gamma}/2  \right)^{m/2}\;
\left(y_{\gamma} \right)^{-(s - \frac{1}{2})}\;  \delta (y -
y_{\gamma}).$$

As in the weight zero case,
 we then have
\begin{equation} \label{LCAL1b} \begin{array}{lll} \tilde{\lcal}(s,\psi_m) &=& \sum_{n=0}^{\infty} B_m(s, n) \int_0^{\infty}
y^{-n} d \Psi (y; s; \psi_m) \\ &&\\   &=& \sum_{n = 0}^{\infty}
B(s, n) \rcal(s + n; \psi_m).
\end{array} \end{equation}

 Since the poles of  $\rcal(s + n,\psi_m)$ are the shifts by
$-n$ of the poles of $\rcal(s, \psi_m)$, and since the non-trivial
poles of $\rcal(s, \psi_m)$ in $\Re e\;s > 0$ lie only at the points
$s = \frac{1}{2} + ir $, only the term $n  = 0$ in the series
contributes non-trivial poles to the critical strip, and as above
this term has  $m = 0 = k_j$ (for all $j =1, \dots, m)$. Thus,
$$Res_{s = \frac{1}{2} + ir} \tilde{\lcal}(s; \psi_m) =  Res_{s = \frac{1}{2} + i r} \rcal(s;
\psi_m). $$

\end{proof}

To complete the proof of the theorem, we now observe that
\begin{equation} \lcal(s, \psi_m) - \tilde{\lcal}(s, \psi_m) = \sum_{\gamma} \left(\frac{\int_{\gamma_0} \psi_m ds}{\sinh
L_{\gamma}/2} \right)\;  \left[ 1 - \left(\tanh L_{\gamma}/2
\right)^{m/2} \right]\; e^{- (s - \frac{1}{2}) L_{\gamma}}.
\end{equation}
Since $ \left[ 1 - \left(\tanh L_{\gamma}/2 \right)^{m/2} \right]
= \ocal(e^{- L_{\gamma}})$ and since
\begin{equation}\sum_{\gamma} \left|  \left(\frac{\int_{\gamma_0} \psi_m ds}{\sinh
L_{\gamma}/2} \right)\;  e^{(-s - 1  + \frac{1}{2}) L_{\gamma}}
\right| < \infty, \;\;\; \Re e\;s > 0,
\end{equation}
by the prime geodesic theorem, it follows that $\lcal(s, \psi_m)$
has the same poles and residues in the critical strip as
$\tilde{\lcal}(\psi_m)$.

This completes the proof of Theorem \ref{NONAME}.

\subsection{The analytic case}

We briefly dicuss the extension of Theorem \ref{main2a} to more
general analytic symbols $\sigma \in C^{\omega}(\Gamma \backslash
G). $

By Proposition \ref{EXPRESSION}, by  a similar  calculation as in
Corollary \ref{BASICINT2}, we have
\begin{equation} \label{STEVESZETAa}  \lcal(s; \sigma)  =
\sum_{r_j} \frac{ \langle \sigma, \Xi_{ir_j} \rangle}{\langle
\phi_{ir_j}, \Xi_{ir_j} \rangle} \lcal(s; \phi_{ir_j})
 +
 \sum_{m, \pm} \frac{ \langle \sigma, \Xi_m^{\pm} \rangle}{\langle \psi_m, \Xi_m^{\pm} \rangle}
\lcal(s; \psi_m), \;\; (\Re s > 1).
 \end{equation}
 Here, we interchanged the summation over $\gamma$ and over $r_j$, which is  justified  by  Proposition \ref{XIEST} and the prime geodesic
 theorem. The significant new question is the analytic continuation
 of (\ref{STEVESZETAa}) to the critical strip.
 If $\sigma$ has only finitely many components, then of
 course the analytic continuation of the sums comes down to that of the individual terms.
 To extend the analysis to infinite series, we recall that the
 coefficients \begin{equation} \label{COEFFSFINAL} \frac{ \langle \sigma, \Xi_{ir_j} \rangle}{\langle
\phi_{ir_j}, \Xi_{ir_j} \rangle}, \;\; \frac{ \langle \sigma,
\Xi_m^{\pm} \rangle}{\langle \psi_m, \Xi_m^{\pm} \rangle}
\end{equation}  are rapidly decaying in  $r_j$, resp. $m$  (see
(\ref{ABSCONVLEM})) as long as $\sigma \in C^{\infty}$. Although
we will not investigate it here, for  $\sigma \in C^{\omega}$ we
expect that the coefficients are of exponential decay.

 To extend the analysis to infinite series, it is necessary to
 have estimates of the  zeta functions $\lcal(s; \phi_{ir_j})$,
resp. $\lcal(s; \psi_{m})$  in the parameters $r_j$, resp. $m$. By
Lemma \ref{ZSUMa} and (\ref{LCAL1b}),  the estimates comes down to
the estimates of $\rcal(s; \phi_{ir_j}), \rcal(s; \psi_m)$ in the
$r_k, m$ parameters in the critical strip. By the trace formulae,
this in turn reduces to estimates of the coefficients
\begin{equation} \label{FINALCOEFFS} \frac{\mu_{ir_k}^c(\frac{1}{2} + ir_n)
}{\mu_{ir_k}^c(s)},\;\;\; \frac{\mu_{m}^d(\frac{1}{2} +
ir_n)}{\mu_m^d(s)} \end{equation} in the joint $(r_k, r_n)$ resp.
$(r_k, m)$ parameters. Proposition \ref{MUINTS} reduces these
estimates to ratios of $\Gamma$ functions. Hence, one can obtain
the analytic continuation properties of $\lcal(s; \sigma)$ as long
as the coefficients (\ref{COEFFSFINAL})
 decay rapidly
enough to offset the growth of the $\Gamma$ factors in
(\ref{FINALCOEFFS}).  But  we will not determine the optimal
statement here.


\begin{thebibliography}{HHHH}


\bibitem[A-P]{A-P} A. Alvarez-Parrilla, Explicit geodesic-flow
invariant distributions using $SL(2, \R)$ ladders, Int. J. Math.
and Math. Sci. 2005: 8 (2005), 1299-1315.




\bibitem[AZ]{AZ} N. Anantharaman and S. Zelditch, in progress.


\bibitem[BR]{BR} J.  Bernstein and A.  Reznikov,
Sobolev norms of automorphic functionals. Int. Math. Res. Not.
2002, no. 40, 2155--2174.

\bibitem[C]{C} F. Chamizo,
Automorphic forms and differentiability properties.  Trans. Amer.
Math. Soc. 356 (2004), no. 5, 1909--1935.

\bibitem[Co]{Co} S. Cosentino, A note on Hlder regularity of invariant
distributions for horocycle flows.  Nonlinearity  18  (2005),  no.
6, 2715--2726.


\bibitem[FF]{FF} L. Flaminio and G. Forni,
 Invariant distributions and time averages for horocycle flows. Duke Math. J. 119 (2003), no. 3, 465--526.

\bibitem[GFa]{GFa}  I. M. Gel'fand and S. V. Fomin, Unitary
representations of Lie groups and geodesic flows on surfaces of
constant negative curvature (in Russian), Dokl. Akad. Nauk SSSR 76
(1951), 771--774.

\bibitem[GF]{GF} I.M.  Gel' fand, S. V.  Fomin,
Geodesic flows on manifolds of constant negative curvature. Amer.
Math. Soc. Transl. (2) 1 (1955), 49--65.




\bibitem[GGP]{GGP} I. M.  Gel\' fand, M. I.  Graev and I. I. Pyatetskii-Shapiro, {\it
Representation theory and automorphic functions}.  W. B. Saunders
Co., Philadelphia, Pa.-London-Toronto, Ont. 1969.

\bibitem[GR]{GR} I. S. Gradshteyn and I.M.  Ryzhik, {\it Table of integrals, series, and products.}
   Academic Press, Inc., San Diego, CA, 2000.

\bibitem[G]{G} V. Guillemin, Lectures on spectral theory of elliptic operators. Duke Math. J. 44 (1977), no. 3, 485--517.

\bibitem[H]{H} S. Helgason, {\it Topics in harmonic analysis on homogeneous spaces. }
Progress in Mathematics, 13. Birkh{\"a}user, Boston, Mass., 1981.

\bibitem[He]{He} S. Helgason, {\it Groups and geometric analysis. Integral geometry, invariant differential operators, and spherical functions. }
Corrected reprint of the 1984 original. Mathematical Surveys and
Monographs, 83. American Mathematical Society, Providence, RI,
2000.

\bibitem[J]{J} A. Juhl, {\it
Cohomological theory of dynamical zeta functions}. Progress in
Mathematics, 194. Birkhäuser Verlag, Basel, 2001.

\bibitem[K]{K} A. W. Knapp,{\it Representation theory of semisimple groups. An overview based on examples.}
 Reprint of the 1986 original. Princeton Landmarks in Mathematics. Princeton University Press, Princeton, NJ,
 2001.

\bibitem[L]{L} S.  Lang,  ${\rm SL}\sb 2(R)$.
 Graduate Texts in Mathematics, 105.
Springer-Verlag, New York, 1985.

\bibitem[Lin]{Lin} E. Lindenstrauss, Invariant measures and arithmetic quantum unique
ergodicity (Annals Math., to appear).

\bibitem[Mar1]{Mar1} B. Marcus, Ergodic properties of horocycle flows for surfaces of negative curvature, Ann. of Math. (2) 105 (1977), 81--105.

\bibitem[Mar2]{Mar2}
B. Marcus, The horocycle flow is mixing of all degrees, Invent.
Math. 46 (1978), 201--209.

\bibitem[Mo]{Mo} T. Morita,  Markov systems and transfer operators associated with
cofinite Fuchsian groups. Ergodic Theory Dynam. Systems 17 (1997),
no. 5, 1147--1181.

\bibitem[McK]{McK} McKean, H. P. Selberg's trace formula as applied to a compact Riemann surface.  Comm. Pure Appl. Math.  25  (1972), 225--246.

\bibitem[MS]{MS} S. D. Miller and W.  Schmid,
 The highly oscillatory behavior of automorphic distributions for $\rm SL(2)$. Lett. Math. Phys. 69 (2004), 265--286.


\bibitem[N]{N} P.J. Nicholls, {\it The Ergodic Theory of Discrete Groups}, London Math. Soc.
Lect. Notes Series 143, Cambridge Univ. Press, Cambridge 143.


\bibitem[O]{O} J.P. Otal,
Sur les fonctions propres du laplacien du disque hyperbolique.
 C. R. Acad. Sci. Paris S{\'e}r.
I Math. 327 (1998), no. 2, 161--166.

\bibitem[P]{P} O. S. Parasyuk, Flows of horocycles on surfaces of constant negative curvature (in Russian), Uspekhi Mat. Nauk 8, no. 3 (1953), 125--126.

\bibitem[Pat0]{Pat0} S. J. Patterson,
The limit set of a Fuchsian group. Acta Math. 136 (1976), no. 3-4,
241--273.


\bibitem[Pat1]{Pat1} S. J.
 Patterson,
Lectures on measures on limit sets of Kleinian groups. Analytical
and geometric aspects of hyperbolic space (Coventry/Durham, 1984),
281--323, London Math. Soc. Lecture Note Ser., 111, Cambridge
Univ. Press, Cambridge, 1987

\bibitem[Pol]{Pol} M. Pollicott,  Some applications of thermodynamic formalism to manifolds with constant negative curvature.  Adv. Math.  85  (1991),  no. 2, 161--192.



\bibitem[RS]{RS} Z.  Rudnick and P. Sarnak,
The behaviour of eigenstates of arithmetic hyperbolic manifolds.
Comm. Math. Phys. 161 (1994), no. 1, 195--213.

\bibitem[Ru87]{Ru87}  D. Ruelle, Resonances for Axiom $A$ flows. J. Differential Geom. 25 (1987), no. 1,
99--116.

\bibitem[Rugh92]{Rugh92} H.H.  Rugh,  The correlation spectrum for hyperbolic analytic maps. Nonlinearity 5 (1992), no. 6,
1237--1263.

\bibitem[Rugh96]{Rugh96} H. H.  Rugh,
Generalized Fredholm determinants and Selberg zeta functions for Axiom A dynamical systems. Ergodic Theory Dynam. Systems 16 (1996), no. 4, 805--819.

\bibitem[Sa1]{Sa1} P. Sarnak,  Determinants of Laplacians. Comm. Math. Phys. 110 (1987), no. 1,
113--120.

\bibitem[Sa2]{Sa2} P. Sarnak, {\it
Some applications of modular forms}. Cambridge Tracts in
Mathematics, 99. Cambridge University Press, Cambridge, 1990.

\bibitem[Sh]{Sh} A. I. Schnirelman, Ergodic properties of eigenfunctions,
{\it Usp.\ Mat.\ Nauk.\ } 29/6 (1974), 181--182.

\bibitem[Se]{Se} C. Series,  The infinite word problem and limit sets in Fuchsian groups, ETDS 1 (1981), 337--360.

\bibitem[SV]{SV} L.  Silberman and A. Venkatesh,  On Quantum unique ergodicity for locally symmetric spaces I
( math.RT/0407413).

\bibitem[Su1]{Su1} D. Sullivan, On the ergodic theory at infinity of an arbitrary
discrete group of hyperbolic motions. Riemann surfaces and related topics:
 Proceedings of the 1978 Stony Brook Conference (State Univ. New York, Stony Brook, N.Y., 1978),
 pp. 465--496, Ann. of Math. Stud., 97, Princeton Univ. Press, Princeton, N.J., 1981.

 \bibitem[Su2]{Su2} D. Sullivan,
The density at infinity of a discrete group of hyperbolic motions.
Inst. Hautes Études Sci. Publ. Math. No. 50, (1979), 171--202


\bibitem[V]{V} A. B. Venkov, {\it Spectral theory of automorphic functions and its applications. }
 Mathematics and its Applications (Soviet Series), 51. Kluwer Academic Publishers Group,
 Dordrecht.

 \bibitem[WW]{WW} E.T. Whittaker and G. N.  Watson, {\it A course of modern analysis} Fourth edition. Reprinted Cambridge University Press, New York
1962.



\bibitem[W]{W} S. A. Wolpert, Semiclassical limits for the hyperbolic plane. Duke Math. J. 108 (2001), no. 3, 449--509.

\bibitem[Z]{Z} S. Zelditch,  Trace formula for compact $\Gamma\backslash {\rm
PSL}\sb 2(R)$ and the equidistribution theory of closed geodesics.
Duke Math. J. 59 (1989), no. 1, 27--81.


\bibitem[Z2]{Z2} S. Zelditch, Uniform distribution of eigenfunctions on compact hyperbolic surfaces,
Duke Math. J. 55 (1987), 919--941.

\bibitem[Z3]{Z3} S. Zelditch, Pseudodifferential analysis on hyperbolic surfaces. J. Funct. Anal. 68 (1986), no. 1,
72--105.






\end{thebibliography}
\end{document}